\documentclass[12pt]{article}
\usepackage{amsmath,amssymb,amsthm,amscd}
\let\ftilde\tilde %Plain
\renewcommand{\tilde}[1]{\mathbf{\ftilde{\mathnormal{#1}}}}
\catcode`\@=11
\renewcommand\section{\@startsection{section}{2}{\z@}%
                                     {-3.25ex\@plus -1ex \@minus -.2ex}%
                                    {1.5ex \@plus .2ex}%
                                     {\normalfont\large\bfseries}}
\catcode`\@=\active

\newcommand{\C}{\mathbb{C}}
\newcommand{\Z}{\mathbb{Z}}

\newcommand{\CP}{\mathbb{P}}

\newcommand{\ob}{\operatorname{Ob}}

\newcommand{\vac}{|0\,\rangle}

\newcommand{\unit}{\mathbf{1}}
\newcommand{\rank}{\operatorname{rank}}

\newcommand{\Hom}{\operatorname{Hom}}

\newcommand{\End}{\operatorname{End}}
\newcommand{\Endc}{\mathcal{E}nd}
\newcommand{\Res}[1]{\underset{#1=0}{\operatorname{Res}}\,}

\newcommand{\frakg}{\mathfrak{g}}
\newcommand{\frakh}{\mathfrak{h}}
\newcommand{\frakU}{\mathcal{U}}

\newcommand{\id}{\operatorname{id}}
\newcommand{\Vir}{\mathcal{V}ir}
\newcommand{\mathO}{\mathcal{O}}
\newcommand{\calo}{\mathO}

\newcommand{\calc}{\mathcal{C}}

\newcommand{\calm}{\mathcal{M}}

\newcommand{\calv}{\mathcal{V}}

\newcommand{\cals}{\mathcal{S}}
\newcommand{\calu}{\mathcal{U}}
\newcommand{\calhw}{\mathcal{HW}}
\newcommand{\calch}{\mathcal{CHW}}

\newcommand{\pairing}[2]{\langle#1,#2\rangle}

\newcommand{\ket}[1]{|#1\rangle}
\newcommand{\bra}[1]{\langle#1|}

\newcommand{\gr}{\operatorname{gr}}
\newcommand{\NO}{\,{\raise0.25em\hbox{$\mathop{\hphantom{\cdot}}%
\limits^{_{\circ}}_{^{\circ}}$}}\,}

\newcommand{\II}{I\!I\;}
\newcommand{\III}{I\!I\!I\;}

\catcode`\@=11

\@addtoreset{equation}{section}
\makeatother
\catcode`\@=\active

\theoremstyle{plain}
  
  \newtheorem{theorem}{Theorem}[subsection]%[section]%
  \newtheorem{corollary}[theorem]{Corollary}
  \newtheorem{lemma}[theorem]{Lemma}
  \newtheorem{proposition}[theorem]{Proposition}
\theoremstyle{remark}
  \newtheorem{definition}[theorem]{Definition}
  \newtheorem{remark}[theorem]{Remark}
  \newtheorem{example}[theorem]{Example}
	\newtheorem{note}[theorem]{Note}

\begin{document}
%\begin{flushright}
%Final 1.2\\
%\today
%\end{flushright}
%For Duke Mathematiucal Journal July 22, 2002
\begin{center}
\begin{large}
\textbf{Conformal field theories associated to regular chiral vertex operator algebras I: theories over the projective line}
\end{large}
\vskip 5ex 
\begin{sc}
Kiyokazu Nagatomo %\footnote{e-mail:nagatomo$@$math.sci.osaka-u.ac.jp} 
and Akihiro Tsuchiya
\end{sc}
\end{center}

\begin{small}
\textbf{Abstract:} 
Based on any chiral vertex operator algebra satisfying a suitable finiteness condition, 
the semisimplicity of the zero-mode algebra as well as a regularity for induced modules, 
we construct conformal field theory over the projective line with the chiral vertex operator algebra as symmetries of the theory. 
We appropriately  generalize the argument in \cite{TUY} so that we are able to define sheaves of conformal blocks  for chiral vertex
operator algebras and study them in detail. We prove the factorization theorem under the fairly general conditions for chiral vertex
operator algebras and the zero-mode algebras.
\end{small}

\renewcommand{\thefootnote}{\arabic{footnote}}
\setcounter{footnote}{0}

\tableofcontents

\section*{Introduction}
In the present paper we study conformal field theory (CFT) defined over the projective line associated to chiral vertex operator
algebras. The theories over compact Riemann surfaces of positive genus will be discussed in the forthcoming papers.

There are two extensively studied examples of conformal field theory defined over compact Riemann surfaces; the one
is  WZNW model associated to integrable highest weight representations of affine Lie algebras (\cite{TUY}), and the other is the
minimal models associated to highest weight representations of the co-called minimal series for the Virasoro algebra (\cite{BFM}).
In both examples sheaves of conformal blocks are constructed over the moduli spaces for $N$-pointed stable curves, and the expected 
properties of them such as the coherency, the existence of $\mathcal{D}$-module structure, and the factorization property along
the boundary of the moduli space  are proved as well as other interesting natures. 
There is another well understood example known as abelian conformal field theory whose symmetry is
the lattice chiral vertex operator algebra associated to a rank one even lattice.  
In \cite{U} the space of conformal blocks attached to a pointed stable curve is constructed and its dimension is determined.
Particularly the space of conformal blocks over the smooth curve is isomorphic to the vector space of
theta functions defined over the Jacobian variety with the multiple of the theta divisor.

The notion of chiral vertex operator algebras was first introduced in order to realize and study the Monster
(the largest sporadic finite simple group)(\cite{Borc}, \cite{FLM}), and in fact there are significant achievements in this
direction.  The set of axioms for a chiral vertex operator algebra consists of algebraic counterparts of the notion of the operator
product expansion in conformal field theory. Later conformal field theories associated to chiral vertex operator algebras over
compact Riemann surfaces started in \cite{Z1} and \cite{Z2} as a generalization of WZNW model and the minimal model. 
Meanwhile Beilinson and Drinfeld \cite{BD} formulated a chiral vertex operator algebra over a curve by using the notion of
$\mathcal{D}$-modules,  which is called a \textit{chiral algebra}\,(cf.~\cite{G}). 
Chiral algebras give rise to a framework of conformal field theory over compact Riemann surfaces of
particularly higher genus. In fact Frenkel and Ben-Zvi (\cite{F}, \cite{FB}) defined sheaves of conformal blocks over curves in
terms of chiral algebras. However, their properties such as the coherency, the local freeness along the boundary of the moduli
space, and the factorization property are not known. In the paper we define sheaves of conformal blocks over \textit{the projective
line}\, following \cite{TUY} and prove the series of  above properties. 
Discussions for \textit{higher genus} compact Riemann surfaces will be reported in subsequent papers.

It must be important to note that a chiral vertex operator algebra does not always provide  conformal field theory with desired
properties (say, ``good'' conformal field theory); several conditions such as rationality, regularity, etc., which likely lead us to
the good theory have been proposed. However it is not known,  if one of these conditions in fact allows us to build a theory, 
and which one is the most natural. In this paper we propound such a sufficient condition, which we call the condition \III\,
(see the explanation below); we also present weaker conditions I and \II\!, which are enough for us to prove several properties such
as the coherency of conformal blocks, etc.  In the paper we call a chiral vertex operator algebras satisfying the condition \III 
``regular''; any regular chiral vertex operator algebra gives rise to a good conformal field theory over the projective line. 
One of the important features of the condition \III is that the completely reducibility of $V$-modules is \textit{not}\, assumed:
previously propounded notions of rationality, or regularity presume that any $V$-module is completely reducible. 
We show that our category of $V$-modules is semisimple, and that the category of $V$-modules and the category of finite dimensional
modules for the zero-mode algebra is equivalent.

The most of the ideas in this paper are found in \cite{TUY}, which are slightly generalized so as to fit our general settings.
The notion of conformal blocks in the paper is equivalent to the one in \cite{BD} though its appearance looks different; 
the definition which is convenient for our purpose is chosen.

\vskip 1.5ex

The paper is organized as follows: Section \ref{sec:vertex} reviews basic concepts of chiral vertex operator algebras 
(\cite{FLM}). The definition of a chiral vertex operator algebra consists of several axioms; the most important one is
the associativity formula. Its physics counterpart is known as the operator product
expansion in conformal field theory; the equivalence between them is also explained here.
Section \ref{sec:modules} gives the notions of the current Lie algebra $\frakg(V)$ 
and the current algebra $\calu(V)$ associated to a chiral vertex operator algebra $V$. 
The current algebra $\calu(V)$ is introduced in order to develop the theory of $V$-module in terms
of suitable associative algebras (cf. recall the role of the universal enveloping algebras in the representation theory for Lie
algebras).  For the purpose the Lie algebra $\frakg(V)$ and its completed universal enveloping algebra
$\widehat{U}(\frakg(V))$ are introduced. The algebra $\calu(V)$ is defined as the quotient algebra of $\widehat{U}(\frakg(V))$ by
the two-sided ideal mostly generated by the associativity formula (Definition \ref{def:current-algebra}). 
A module for $V$ or $V$-module is a finitely generated $\calu(V)$-module satisfying a certain finiteness condition;
the category of $V$-modules is denoted by $\mathcal{M}od\,(V)$.
Section \ref{sec:fin} includes consequences under the condition I, i.e., under Zhu's $C_2$-finite condition.
Using the argument given in \cite{GN} being slightly generalized for our situation we prove that any $V$-module is linearly
spanned by vectors obtained by applying negative Fourier modes without same mode iterations (Theorem \ref{th:basis}); this fact is
called the \textit{fermionic property} of $V$-modules. In particular a certain finiteness property called $C_n$-finiteness ($n\geq2$)
is proved for arbitrary $V$-module when $V$ satisfies $C_2$-condition. The coherency of conformal blocks over the projective line is
assured by the $C_2$-finiteness of $V$-modules. In the higher genus case the coherency of conformal blocks can be proved by using
$C_n$-finiteness of $V$-modules for all $n\geq2$. The fermionic property of $V$-modules is a consequence of the associativity
formula, and this itself is a very interesting phenomena which occurs in the theory.  The other important ingredient in the section
is the zero-mode algebra $A_0(V)$ which is defined as a subquotient algebra of the current algebra. The definition is given so as to
have a clear meaning in $V$-modules theory.  In the appendix it is shown that the zero-mode algebra is isomorphic to the so-called
Zhu's algebra, so that the role of Zhu's algebra in the theory of chiral vertex operator algebras is clarified. 
The category of finite dimensional $A_0(V)$-modules is denoted by $\mathcal{M}od\,(A_0(V))$.  The definition of functors
$\mathcal{HW}:\mathcal{M}od\,(V)\rightarrow\mathcal{M}od\,(A_0(V))$ and
$M:\mathcal{M}od\,(A_0(V))\rightarrow\mathcal{M}od\,(V)$ are also given, where $\mathcal{HW}(M)$  is a set of highest weight vectors
of a $V$-module $M$ and $M(W)=\calu(V)\otimes_{F^0\calu(V)}W$ is the induced $V$-module for an $A_0(V)$-module $W$.  
Section \ref{sec:duality} reviews the duality of $V$-modules; the content is essentially due to \cite{FHL}. 
Section \ref{sec:conformalblocks} presents the definitions of
space of conformal blocks and system of current correlation functions over the projective line, and the one to one
correspondence between them (Theorem \ref{th:syscorr}) is established. One of the main theorems that any space of conformal blocks is
finite dimensional is proved (Theorem \ref{th:finiteblock}). 
In section \ref{sec:connection} the sheaf of conformal blocks is defined over the moduli space of $N$-pointed projective line,
and an integrable connection on it is constructed. 
Section \ref{sec:regularity} states the conditions I\!I ($C_2$-finite condition and the semisimplicity of the zero-mode algebra) and
I\!I\!I\,(the condition \II and the simpleness of the induced module $M(W)$ for a simple $A_0(V)$-module $W$). 
The most important contributions of the present work is in section \ref{sec:factorization}.  In this section we always work with
the condition \III\!. Sheaves of conformal blocks are extended to sheaves over the partial compactification of the moduli space of
$N$-pointed projective line (Proposition \ref{prop:extension}), whose boundary consists of $N$-pointed stable curves with ordinary
double points. It is verified that the connection for vector fields tangent to the boundary is well-defined (Theorem
\ref{th:extension-connection}).  The behavior  of sheaves of conformal blocks along the boundary of the moduli space is studied 
under the conditions \III\!.  It is shown that sheaves of conformal blocks over the moduli space are locally free (Theorem
\ref{th:locally-free}) and have the factorization property along the boundary (Theorem \ref{th:dstr}). 
Section \ref{sec:category} shows that if $V$ is a chiral vertex operator algebra satisfying the conditions \III\!
then the functor $\mathcal{HW}:\mathcal{M}od\,(V)\rightarrow\mathcal{M}od\,(A_0(V))$ gives an equivalence of categories.
Particularly $\mathcal{M}od\,(V)$ is a semisimple category.
In the appendix we include two example of regular chiral vertex operator algebras; the affine and the Virasoro chiral vertex
operator algebras, which gives WZNW-model and the minimal model, respectively.
These two example would give more concrete images of notions discussed at every stage of the paper.
The isomorphism between the zero-mode algebra and Zhu's algebra is also given here. 

\section{Chiral vertex operator algebras}
\label{sec:vertex}
We give definition of chiral vertex operator algebras and characterize the axioms for chiral vertex operator algebras in terms of
2-point functions.

\subsection{Graded vector spaces and endomorphism rings}

Let $M=\bigoplus_{n\in\Z}M_n$ be a graded vector space over the complex number field $\C$ such that $\dim\,M_n<\infty$ 
for all $n\in\Z$ and $M_n=0$ for all sufficiently small $n$. We define the filtration $F_pM\,(p\in\Z)$ by setting
$F_pM=\bigoplus_{n\leq p}M_n$,  which is an increasing filtration on $M$;
\[
\dots\subseteq F_pM\subseteq F_{p+1}M\subseteq\dots,\quad
\cap_{p\in\Z}F_pM=\{0\}\mbox{ and } M=\cup_{p\in\Z}F_pM.
\]
Setting $F^pM=F_{-p}M$ for all $p\in\Z$ we get the decreasing filtration
\[
\dots\supseteq F^pM\supseteq F^{p+1}M\supseteq\dots,\quad
\cap_{p\in\Z}F^pM=\{0\}\mbox{ and } M=\cup_{p\in\Z}F^pM.
\]

Let $\End\,M$ be the endomorphism ring of $M$. 
We set $\Endc\,M=\cup_{p\in\Z}F_p\Endc\,M$ where
\[
F_p\Endc\,M=\{\;\varphi\in\End\,M\;|\;\varphi(M_n)\subseteq F_{n+p}M\mbox{ for all $n\in\Z$}\;\}\quad(p\in\Z).
\]
We also define $F^p\Endc\,M=F_{-p}\Endc\,M$ for all $p\in\Z$. 
The $F_p\Endc\,M\,(p\in\Z)$ define an increasing filtration on $\Endc\,M$
\[
\ldots \subseteq
F_p\Endc\,M\subseteq F_{p+1}\Endc\,M\subseteq\ldots,\quad\cap_{p\in\Z}F_p\Endc\,M=\{0\},
\]
and the $F^p\Endc\,M\,(p\in\Z)$ define a decreasing filtration.
The vector space $\Endc\,M$ is Hausdorff with respect to the linear topology
induced by the filtration $F_p\Endc\,M\,(p\in\Z)$ and is complete. 

We set $\Endc^{\,f}\,M=\bigoplus_{p\in\Z}\Endc_p\,M\subset\Endc\,M$ where
\[
\Endc_p\,M=\{\;\varphi\in\End\,M\;|\;\varphi(M_n)\subseteq M_{n+p}\mbox{ for all $n\in\Z$}\;\}\quad(p\in\Z).
\] 
The vector space $\Endc^{\,f}\,M$ is endowed with the increasing filtration 
$F_p\Endc^{\,f}\,M=F_p\Endc\,M\cap\Endc^{\,f}\,M=\bigoplus_{n\leq p}\Endc_n\,M$. We see that $\Endc\,M$ is the formal completion of
$\Endc^{\,f}\,M$.

\subsection{The Virasoro algebra}
The Virasoro algebra $\Vir$ is a vector
space
$\Vir=\bigoplus_{n\in\Z}\C T(n)\oplus\C C$ with commutation relations
\[
[T(m),T(n)]=(m-n)T(m+n)+\frac{m^3-m}{12}\delta_{m+n,0}\,C,\quad [T(m),C]=0.
\]
Let $M$ be a $\Vir$-module such that the center $C$ acts on $M$ as a complex number $c_M$. 
Then $c_M$ is called \textit{central charge}.

\subsection{Chiral vertex operator algebras - definition}
Let us recall definition of chiral vertex operator algebras (\cite{Borc}, \cite{FLM}).

\begin{definition}\label{def:chiral}
The triplet $(V,\,J,\,T\,)$ is called a \textit{chiral vertex operator algebra} if it satisfies the following conditions:

\vskip 1.5ex
\noindent
(a) $V$ is a graded vector space $V=\bigoplus_{\Delta=0}^\infty V_\Delta$ such that $\dim V_0=1$ and $\dim V_\Delta<\infty$
for all $\Delta$. Any element $v$ of $V_\Delta$ is called \textit{homogeneous} vectors with weight $\Delta$;
we denote $\Delta=|v|$.

\vskip 1ex
\noindent
(b) There exist two special nonzero elements $\vac\in V_0$ and $T\in V_2$, 
which are called the \textit{vacuum} and the \textit{Virasoro vector} respectively.

\vskip 1ex
\noindent
(c) For any integer $n$ there exists a linear map $J_n:V\rightarrow \Endc_{-n}\,V$ such that $J_0(\vac)= \id_V,\,J_n(\vac)=0\,(n\neq 0)$,
and $J_{-\Delta}(v)\ket{0}=v,\, J_n(v)\ket{0}=0\,(n>-\Delta)$ for all $v\in V_\Delta$.
We set $T(n)=J_n(T)$ for all $n\in\Z$. The $T(n)\,(n\in\Z)$ and $J_0(\vac)=\id_V$ form a representation of the Virasoro algebra 
on $V$ with some central charge $c_V\in\C$, and $T(0)$ is the grading operator, i.e., $T(0)v=\Delta v$ for any $v\in V_\Delta$.

\vskip 1ex
\noindent
(d) $J_n(T(-1)v)=-(n+\Delta)J_{n}(v)$ for any $v\in V_\Delta$.

\vskip 1ex
\noindent
(e) Let $v_1\in V_{\Delta_1},\,v_2\in V_{\Delta_2}$ and $m,\,n\in\Z$. Then
\[
[J_m(v_1),J_n(v_2)]=\sum_{j=0}^{\Delta_1+\Delta_2-1}\binom{m+\Delta_1-1}{j}J_{m+n}(J_{j-\Delta_1+1}(v_1)v_2).
\]

\vskip 1ex
\noindent
(f) Let $v_1\in V_{\Delta_1},\,v_2\in V_{\Delta_2}$ and $m,\,n\in\Z$. Then
\begin{multline*}
J_m(J_n(v_1)v_2)=
\sum_{j=0}^\infty(-1)^j\binom{n+\Delta_1-1}{j}\\
\times
\left(
J_{n-j}(v_1)J_{m-n+j}(v_2)
-(-1)^{n+\Delta_1-1}J_{m+\Delta_1-j-1}(v_2)J_{j-\Delta_1+1}(v_1)
\right).
\end{multline*}
\end{definition}

\begin{note}
The identities (e) and (f) is called the \textit{commutator formula} and the \textit{associativity formula}
respectively in the literature.
\end{note}

The following proposition is well-known (\cite{Borc}, \cite{FLM}).

\begin{proposition}\label{prop:for}
{\rm (1)} $T(1)T=0$ and $\displaystyle{T(2)T=\frac{c_V}{2}\vac}$.

\vskip 1ex
\noindent
{\rm(2)} For any $v\in V_\Delta$ and $n\in\Z$,
\[
[T(-1),J_n(v)]=-(n+\Delta-1)J_{n-1}(v),\quad
[T(0),J_n(v)]=-nJ_{n}(v).
\]

\vskip 1ex
\noindent
{\rm(3)} For any $v_1\in V_{\Delta_1},\,v_2\in V_{\Delta_2}$ and $n\in\Z$,
\[
J_n(v_1)v_2=\sum_{j=0}^\infty\frac{(-1)^{n+\Delta_1+j}}{j!}T(-1)^j
\left(
J_{n+\Delta_1-\Delta_2+j}(v_2)v_1
\right).
\]
\end{proposition}

\begin{note}
The first formula in the proposition (2) is called the \textit{derivation property}, and the one in the proposition (3) is
called the
\textit{skew symmetry}.
\end{note}

\subsection{$2$-point functions}
\label{sub:2point}
In this subsection we present another set of axioms for a chiral vertex operator algebra using the notion of the operator product
expansion.

\begin{definition}
Let $V=\bigoplus_{\Delta=0}^\infty V_\Delta$ be a chiral vertex operator algebra. For any $v\in V_\Delta$ we associate a formal power
series with the formal variable $z$ and $z^{-1}$ with coefficients in $\Endc^{\,f}\,V$ by
\[
J(v,z)=\sum_{n\in\Z}J_n(v)z^{-n-\Delta}\in
(\Endc^{\,f}\,V)[[z,z^{-1}]],
\]
which is called the \textit{nonabelian current} associated to the vector $v$. For any $v_1\in V_{\Delta_1}$ and $v_2\in V_{\Delta_2}$
the product $J(v_1,z_1)J(v_2,z_2)$ is well-defined as an element in $(\Endc^{\,f}\,V)[[z_1,z_1^{-1},z_2,z_2^{-1}]]$.
\end{definition}

\begin{definition}
Let $V=\bigoplus_{\Delta=0}^\infty V_\Delta$ be a chiral vertex operator algebra. We denote by $V^*$ the graded dual of $V$, i.e., 
$V^*=\bigoplus_{\Delta=0}^\infty V_\Delta^*$ where
$V_\Delta^*=\Hom_\C(V_\Delta,\C)$, and by $\pairing{\;\,}{\;}$ the natural dual pairing between $V^*$ and $V$.
\end{definition}

\begin{proposition}
Let  $v_1\in V_{\Delta_1}$ and $v_2\in V_{\Delta_2}$.

\vskip 1.5ex
\noindent
{\rm (1)} The formal $2$-point function $\bra{v^*}J(v_1,z_1)J(v_2,z_2)\ket{v}\in\C[[z_1,z_1^{-1},z_2,z_2^{-1}]]$
absolutely converges on the domain $|z_1|>|z_2|>0$ and is analytically continued to a rational function on $\CP^1\times\CP^1$
with possible poles on $z_1=z_2$, $z_i=0\,(i=1,2)$ and $z_i=\infty\,(i=1,2)$; this rational function is denoted by the same notation.  

\vskip 1ex
\noindent
{\rm (2)} The formal power series 
$\bra{v^*}J(J(v_1,z_1-z_2)v_2,z_2)\ket{v}\in\C[[z_1-z_2,(z_1-z_2)^{-1},z_2,z_2^{-1}]]$
absolutely converges on the domain $0<|z_1-z_2|<|z_2|$ and is analytically continued to a rational function on $\CP^1\times\CP^1$
with possible poles on $z_1=z_2$, $z_i=0\,(i=1,2)$ and $z_i=\infty\,(i=1,2)$;
this rational function is also denoted  by the same notation.

\vskip 1ex
\noindent
{\rm (3)} The following relations hold
\begin{align*}
\bra{v^*}J(v_1,z_1)J(v_2,z_2)\ket{v}=&\bra{v^*}J(v_2,z_2)J(v_1,z_1)\ket{v}\\
=&\bra{v^*}J(J(v_1,z_1-z_2)v_2,z_2)\ket{v}\\
=&\bra{v^*}J(J(v_2,z_2-z_1)v_1,z_1)\ket{v}.
\end{align*}

\vskip 1ex
\noindent
{\rm (4)} For any $v_1\in V_{\Delta_1}$
\[
\frac{d}{dz}\bra{v^*}J(v_1,z)\ket{v}=\bra{v^*}J(T(-1)v_1,z)\ket{v}.
\]
\end{proposition}
\begin{proof}
For $v\in V_\Delta$ we set 
\[
J^{>0}(v,z)=\sum_{n\geq-\Delta+1}J_n(v)z^{-n-\Delta}\quad\mbox{and}\quad J^{\leq0}(v,z)=\sum_{n\leq-\Delta}J_n(v)z^{-n-\Delta}.
\]
Note that 
\begin{align}
&\bra{v^*}J(v_1,z_1)J(v_2,z_2)\ket{v}=\bra{v^*}[J^{>0}(v_1,z_1),J(v_2,z_2)]\ket{v}\notag\\
&\phantom{\bra{v^*}J(v_1,z_1)J(v_2,z_2)}
+\bra{v^*}J(v_2,z_2)J^{>0}(v_1,z_1)\ket{v}+\bra{v^*}J^{\leq0}(v_1,z_1)J(v_2,z_2)\ket{v},\label{eqn:2fuunc1}\\
&\bra{v^*}J(v_2,z_2)J(v_1,z_1)\ket{v}=\bra{v^*}[J(v_2,z_2),J^{\leq0}(v_1,z_1),]\ket{v}\notag\\
&\phantom{\bra{v^*}J(v_1,z_1)J(v_2,z_2)}+\bra{v^*}J(v_2,z_2)J^{>0}(v_1,z_1)\ket{v}+\bra{v^*}J^{\leq0}(v_1,z_1)J(v_2,z_2)\ket{v},
\label{eqn:2fuunc2}
\end{align}
where the second and the third term in the right hand sides are Laurent polynomials of $z_1$ and $z_2$.
The commutator formula in Definition \ref{def:chiral} is equivalent to the following set of equations, each of which holds as an
element in 
$(\Endc^{\,f}\,V)[[z_1^{-1},z_2,z_2^{-1}]]$ and $(\Endc^{\,f}\,V)[[z_1,z_2,z_2^{-1}]]$, respectively;
\begin{align}
&[J^{>0}(v_1,z_1),J(v_2,z_2)]=\sum_{j=0}^{\Delta_1+\Delta_2-1}J(J_{j-\Delta_1+1}(v_1)v_2,z_2)
\left.(z_1-z_2)^{-j-1}\right|_{|z_1|>|z_2|},\label{eqn:com101}\\
&[J(v_2,z_2),J^{\leq0}(v_1,z_1)]=\sum_{j=0}^{\Delta_1+\Delta_2-1}J(J_{j-\Delta_1+1}(v_1)v_2,z_2)
\left.(z_1-z_2)^{-j-1}\right|_{|z_1|<|z_2|},\label{eqn:com102}
\end{align}
where
\begin{align*}
&\left.(z_1-z_2)^{-j-1}\right|_{|z_1|>|z_2|}=\sum_{m=0}^\infty\binom{m}{j}z_1^{-m-1}z_2^{m-j},\\
&\left.(z_1-z_2)^{-j-1}\right|_{|z_1|<|z_2|}=-\sum_{m=0}^\infty\binom{-m-1}{j}z_1^{m}z_2^{-m-1-j}.
\end{align*}

Using (\ref{eqn:com101}) and (\ref{eqn:com102}) we see that the expansion of
$J(J^{>0}(v_1,z_1-z_2)v_2,z_2)\in(\Endc^{\,f}\,V)[[z_1-z_2,(z_1-z_2)^{-1},z_2,z_2^{-1}]]$ on the domain 
$|z_1|>|z_2|>0$ and $|z_2|>|z_1|>0$ is respectively given by 
\begin{equation}\label{eqn:ope201}
J(J^{>0}(v_1,z_1-z_2)v_2,z_2)=
\begin{cases}
[J^{>0}(v_1,z_1),J(v_2,z_2)]&\quad\mbox{on $|z_1|>|z_2|>0$},\\
{}[J(v_2,z_2),J^{\leq0}(v_1,z_1)]&\quad\mbox{on $|z_2|>|z_1|>0$}.
\end{cases}
\end{equation}
Finally using the associativity formula we see that 
\begin{equation}\label{eqn:asso01}
J(J^{\leq0}(v_1,z_1-z_2)v_2,z_2)
=J^{\leq0}(v_1,z_1)J(v_2,z_2)+J(v_2,z_2) J^{>0}(v_1,z_1)
\end{equation}
as elements in $(\Endc^{\,f}\,V)[[z_1-z_2,z_2,z_2^{-1}]]$, where in the right hand side we use the expansion 
\[
z_1^{-j-1}=\sum_{m=0}^\infty\binom{-j-1}{m}(z_1-z_2)^mz_2^{-m-j-1}.
\] 

Using  (\ref{eqn:com101}) we see that the first term of (\ref{eqn:2fuunc1}) absolutely converges on the domain $|z_1|>|z_2|>0$ to a
rational function with the prescribed pole conditions, and that 
$\bra{v^*}J(v_1,z_1)J(v_2,z_2)\ket{v}=\bra{v^*}J(v_2,z_2)J(v_1,z_1)\ket{v}$ by (\ref{eqn:2fuunc1})-(\ref{eqn:com102}). The formal
power series
$\bra{v^*}J(J(v_1,z_1-z_2)v_2,z_2)\ket{v}$ also absolutely converges on the domain $0<|z_1-z_2|<|z_2|$ to
$\bra{v^*}J(v_1,z_1)J(v_2,z_2)\ket{v}$ by (\ref{eqn:ope201}) and (\ref{eqn:asso01}).

Finally (4) follows from the condition (d) for the chiral vertex operator algebra.
\end{proof}

\begin{note}
The first and second equality of the proposition (3) is called the \textit{$\mathcal{S}_2$-symmetry} and
the \textit{operator product property}, respectively. 
\end{note}

\begin{remark}
The $\cals_2$-symmetry and the operator product property give the commutator formula and the associativity formula in Definition
\ref{def:chiral}. The third equality of the proposition (3) gives the skew-symmetry. 
The statement (4) is nothing but (d) in Definition \ref{def:chiral}
\end{remark}

\section{Modules for chiral vertex operator algebras}\label{sec:modules}
\subsection{Current Lie algebras}\label{current}
Let $(V,\,J,\,T\,)$  be a chiral vertex operator algebra with a grading $V=\bigoplus_{\Delta=0}^\infty V_\Delta$. 
We set $\widetilde{V}^{(1)}=\bigoplus_{\Delta=0}^\infty V_\Delta\otimes\C((\xi))(d\xi)^{1-\Delta}$ and
$\widetilde{V}^{(0)}=\bigoplus_{\Delta=0}^\infty V_\Delta\otimes\C((\xi))(d\xi)^{-\Delta}$. 
Let $\nabla:\widetilde{V}^{(0)}\rightarrow\widetilde{V}^{(1)}$ be the linear map defined by
\[
v\otimes f(\xi)(d\xi)^{-\Delta}\longmapsto T(-1)v\otimes f(\xi)(d\xi)^{-\Delta}+v\otimes \frac{df(\xi)}{d\xi}(d\xi)^{1-\Delta}.
\]
We set $\frakg(V)=\widetilde{V}^{(1)}/\nabla\widetilde{V}^{(0)}$ and denote by $J(v,f)$ the image of $v\otimes f(d\xi)^{1-\Delta}\,(v\in V_\Delta)$
under the canonical  projection  $\widetilde{V}^{(1)}\rightarrow\frakg(V)=\widetilde{V}^{(1)}/\nabla\widetilde{V}^{(0)}$.

The vector space $\widetilde{V}^{(1)}$ is filtered by the decreasing filtration
\begin{align*}
&F^p\widetilde{V}^{(1)}=\bigoplus_{\Delta=0}^\infty V_\Delta\otimes\C[[\xi]]\xi^{p+\Delta-1}(d\xi)^{1-\Delta}\quad(p\in\Z),\\
&\ldots\supseteq F^p\widetilde{V}^{(1)}\supseteq F^{p+1}\widetilde{V}^{(1)}\supseteq\ldots,
\quad\underset{p}{\cap}F^p\widetilde{V}^{(1)}=0
\quad\mbox{and}\quad\widetilde{V}^{(1)}=\underset{p}{\cup}F^p\widetilde{V}^{(1)},
\end{align*}
which defines a Hausdorff linear topology on the vector space $\widetilde{V}^{(1)}$.
There also exists a decreasing filtration on $\widetilde{V}^{(0)}$;
\[
F^p\widetilde{V}^{(0)}=\bigoplus_{\Delta=0}^\infty V_\Delta\otimes\C[[\xi]]\xi^{p+\Delta-1}(d\xi)^{-\Delta}.
\]
Since $\nabla:F^p\widetilde{V}^{(0)}\rightarrow F^{p-1}\widetilde{V}^{(1)}$ for all $p\in\Z$
the filtration $F^p\,(p\in\Z)$ on $\tilde{V}^{(1)}$ induces a filtration on the Lie algebra $\frakg(V)$, which we denote by
$F^p\frakg(V)\,(p\in\Z)$.

We introduce a bilinear operation 
$[\;\,,\;]:\widetilde{V}^{(1)}\times\widetilde{V}^{(1)}\rightarrow\widetilde{V}^{(1)}$ by
\begin{multline*}
[v_1\otimes f_1(d\xi)^{1-\Delta_1}, v_2\otimes f_2(d\xi)^{1-\Delta_2}]\\
=\sum_{m=0}^{\Delta_1+\Delta_2-1}\frac{1}{m!}
J_{m-\Delta_1+1}(v_1)v_2\otimes\frac{d^mf_1}{d\xi^m}f_2\,(d\xi)^{m+2-\Delta_1-\Delta_2}.
\end{multline*}
\begin{proposition}\label{prop:lie}
The bilinear operation $[\;\,,\;]$ induces a Lie bracket on $\frakg(V)$, and the Lie algebra $\frakg(V)$ is filtered by
$F^p\frakg(V)\,(p\in\Z)$, i.e., $[F^p\frakg(V),F^q\frakg(V)]\subseteq F^{p+q}\frakg(V)$.
\end{proposition}
\begin{proof}
In order to show that the bilinear operation $[\;,\,]$ on $\tilde{V}^{(1)}$ induces a bilinear map on $\frakg(V)$ we first prove that
\[
[v_1\otimes f_1, v_2\otimes f_2]+[v_2\otimes f_2, v_1\otimes f_1]\equiv0\mod\nabla\widetilde{V}^{(0)}.
\]
for all $v_i\in V_{\Delta_i}$ and $f_i\in\C((\xi))(d\xi)^{1-\Delta_i}\,(i=1,2)$.
Hereafter we omit the symbol $d\xi$ for convenience.
To do this it suffices to see that
\[
[v_1\otimes f_1, v_2\otimes f_2]=\sum_{j=0}^\infty\sum_{n=0}^\infty\frac{(-1)^{j+1}}{j!n!}\nabla^j
\left(J_{n+j-\Delta_2+1}(v_2)v_1\otimes\frac{d^nf_2}{d\xi^n}f_1\right)=:D_R.
\]
The key is the following identity for differential polynomials;
\[
f_2\frac{d^if_1}{d\xi^i}=\sum_{k=0}^i
(-1)^{i+k}\binom{i}{k}\frac{d^k}{d\xi^k}\left(f_1\frac{d^{i-k}f_2}{d\xi^{i-k}}\right).
\]
Using the Leibnitz rule for $\nabla=T(-1)\otimes \id+\id\otimes d/d\xi$ we see that
\[
D_R=\sum_{j=0}^\infty\sum_{n=0}^\infty\frac{(-1)^{j+1}}{j!n!}\sum_{k=0}^j
\binom{j}{k}T(-1)^{j-k}J_{n+j-\Delta_2+1}(v_2)v_1\otimes\frac{d^k}{d\xi^k}\left(\frac{d^nf_2}{d\xi^n}f_1\right).
\]
We now set $m=j-k,\,i=n+k$, and eliminate $j$ and $n$;
\begin{align*}
D_R
&=\sum_{i,\,m=0}^\infty\frac{(-1)^{i+m+1}}{i!m!}
T(-1)^{m}J_{m+i-\Delta_2+1}(v_2)v_1\otimes\sum_{k=0}^i(-1)^{i+k}\binom{i}{k}\frac{d^k}{d\xi^k}\left(\frac{d^{i-k}f_2}{d\xi^{i-k}}f_1\right),\\
&=\sum_{i,\,m=0}^\infty\frac{(-1)^{i+m+1}}{i!m!}
T(-1)^{m}J_{m+i-\Delta_2+1}(v_2)v_1\otimes\frac{d^{i}f_1}{d\xi^{i}}f_2\,,\\
&=[v_1\otimes f_1, v_2\otimes f_2]
\end{align*}
by the skew symmetry.

Using the property (d) of chiral vertex operator algebras we obtain $[\nabla\widetilde{V}^{(0)},\tilde{V}^{(1)}]=0$, and
then $[\tilde{V}^{(1)},\nabla\widetilde{V}^{(0)}]\subset\nabla\widetilde{V}^{(0)}$ by the $\operatorname{mod}\;\nabla\widetilde{V}^{(0)}$-skew
symmetry. Now the bilinear operation $[\;\,,\;]$ on $\widetilde{V}^{(1)}$ induces a bilinear map on $\frakg(V)$ such that 
$[F^p\frakg(V),F^q\frakg(V)]\subseteq F^{p+q}\frakg(V)$.

To see that this bilinear map is a Lie bracket it suffices to show  that the $[\;\,,\;]$ satisfies the Jacobi identity, i.e.,
\begin{multline*}
\left[[v_1\otimes f_1,v_2\otimes f_2],v_3\otimes f_3\right]\\
=\left[v_1\otimes f_1,[v_2\otimes f_2,v_3\otimes f_3]\right]
-\left[v_2\otimes f_2,[v_1\otimes f_1,v_3\otimes f_3]\right]
\end{multline*}
for all $v_i\in V_{\Delta_i}$ and $f_i\in\C((\xi))(d\xi)^{1-\Delta_i}\,(i=1,2,3)$.
Let us denote the right hand side of the Jacobi identity by $J_R$. 
By definition of the bracket we see that
\[
J_R=\sum_{m,\,n=0}^\infty\frac{1}{m!\,n!}
[J_{m-\Delta_1+1}(v_1),J_{n-\Delta_2+1}(v_2)]v_3\otimes\frac{d^mf_1}{d\xi^m}\frac{d^nf_2}{d\xi^n}f_3.
\]
We now apply the commutator formula to the right hand side so that
\[
J_R
=\sum_{m,\,n=0}^\infty\sum_{j=0}^m\frac{1}{m!\,n!}\binom{m}{j}
J_{m+n-\Delta_1-\Delta_2+2}\left(J_{j-\Delta_1+1}(v_1)v_2\right)v_3\otimes\frac{d^mf_1}{d\xi^m}\frac{d^nf_2}{d\xi^n}f_3.
\]
Set $i=m+n-j$ and eliminate $m$;
\begin{align*}
J_R&=\sum_{i,\,j=0}^\infty\frac{1}{i!\,j!}
J_{i+j-\Delta_1-\Delta_2+2}\left(J_{j-\Delta_1+1}(v_1)v_2\right)v_3
\otimes\sum_{n=0}^i\binom{i}{n}\frac{d^{j+i-n}f_1}{d\xi^{j+i-n}}\frac{d^nf_2}{d\xi^n}f_3\\
&=\sum_{i,\,j=0}^\infty\frac{1}{i!\,j!}
J_{i+j-\Delta_1-\Delta_2+2}\left(J_{j-\Delta_1+1}(v_1)v_2\right)v_3
\otimes\frac{d^{i}}{d\xi^{i}}\left(\frac{d^jf_1}{d\xi^j}f_2\right)f_3\\
&=\left[\sum_{j=0}^\infty\frac{1}{j!}
J_{j-\Delta_1+1}(v_1)v_2\otimes\frac{d^jf_1}{d\xi^j}f_2, v_3\otimes f_3\right]\\
&=\left[[v_1\otimes f_1,v_2\otimes f_2],v_3\otimes f_3\right].
\end{align*}
\end{proof}

\begin{definition}
Let $V$ be a chiral vertex operator algebra. The Lie algebra $\frakg(V)$ with the Lie algebra structure 
\[
[J(v_1,f_1), J(v_2,f_2)]
=\sum_{m=0}^{\Delta_1+\Delta_2-1}\frac{1}{m!}
J(J_{m-\Delta_1+1}(v_1)v_2,\frac{d^mf_1}{d\xi^m}f_2)
\]
is called the \textit{current Lie algebra}\, associated to the chiral vertex operator algebra $V$.

For any $v\in V_\Delta$ and $n\in\Z$ we set $J_n(v)=J(v,\xi^{n+\Delta-1})$.
Specifically for $v=T$ we denote $T(n)=J(T,\xi^{n+1})$.
\end{definition}

We could also do the whole story for $\bigoplus_{\Delta=0}^\infty V_\Delta\otimes\C[\xi,\xi^{-1}](d\xi)^{1-\Delta}$;
then we get a graded Lie algebra $\frakg(V)^f=\bigoplus_{n\in\Z}\frakg(V)^f_n$ where $\frakg(V)^f_n$ is linearly spanned by 
$J_n(v)\,(v\in V)$. We see that $\frakg(V)^f$ is a Lie subalgebra of $\frakg(V)$ and that $J_n(T(-1)v)=-(n+\Delta)J_n(v)$ 
for all $v\in V_\Delta$ and $n\in\Z$.

\begin{remark}
{\rm (1)}
$J_n(\vac)=0$ for all $n\neq0$.

\vskip 1ex
\noindent
{\rm (2)} Let $v_1\in V_{\Delta_1}$ and $v_2\in V_{\Delta_2}$. Then
\[
[J_m(v_1),J_n(v_2)]=\sum_{j=0}^{\Delta_1+\Delta_2-1}\binom{m+\Delta_1-1}{j}J_{m+n}(J_{j-\Delta_1+1}(v_1)v_2),
\]
and in particular the linear map $\frakg(V)^f\rightarrow\Endc^{\,f}\,V$ defined by $J_n(v)=J(v,\xi^{n+\Delta-1})\mapsto J_n(v)$
for all $v\in V_\Delta$ and $n\in\Z$ is a Lie algebra homomorphism. 

\vskip 1ex
\noindent
{\rm (3)} The vector space $\Vir(V)=\bigoplus_{n\in\Z}\C T(n)\oplus\C J_0(\,\ket{0})$ is a graded Lie subalgebra of $\frakg(V)^f$,
and has the bracket relations
\[
[T(m),T(n)]=(m-n)T(m+n)+\frac{m^3-m}{12}\delta_{m+n,0}c_VJ_0(\,\ket{0}),\quad
[T(n),J_0(\,\ket{0})]=0.
\]
\end{remark}

\subsection{Current algebras}
We set $\frakg(V)_{\geq0}=F^0\frakg(V)(=F_0\frakg(V))$ and $\frakg(V)_{<0}=\bigoplus_{n<0}\frakg(V)^{\,f}_n$.
The vector spaces $\frakg(V)_{\geq0}$ and $\frakg(V)_{<0}$ are Lie subalgebras of $\frakg(V)$
and $\frakg(V)=\frakg(V)_{<0}\oplus\frakg(V)_{\geq0}$. 
Since $\frakg(V)_{<0}$ is a Lie subalgebra of $\frakg(V)^f$, the notation $\frakg(V)^f_{<0}$ is also used.
Let $U(\frakg(V))$ be the universal enveloping algebra of the Lie algebra $\frakg(V)$;
there is a canonical vector space isomorphism
$U(\frakg(V))\cong U(\frakg(V)_{<0})\otimes U(\frakg(V)_{\geq0})$.  
For any $p\in\Z$ we define a vector subspace of $U(\frakg(V))$ by
\[
F^p\,U(\frakg(V))=\sum_{p_1+\cdots+p_r\geq p}F^{p_1}\frakg(V)\cdots F^{p_r}\frakg(V),
\]
where we suppose that $\C\cdot1\subseteq F^p\,U(\frakg(V))$ for all $p\leq0$, which defines a decreasing filtration.
Then $U(\frakg(V))$ becomes a filtered algebra and a filtered Lie algebra with the canonical Lie
algebra structure; for all $p,\,q\in\Z$
\begin{align*}
&F^p\,U(\frakg(V))\cdot F^q\,U(\frakg(V))\subseteq  F^{p+q}\,U(\frakg(V)),\\
&[ F^p\,U(\frakg(V)), F^q\,U(\frakg(V))]\subseteq  F^{p+q}\,U(\frakg(V)).
\end{align*}
The filtration $F^p\,(p\in\Z)$ on $U(\frakg(V))$ induces a decreasing filtration on $U(\frakg(V)_{<0})$ and $U(\frakg(V)_{\geq0})$
such that $F^1\,U(\frakg(V)_{<0})=0$ and $F^0\,U(\frakg(V)_{\geq0})=U(\frakg(V)_{\geq0})$, respectively. 
We note that for every $p\in\Z$
\[
F^p\,U(\frakg(V))=\sum_{\substack{p_1+p_2=p\\p_1\leq0,\,p_2\geq0}}F^{p_1}\,U(\frakg(V)_{<0})\cdot F^{p_2}\,U(\frakg(V)_{\geq0}).
\]

For any fixed $p\in\Z$ we introduce a decreasing filtration on  $F^p\,U(\frakg(V))$ by setting
\[
F^p_N\,U(\frakg(V))=\sum_{\substack{p_1+p_2=p\\p_1\leq0,\,p_2\geq N}}F^{p_1}\,U(\frakg(V)_{<0})\cdot F^{p_2}\,U(\frakg(V)_{\geq0})
\quad\mbox{for}\quad N\in\Z_{\geq0}.
\]
Then by definition
$F^{p_1}_{N_1}\,U(\frakg(V))\cdot F^{p_2}_{N_2}\,U(\frakg(V))\subseteq F^{p_1+p_2}_{N_2}\,U(\frakg(V))$
for all $p_1,\,p_2\in\Z$ and $N_1,\,N_2\in\Z_{\geq0}$, and if $N_1+p_2\geq0$ then 
$F^{p_1}_{N_1}\,U(\frakg(V))\cdot F^{p_2}_{N_2}\,U(\frakg(V))\subseteq F^{p_1+p_2}_{N_1+p_2}\,U(\frakg(V))$
by the commutator formula.
We now set
\[
F^p\,\widehat{U}(\frakg(V))=\underset{\underset{N}{\longleftarrow}}{\lim}\;F^p\,U(\frakg(V))/F^p_N\,U(\frakg(V)),
\]
i.e., $F^p\,\widehat{U}(\frakg(V))$ is the formal completion of $F^p\,U(\frakg(V))$ with respect to the linear topology
induced by the filtration $F^p_N\,(N\in\Z_{\geq0})$.
Since $F^{p+1}_N\,U(\frakg(V))=F^{p+1}\,U(\frakg(V))\cap F^{p}_N\,U(\frakg(V))$ for all $p\in\Z$ and $N\in\Z_{\geq0}$ 
we see that
\[
\ldots\supseteq F^p\,\widehat{U}(\frakg(V))\supseteq F^{p+1}\,\widehat{U}(\frakg(V))\supseteq\ldots\quad\mbox{and}\quad
\cap_{p}F^p\,\widehat{U}(\frakg(V))=0.
\]
We set $\widehat{U}(\frakg(V))=\cup_pF^p\,\widehat{U}(\frakg(V))$, which is a complete algebra over $\C$ and satisfies
\begin{align*}
&F^p\,\widehat{U}(\frakg(V))\cdot F^q\,\widehat{U}(\frakg(V))\subseteq F^{p+q}\,\widehat{U}(\frakg(V)),\\
&[F^p\,\widehat{U}(\frakg(V)),F^q\,\widehat{U}(\frakg(V))]\subseteq F^{p+q}\,\widehat{U}(\frakg(V)).
\end{align*}

\begin{definition}\label{def:current-algebra}
 Let $\mathfrak{B}(V)$ be the two-sided ideal of $\widehat{U}(\frakg(V))$ generated by
$J_0(\ket{0})-1$ and the associativity relation for all $v_1\in V_{\Delta_1},\,v_2\in V_{\Delta_2}$ and $m,\,n\in\Z$, i.e., 
\begin{multline*}
J_m(J_n(v_1)v_2)-
\sum_{j=0}^\infty(-1)^j\binom{n+\Delta_1-1}{j}\\
\times
\left(
J_{n-j}(v_1)J_{m-n+j}(v_2)
-(-1)^{n+\Delta_1-1}J_{m+\Delta_1-j-1}(v_2)J_{j-\Delta_1+1}(v_1)
\right).
\end{multline*}
We set $\calu(V)=\widehat{U}(\frakg(V))/\mathfrak{B}(V)$, which is called
the \textit{current algebra} associated to the chiral vertex operator algebra $V$. 
\end{definition}

\begin{remark}
{\rm (1)} Let $F^p\,\calu(V)\,(p\in\Z)$ be the filtration induced by $F^p\,\widehat{U}(\frakg(V))$.
Then the current algebra $\calu(V)$ is a filtered algebra and a filtered Lie algebra.

\vskip 1ex
\noindent
{\rm (2)} The commutator and associativity formulas hold on $\calu(V)$ as well as $J_n(T(-1)v)=-(n+\Delta)J_n(v)$ and
$J_n(\ket{0})=\delta_{n,0}1$ for all $v\in V_\Delta$ and $n\in\Z$.
\end{remark}

\subsection{Modules for chiral vertex operator algebras}

\begin{definition}
Let $V$ be a chiral vertex operator algebra. A \textit{module} $M$ for $V$ is a $\calu(V)$-module such that

\vskip 1.5ex
\noindent
(a) $M$ is a finitely generated $\calu(V)$-module.

\vskip 1ex
\noindent
(b) For any $m\in M$ the vector space $F^0\calu(V)m$ is finite dimensional.

\vskip 1.5ex
\noindent
We denote by $\calm\!\operatorname{od}(V)$ the category of $V$-modules.
The chiral vertex operator algebra $V$ is a $\calu(V)$-module, and is generated by the vacuum $\ket{0}$.
Since $J_n(v)\in\Endc^{\,f}_{-n}\,V$ we see that  $F^0\calu(V)v\subseteq\bigoplus_{\Delta=0}^{|v|}V_\Delta$, in particular,
$F^0\calu(V)v$ is finite dimensional for any $v\in V$ so that
$V$ is a $V$-module, which is called the \textit{vacuum module}\,.
\end{definition}

Let $M$ be a $V$-module and $\{m_i\}$ be a \textit{finite} set of generators as a $\calu(V)$-module. 
The vector space $W=\sum_{i}F^0\mathcal{U}(V)m_i$ is finite dimensional by the condition (b). 
The set $F^0\mathcal{U}(V)$ is a subalgebra of $\calu(V)$, and $W$ is a $F^0\mathcal{U}(V)$-module. In particular, $W$ is
$T(0)$-invariant and is decomposed into a direct sum of generalized eigenspaces for the operator $T(0)$. 
Since $M=U(\frakg(V)_{<0}^f)W$  the $V$-module $M$ is also decomposed into a direct sum of generalized eigenspaces for
$T(0)$. For any $h\in\C$ we set 
\[
M_{(h)}=\{m\in M\;|\; (T(0)-h)^km=0\mbox{ for a positive integer $k$}\}.
\]
We note that $J_n(v):M_{(h)}\rightarrow M_{(h-n)}$ for all $v\in V_\Delta$ and $n\in\Z$.

\begin{proposition}\label{prop:cg}
Let $V$ be a chiral vertex operator algebra and $M$ be a $V$-module. 
Then there exists a finite set of complex numbers
$\{h_i\}$ such that $h_i-h_j\not\in\Z\,(i\neq j)$ and $M=\bigoplus_{h\in P(M)}M_{(h)}$, where $P(M)=\cup_{i}(h_i+\Z_{\geq0})$.
\end{proposition}

\begin{definition}
Given any $v\in V_\Delta$ we associate a formal power series 
\[
J^M(v,z)=\sum_{n\in\Z}J_n(v)z^{-n-\Delta}\in (\End\,M)[[z,z^{-1}]],
\]
which is called the \textit{nonabelian current} over $M$ associated to the vector $v$. 
\end{definition}

\section{Finiteness theorem for $V$-modules}
\label{sec:fin}
\subsection{Finiteness condition}

\begin{definition}
Let $V$ be a chiral vertex operator algebra and $M$ be a $V$-module. For any integer $n\geq2$ we denote by $C_n(M)$ the vector
subspace of $M$, which is linearly spanned by vectors $J_{-\Delta-p}(v)m$ 
for all $v\in V_{\Delta}\,(\Delta\geq0),\,m\in M$ and $p\geq n-1$.  If $\dim M/C_n(M)<\infty$ we say that the $V$-module $M$
satisfies the \textit{$C_n$-condition} or that $M$ is \textit{$C_n$-finite}\,. 
\end{definition}

We see that $C_2(M)\supseteq C_3(M)\supseteq\ldots\supset C_n(M)\supseteq\ldots$, 
and that if $M$ is $C_n$-finite for an $n\geq2$ then it is $C_2$-finite. Gaberdiel and Neitzke (\cite[Theorem 11]{GN}) showed that
that if a chiral vertex operator algebra $V$ is $C_2$-finite then it is $C_n$-finite for all $n\geq 2$. We shall prove that if $V$ is
$C_2$-finite then any $V$-module is $C_n$-finite for all $n\geq2$ (see Corollary \ref{cor:finiteproperty}). 

\subsection{Fermionic property of modules}
The current Lie algebra $\frakg(V)$ is a filtered Lie algebra with the filtration $F^p\frakg(V)\,(p\in\Z)$.
In this subsection we use another filtration which we denote by $G_p\frakg(V)\,(p\in\Z)$.

\begin{definition}
For any integer $p$ we set
\[
G_p\tilde{V}^{(1)}=\bigoplus_{\Delta\leq p}V_\Delta\otimes\C((\xi))(d\xi)^{1-\Delta}
\]
and denote by $G_p\frakg(V)$ the image of $G_p\tilde{V}^{(1)}$ under the canonical projection $\tilde{V}^{(1)}\rightarrow\frakg(V)$.
\end{definition}

We see that $G_p\frakg(V)\,(p\in\Z)$ defines an increasing filtration on $\frakg(V)$;
\[
0=G_{-1}\frakg(V)\subseteq\ldots\subseteq G_p\frakg(V)\subseteq G_{p+1}\frakg(V)\subseteq\ldots,
\quad\cup_p G_p\frakg(V)=\frakg(V).
\]
The commutator formula shows:
\begin{lemma}
$[G_p\frakg(V),G_q\frakg(V)]\subset G_{p+q-1}\frakg(V)$ for any $p,\,q\in\Z$, and in particular 
$\frakg(V)$ is a filtered Lie algebra with the filtration $G_p\frakg(V)\,(p\in\Z)$.
Then $\gr_\bullet^G\frakg(V)=\bigoplus_{p=0}^\infty\gr_p^G\frakg(V)$ is a graded abelian Lie algebra where
$\gr_p^G\frakg(V)=G_p\frakg(V)/G_{p-1}\frakg(V)$.
\end{lemma}

We set 
\[
G_p\,U(\frakg(V))=
\begin{cases}
\C&\quad (p<0),\\
\sum_{p_1+\cdots+p_r\leq p}G_{p_1}\frakg(V)\cdots G_{p_r}\frakg(V)&\quad (p\geq0),
\end{cases}
\]
which defines an increasing filtration on $U(\frakg(V))$. This filtration induces an increasing filtration on
$\widehat{U}(\frakg(V))$ since the filtration $G_p\,U(\frakg(V))\,(p\in\Z)$ is compatible with the filtration $F_N^qU(\frakg(V))$.
Let $G_p\,\calu(V)\,(p\in\Z)$ be the filtration on $\calu(V)$, which is induced by $G_p\widehat{U}(\frakg(V))$. 

\begin{proposition} 
The current algebra $\calu(V)$ is a filtered algebra, i.e., $G_p\,\calu(V)\cdot G_q\,\calu(V)\subseteq
G_{p+q}\,\calu(V)$  for any integers $p$ and $q$. Moreover  $[G_p\,\calu(V),G_q\,\calu(V)]\subseteq G_{p+q-1}\,\calu(V)$,
and in particular $\gr_\bullet^G\calu(V)=\bigoplus_{p=0}^\infty\gr_p^G\calu(V)$ is a commutative graded algebra 
where $\gr_p^G\calu(V)=G_p\calu(V)/G_{p-1}\calu(V)$. 
\end{proposition}

Let $G_pU(\frakg(V)^f)\,(p\in\Z)$ be the induced filtration, i.e., $G_pU(\frakg(V)^f)=G_p\widehat{U}(\frakg(V))\cap U(\frakg(V)^f)$.
Then $\gr_{\bullet}^GU(\frakg(V)^f)$ is a commutative graded algebra, and there is a canonical isomorphism
\[
\gr_{\bullet}^GU(\frakg(V)^f)\cong\gr_{\bullet}^GU(\frakg(V)^f_{<0})\otimes\gr_{\bullet}^GU(\frakg(V)^f_{\geq0}).
\]
The filtration $F^p\calu(V)\,(p\in\Z)$ induces a filtration on $\gr_{\bullet}^G\calu(V)$ so that $\gr_{\bullet}^G\calu(V)$ becomes
a filtered graded commutative algebra.

Let us introduce a \textit{good filtration} on a $V$-module $M$. 
By definition of a $V$-module there exists a finite dimensional $F^0\calu(V)$-submodule $W$ 
such that the canonical $\calu(V)$-module homomorphism
\[
\calu(V)\otimes_{F^0\calu(V)}W\rightarrow M\rightarrow 0
\]
is exact.  We set $G_{-1}M=0$ and $G_pM=G_p\,\calu(V)W$ for all nonnegative integer $p$, which defines an increasing filtration on
$M$ such that $G_p\,\calu(V)G_qM\subseteq G_{p+q}M$. In particular $\gr_{\bullet}^GM$ is a $\gr_{\bullet}^G\calu(V)$-module. 
The vector space $\gr_{\bullet}^GM$ is also a $\gr_{\bullet}^GU(\frakg(V)^f)$-module.
Note that $G_pU(\frakg(V)^f)W=G_pM$ for all $p$.

\begin{proposition}
The graded $\gr_{\bullet}^G\calu(V)$-module $\gr_{\bullet}^GM$ is generated by $\gr_0^GM=W$ as a
$\gr_\bullet^GU(\frakg(V)_{<0}^f)$-module.
\end{proposition}

In other words the canonical $\gr_\bullet^GU(\frakg(V)_{<0}^f)$-module homomorphism
\[
\gr_\bullet^GU(\frakg(V)_{<0}^f)\otimes_{\C}W\rightarrow\gr_{\bullet}^GM\rightarrow0
\]
is exact. Note that there exists a canonical graded algebra homomorphism
\[
\gr_\bullet^GU(\frakg(V)_{<0}^f)\cong S(\gr_\bullet^G\frakg(V)_{<0}^f)
\]
where $S(R)$ denotes the symmetric algebra of over a graded vector space $R$\,.

\begin{lemma}
For any $v\in C_2(V)$ the element of $\gr_\bullet^G\calu(V)$, whose representative is $J_{n}(v)\,(n\in\Z)$, trivially acts on
$\gr_\bullet^GM$.
\end{lemma}
\begin{proof}
We can assume that $v=J_{-\Delta_1-p}(v_1)v_2\,(p\geq1)$ with $v_1\in V_{\Delta_1}$ and
$v_2\in V_{\Delta_2}$. On the one hand $|v|=\Delta_1+\Delta_2+p$ and $J_{n}(v)G_qM\subset G_{\Delta_1+\Delta_2+p+q}M$, while
the right hand side of the associativity formula for $J_{n}(v)$ belongs to $G_{\Delta_1+\Delta_2}\calu(V)$,
and then $J_{n}(v)G_qM\subset G_{\Delta_1+\Delta_2+q}M$. 
\end{proof}

Let $\frakg(V)^{f,\,C_2}$ be a graded vector subspace of $\frakg(V)^f$, which is linearly spanned by $J_{n}(v)$ where $v\in C_2(V)$
and $n\in\Z$.  We set $\mathfrak{u}=\gr_\bullet^G\frakg(V)^f/\gr_\bullet^G\frakg(V)^{f,\,C_2}$. 
Then the action of $\gr_\bullet^GU(\frakg(V)^f)\cong S(\gr_\bullet^G\frakg(V)^f)$ on $\gr_\bullet^G M$ induces an action of 
the commutative algebra $S(\mathfrak{u})$ on
$\gr_\bullet^GM$. We now set $\mathfrak{u}_{<0}=\gr_\bullet^G\frakg(V)_{<0}/\gr_\bullet^G\frakg(V)^{f,\,C_2}_{<0}$ 
where $\frakg(V)^{f,\,C_2}_{<0}=\frakg(V)^{f,\,C_2}\cap\frakg(V)_{<0}$; the set $\mathfrak{u}_{<0}$ is canonically a vector subspace
of $\mathfrak{u}$ and $\gr_\bullet^GM=S(\mathfrak{u}_{<0})\gr_0^GM$.

Let $U$ be a graded subspace of $V$ such that $V=U\oplus C_2(V)\oplus\C\ket{0}$. 
Then the graded vector space $\mathfrak{u}_{<0}$ is linearly spanned by vectors $[J_n(v)]\,(n<0)$ for $v\in U$.
We note that $\gr_{-1}^GM=0$ and $\gr_0^GM=W$, while for any positive integer $p$ 
the space $\gr_p^GM$ is linearly spanned by vectors
\[
J_{-n_1}(v_1)\cdots J_{-n_r}(v_r)w\quad(w\in W,\,n_1\geq n_2\geq\ldots\geq n_r>0)
\]
where $v_1,\ldots,v_r\in U$ such that $\sum_{i=1}^r|v_i|=p$.

\begin{lemma}\label{lem:square}
Let $v_1\in V_{\Delta_1}$ and $v_2\in V_{\Delta_2}$.
For any integer $m$ there exist integers $c_{ij}\,(i,j\in\Z)$ such that 
\[
J_{-m}(v_1)J_{-m}(v_2)=J_{-2m}(J_{-\Delta_1}(v_1)v_2)
+\sum_{\substack{i+j=2m\\i\neq j}}c_{ij}J_{-i}(v_1)J_{-j}(v_2)\quad\mbox{in $\gr_\bullet^G\calu(V)$}.
\]
\end{lemma}
\begin{proof}
Since the associativity formula 
\[
J_{-2m}(J_{-\Delta_1}(v_1)v_2)
=\sum_{j=0}^\infty
\left(
J_{-\Delta_1-j}(v_1)J_{-2m+\Delta_1+j}(v_2)+J_{-2m+\Delta_1-j-1}(v_2)J_{j-\Delta_1+1}(v_1)
\right)
\]
holds for all integer $m$ we get for $m\geq\Delta_1$ 
\begin{multline*}
J_{-m}(v_2)J_{-m}(v_1)=J_{-2m}(J_{-\Delta_1}(v_1)v_2)\\
-\sum_{\substack{j\geq0\\j\neq m-\Delta_1}}J_{-\Delta_1-j}(v_1)J_{-2m+\Delta_1+j}(v_2)
-\sum_{j=0}^\infty J_{-2m+\Delta_1-j-1}(v_2)J_{j-\Delta_1+1}(v_1),
\end{multline*}
while for $m\leq\Delta_1-1$
\begin{multline*}
J_{-m}(v_1)J_{-m}(v_2)=J_{-2m}(J_{-\Delta_1}(v_1)v_2)\\
-\sum_{j=0}^\infty J_{-\Delta_1-j}(v_1)J_{-2m+\Delta_1+j}(v_2)
-\sum_{\substack{j\geq0\\j\neq\Delta_1-m-1}}J_{-2m+\Delta_1-j-1}(v_2)J_{j-\Delta_1+1}(v_1).
\end{multline*}
\end{proof}

The fermionic property of $V$-modules means:
\begin{theorem}\label{th:basis}
Let $V$ be a chiral vertex operator algebra and $M$ be a $V$-module.
Let $U$ be a graded vector subspace of $V$ such that $V=U\oplus C_2(V)\oplus\C\ket{0}$. 
For any nonnegative integer $p$ the vector space $\gr_p^GM$ is linearly spanned by vectors
\[
J_{-n_1}(v_1)\cdots J_{-n_r}(v_r)w\quad (w\in W,\,n_1>n_2>\ldots>n_r>0)
\]
where $v_1,\ldots,v_r\in U$ and $\sum_{i=1}^r|v_i|=p$.
\end{theorem}

In order to prove the theorem  we first show that the following statement S$(p,N)$ holds for any nonnegative integers $p$ and $N$;

\vskip 2ex
\noindent
\textbf{S$(p,N)$:} The vector space $\gr_p^GM$ is linearly spanned by vectors
\[
J_{-n_1}(v_1)\cdots J_{-n_r}(v_r)w\quad(w\in W, n_1\geq n_2\geq\ldots\geq n_r>0)
\]
where $v_1,\ldots,v_r\in U$ such that $\sum_{i=1}^r|v_i|=p$, and $n_i\neq n_j$ if $n_i,\,n_j\leq N\,(i\neq j)$.
\vskip 2ex

We will prove the statement S$(p,N)$ by induction with respect to the lexicographic order on $(p,N)\in(\Z_{\geq0})^2$.
The statement S$(0,N)$ holds for all $N$. We suppose that the statement holds for any pair which is strictly smaller than $(p,N)$.
Then by S$(p,N-1)$ any vector in $\gr_p^GM$ is a linear combination of vectors
\[
J_{-n_1}(u_1)\cdots J_{-n_r}(u_r)\;\mathbb{J}^s\;J_{-m_{1}}(w_1)\cdots J_{-m_t}(w_t)w,\quad w\in W
\]
where $n_1\geq n_2\geq\ldots\geq n_r>N>m_{1}>\ldots>m_t>0$ and $\mathbb{J}^s=J_{-N}(v_1)\cdots J_{-N}(v_s)$.
We set $\mathbb{I}^r=J_{-n_1}(u_1)\cdots J_{-n_r}(u_r)$ and $\mathbb{K}^t=J_{-m_{1}}(w_1)\cdots J_{-m_t}(w_t)$.
Applying Lemma \ref{lem:square} to $\mathbb{J}^s$ we see that
\[
\mathbb{I}^r\mathbb{J}^s\mathbb{K}^tw
=\mathbb{I}^r\mathbb{J}^{s-2}J_{-2N}(J_{-\Delta_{s-1}}(v_{s-1})v_s)\mathbb{K}^tw
+\sum_{\substack{i+j=2N\\i\neq j}}c_{ij}\mathbb{I}^r\mathbb{J}^{s-2}J_{-i}(v_{s-1})J_{-j}(v_{s})\mathbb{K}^tw.
\]
In the second term in the right hand side of the equation in the previous text we see that either $i>N$ or $j>N$. We
then can apply S$(p-1,N)$ to either the term $\mathbb{J}^{s-2}J_{-j}(v_{s})\mathbb{K}^tw$ or
$\mathbb{J}^{s-2}J_{-i}(v_{s-1})\mathbb{K}^tw$ in order to get the desired expression. 
Now in the first term we replace $J_{-\Delta_{s-1}}(v_{s-1})v_s$ by $u+C_2(V)$ for a suitable $u\in U$ so that  
$\mathbb{I}^r\mathbb{J}^{s-2}J_{-2N}(J_{-\Delta_{s-1}}(v_{s-1})v_s)\mathbb{K}^tw
=\mathbb{I}^rJ_{-2N}(u)\mathbb{J}^{s-2}\mathbb{K}^tw$ in $\gr_p^GM$. Now the induction on $s$ can be applied to prove 
the statement S$(p,N)$.

\vskip 1.5ex

We now give a proof of Theorem \ref{th:basis}.
Let $M=\bigoplus_{h\in P(M)}M_{(h)}$ and set $\alpha=\min\,\{\operatorname{Re}h\,|\,h\in P(M)\}$. 
We will show that for any $h\in P(M)$ any vector from $\gr_p^GM_{(h)}$ has the desired expression. We suppose that
\[
J_{-n_1}(v_1)\cdots J_{-n_r}(v_r)w\in \gr_p^GM_{(h)}\quad\mbox{and}\quad w\in W\cap M_{(h')},
\]
where $n_1\geq n_2\geq\ldots\geq n_r>0$.
Using S$(p,N)$ for a positive integer $N$ such that $N\geq\operatorname{Re}h-\alpha+1$ 
we can assume that $n_i\neq n_j\,(i\neq j)$ if $n_i,\,n_j\leq N$. 
Now $n_1+\cdots+n_r=h-h'$, and in particular $n_1\leq\operatorname{Re}(h-h')<N$.  So we get $n_i\neq n_j$ for all $i\neq j$.

\begin{corollary}\label{cor:finiteproperty}
Let $V$ be a $C_2$-finite chiral vertex operator algebra and $M$ be a $V$-module. 

\vskip 1.5ex
\noindent
{\rm (1)} $M$ is $C_n$-finite for all $n\geq2$.

\vskip 1ex
\noindent
{\rm (2)} Let $M=\bigoplus_{h\in P(M)}M_{(h)}$. Then $M_{(h)}$ is finite dimensional for all $h\in P(M)$.
\end{corollary} 
\begin{proof}
Let $U$ be a finite dimensional graded vector subspace of $V$ such that $V=U\oplus C_2(V)\oplus\C\ket{0}$. Any element of $M$ is a
linear combination of vectors
\[
v(\vec{n},w):=J_{-n_1}(v_1)\cdots J_{-n_r}(v_r)w\quad (w\in W,\,v_1,\ldots,v_r\in U,\,n_1>n_2>\ldots>n_r>0).
\]
If $n_1\geq|v_1|+n-1$ then $v(\vec{n},w)\in C_n(M)$. Otherwise $|v_1|+n\geq n_1>n_2>\ldots>n_r>0$, 
and for fixed $v_1,\ldots,v_r$ the choice of such
$(n_i)$ is finite. Now since $U$ is finite dimensional we see that $M/C_n(M)$ is finite dimensional.

Suppose that $w\in W\cap M_{(h_0)}$ and $v(\vec{n},w)\in M_{(h)}$.
Then we see that $h_0+n_1+\cdots+n_r=h$\,; only a finite number of $n_1>n_2>\ldots>n_r>0$ satisfy this condition so that 
the vector space $M_{(h)}$ is finite dimensional since $U$ is finite dimensional.
\end{proof}

We can characterize $V$-modules in terms of nonabelian currents as we have done for a chiral vertex operator algebra in
\S\,\ref{sub:2point}.  Let $M^*=\bigoplus_{h\in P(M)}M_{(h)}^*$ be the graded dual of a $V$-module $M$, and $\pairing{\;\,}{\;}$ be
the dual pairing between
$M^*$ and $M$. Fix $m\in M$ and $m^*\in M^*$. For any $v_1\in V_{\Delta_1}$ and $v_2\in V_{\Delta_2}$ the formal power series
$\bra{m^*}J^M(v_1,z_1)J^M(v_2,z_2)\ket{m}$
absolutely converges on the domain $|z_1|>|z_2|>0$ and is analytically continued to a rational function on $\CP^1\times\CP^1$ with
possible poles  on $z_1=z_2$, $z_i=0\,(i=1,2)$ and $z_i=\infty\,(i=1,2)$.   We use the same  notation to denote this rational
function. The formal power series
$\bra{m^*}J^M(J(v_1,z_1-z_2)v_2,z_2)\ket{m}$ absolutely converges on the domain $0<|z_1-z_2|<|z_2|$ and is analytically continued to
a rational function on $\CP^1\times\CP^1$: we have the $\cals_2$-symmetry and the operator product property
\begin{align*}
\bra{m^*}J^M(v_1,z_1)J^M(v_2,z_2)\ket{m}&=\bra{m^*}J^M(v_2,z_2)J^M(v_1,z_1)\ket{m}\\
&=\bra{m^*}J^M(J(v_1,z_1-z_2)v_2,z_2)\ket{m}\\
&=\bra{m^*}J^M(J(v_2,z_2-z_1)v_1,z_1)\ket{m}.
\end{align*}

\begin{remark}\label{rem:extension}
 Let $M$ be a graded $\frakg(V)^f$-module.
For $v\in V_\Delta$ we set $J^M(v,z)=\sum_{n\in\Z}J_n(v)z^{-n-\Delta}\in(\Endc^{\,f}\,M)[[z,z^{-1}]]$.
Then $\frakg(V)^f$-module structure is uniquely extended to a $\calu(V)$-module structure if and only if all nonabelian currents
$J^M(v,z)$ have the
$\cals_2$-symmetry, the operator product property, and $J_0(\ket{0})=\id_M$.
\end{remark}

\subsection{Zero-mode algebras}
We set $\mathfrak{J}=\calu(V)F^1\calu(V)$ and $\mathfrak{J}^0=\mathfrak{J}\cap F^0\calu(V)$. 
The  set $\mathfrak{J}$ is a left ideal of $\calu(V)$ and $\mathfrak{J}^0$ is a two-sided ideal of $F^0\calu(V)$. 
We note that $\mathfrak{J}^0=\sum_{n\geq1}F^{-n}\frakU(V)\cdot F^{n}\frakU(V)$.
We set $A_0(V)=F^0\calu(V)/\mathfrak{J}^0$, which is an associative algebra over $\C$ with a unit $1$.

\begin{definition}
Let $V$ be a chiral vertex operator algebra. We call $A_0(V)$ the \textit{zero-mode algebra} associated to $V$. 
We denote by $[J_0(v)]\,(v\in V)$ the element of $A_0(V)$ being represented by $J_0(v)$. 
The category of finite dimensional $A_0(V)$-modules is denoted by $\calm\!\operatorname{od}(A_0(V))$ . 
\end{definition}

\begin{proposition}\label{prop:zero-trivial}
{\rm (1)}  Let $U$ be a graded vector subspace of $V$ such that $V=U\oplus C_2(V)$. 
Then $A_0(V)$ is linearly spanned by $\{\,[J_0(v)]\,|\, v\in U\,\}$. 

\vskip 1ex
\noindent
{\rm (2)} Let $\mathfrak{A}$ be the vector subspace of $A_0(V)$, which is linearly spanned by $\{\,[J_0(v)]\,|\, v\in
V_\Delta\,(\Delta\geq1)\,\}$. Then $\mathfrak{A}$ is a two-sided ideal of $A_0(V)$, and the quotient algebra $A_0(V)/\mathfrak{A}$
is one dimensional $A_0(V)$-module which is generated by $[1]$ where $1=J_0(\ket{0})$.
\end{proposition}
\begin{proof}
By using the associativity and commutator formulas
we see that for any $v_1\in V_{\Delta_1}$ and $v_2\in V_{\Delta_2}$
\begin{multline}\label{eqn:zero11}
J_0(v_1)J_0(v_2)
\equiv J_0(J_{-\Delta_1}(v_1)v_2)\\
-\sum_{1-\Delta_1\leq j<0}\sum_{k=0}^{\Delta_1+\Delta_2-1}
\binom{\Delta_2-j-1}{k}J_0(J_{k-\Delta_2+1}(v_2)v_1)\mod\mathfrak{J}^0\,.
\end{multline}
Since $A_0(V)$ is linearly spanned by vectors $[J_0(v_1)]\cdots[J_0(v_n)]$ the induction on the number of products shows that
$A_0(V)$ is linearly spanned by vectors $[J_0(v)],\,v\in V$.

Now for any $v_1\in V_{\Delta_1}$ and $v_2\in V_{\Delta_2}$ and $p\geq1$
\begin{multline*}
J_0(J_{-\Delta_1-p}(v_1)v_2)\equiv (-1)^{p+\Delta_1+1}\binom{-p-1}{\Delta_1-1}J_0(v_2)J_0(v_1)\\
+\sum_{\substack{0<j<\Delta_1-1\\k\geq0}}(-1)^{p+j+\Delta_1-1}\binom{-p-1}{-j+\Delta_1-1}
\binom{\Delta_1+j-1}{k}J_0(J_{k-\Delta_2+1}(v_2)v_1)\mod\mathfrak{J}^0.
\end{multline*}
The induction on the filtration on $A_0(V)$ induced by  $F_pV=\oplus_{\Delta\leq p}V_\Delta\,(p\in\Z)$ 
applied to the equation in the previous text as well as (\ref{eqn:zero11}) shows that 
$A_0(V)$ is linearly spanned by $\{\,[J_0(v)]\,|\, v\in U\,\}$.

If $\Delta_1\geq 1$ and $\Delta_2\geq 1$ the nonzero multiple of $J_0(\ket{0})$ never appears in the right hand side of
(\ref{eqn:zero11}) so that the set $\mathfrak{A}$ is a two-sided ideal of $A_0(V)$.  
We now note that $\dim A_0(V)/\mathfrak{A}\leq1$ and $V_0=\C\ket{0}$ is an $A_0(V)$-module by (1). Then we see that
$A_0(V)/\mathfrak{A}\cong\C$.
\end{proof}

\begin{definition}
Let $M$ be a $V$-module. We set $\mathcal{HW}(M)=\{m\in M\,|\, \frakg(V)^f_{>0}m=0\}$ and
$\mathcal{CHW}(M)=M/\frakg(V)^f_{<0}M$, which is called the space of \textit{highest weight vectors} of $M$ and
the space of \textit{co-highest weight vectors} of $M$, respectively. 
\end{definition}

\begin{lemma}\label{lem:hwe} 
Let $M\neq0$ be a $V$-module. Then $\mathcal{HW}(M)\neq0$, and it is canonically a left $A_0(V)$-module.
\end{lemma}

\begin{theorem}\label{th:module}
Let $V$ be a $C_2$-finite chiral vertex operator algebra. 

\vskip 1.5ex
\noindent
{\rm (1)} The zero-mode algebra $A_0(V)$ is finite dimensional.

\vskip 1ex
\noindent
{\rm (2)} The number of simple objects of $\calm\!\operatorname{od}(A_0(V))$ is finite.

\vskip 1ex
\noindent
{\rm (3)}
The $A_0(V)$-module $\mathcal{HW}(M)$ is finite dimensional. In particular
\[
\mathcal{HW}:\mathcal{M}od\,(V)\rightarrow\mathcal{M}od\,(A_0(V))
\]
is a functor.

\vskip 1ex
\noindent
{\rm (4)} Let $W$ be a simple $A_0(V)$-module and set $M(W)=\calu(V)\otimes_{F^0\calu(V)}W$. Then $M(W)$ is a $V$-module,
 and there exists a unique maximal proper $V$-submodule $N(W)$ so that $L(W)=M(W)/N(W)$ is a simple $V$-module.

\vskip 1ex
\noindent
{\rm (5)} The map $W\mapsto L(W)$ define a bijective map between the set of equivalent classes of finite dimensional simple
$A_0(V)$-modules and the set of equivalent classes of simple $V$-modules.

\vskip 1ex
\noindent
{\rm (6)} Let $M$ be a $\calu(V)$-module such that $M=\bigoplus_{h\in P(M)}M_{(h)}$ where $M_{(h)}=\{m\in
M\;|\;(T(0)-h)^nm=0\mbox{ for some } n\in\Z_{>0}\}$.  
Then $M$ is a $V$-module if and only if $\dim M_{(h)}<\infty$ for all $h$.

\vskip 1ex
\noindent
{\rm (7)} For any $V$-module $M$ there exists a composition series\,{\rm ;}
\[
0=F_0M\subset F_1M\subset\ldots\subset F_nM=M
\]
where $F_iM\,(0\leq i\leq n)$ are $V$-modules and $F_iM/F_{i-1}M\,(1\leq i\leq n)$ are simple $V$-modules.

\end{theorem}
\begin{proof}
The statements (1) and (2) follows from Proposition \ref{prop:zero-trivial} (1). We give proofs for the statement (3)-(6). The
existence of a composition series (7) can be proved by the standard argument.

(3) Suppose that $\mathcal{HW}(M)\cap M_{(h)}\neq0$. 
Then there exists a simple $A_0(V)$-submodule of $\mathcal{HW}(M)$ on which $T(0)$ acts as
$h\id$. However the number of such an $h$ is finite since $A_0(V)$ is finite dimensional.

(4) Since $W$ is simple $A_0(V)$-module $T(0)=h\id_W$ for some $h\in\C$. Then $M(W)=\bigoplus_{n=0}^\infty M_{h+n}$
where $M_{h'}=\{m\in M\;|\;T(0)m=h'm\;\}$.
Let $N$ be a proper $V$-submodule of $M(W)$. Then $N=\bigoplus_{n\geq1}N_{h+n}$ otherwise $N=M(W)$ since $W$ is simple, so that
the sum of proper $V$-submodules is again proper. We now set $N(W)$ to be the sum of all proper
$V$-submodules of $M(W)$. Then $N(W)$ is the unique maximal proper $V$-submodule of $M(W)$.

(5) Let $M$ be a simple $V$-module and $W$ be a simple $A_0(V)$-module in $M$; such a simple module exists as shown in (3).
Since $M$ is simple we see that the canonical $\calu(V)$-module homomorphism $\calu(V)\otimes_{F^0\calu(V)}W\rightarrow
M\rightarrow0$ is exact. The kernel of this sequence is the maximal proper $\calu(V)$-submodule since $M$ is simple, and then
$L(W)\cong M$.

(6) We suppose that $M=\bigoplus_{h\in P(M)}M_{(h)},\,\dim M_{(h)}<\infty$.
For any $m\in M_{(h)}$ we see that $F^0\calu(V)m\subset\bigoplus_{n=0}^\infty M_{(h-n)}$, 
and that $F^0\calu(V)m$ is finite dimensional since $M_{(h-n)}=0$ for sufficiently large $n$.

Let $S$ be the maximal finite set of complex numbers such that each of whose elements is a generalized eigenvalue of $T(0)$ on
a finite dimensional
$A_0(V)$-module. We set $W=\bigoplus_{h\in P(M)\cap S}M_{(h)}$, and let $N$ be the $\calu(V)$-submodule of $M$ generated by $W$.
Since the $V$-submodule $N$ is finitely generated it suffices to show that $M=N$.
We now suppose that $N$ is a proper $\calu(V)$-submodule of $M$.
By definition $M/N=\bigoplus_{h\in P(M)\setminus S}(M/N)_{(h)}$. 
Then there exists a nonzero $A_0(V)$-module in $M/N$, and we see that $P(M/N)\cap S\neq\emptyset$ which is a contradiction.
\end{proof}

\section{Duality}\label{sec:duality}
Let $M$ be a $V$-module and $M^*$ be its graded dual. Following Frenkel-Huang-Lepowsky \cite{FHL} the $\frakg(V)^f$-module
structure on $M^*$ is introduced. It is proved that the graded dual $M^*$ becomes a $V$-module in our sense if $V$ is $C_2$-finite.

\subsection{Anti-automorphism for current Lie algebras}
\label{sub:theta}
We define a linear map $\theta:\frakg(V)^f\rightarrow\frakg(V)^f$ by setting
\[
\theta(J_n(v))=(-1)^\Delta\sum_{j=0}^\infty\frac{1}{j!}J_{-n}(T(1)^jv)\quad\mbox{for $v\in V_\Delta$ and $n\in\Z$}.
\]

\begin{proposition}
The linear map $\theta:\frakg(V)^f\rightarrow\frakg(V)^f$ is an anti-Lie algebra involution.
\end{proposition}

In the following we will prove the proposition. The fact that $\theta$ defines a well-defined linear map is not obvious; 
one has to show that $\theta$ preserves the relation $J_n(T(-1)v)= -(n+\Delta)J_n(v)$ for all 
$v\in V_\Delta$ and $n\in\Z$. But it follows from the formula $[T(1)^j,T(-1)]=2jT(1)^{j-1}T(0)-j(j-1)T(1)^{j-1}$;
note that the set $\{T(-1),\,T(0),\, T(1)\}$ forms the Lie algebra $sl_2(\C)$.

We set $J^\theta(v,z)=\sum_{n\in\Z}\theta(J_n(v))z^{-n-\Delta}$ so that 
$J^\theta(v,z)=J(e^{zT(1)}(-z^{-2})^{T(0)}v,z^{-1})$.
The fact that $\theta$ is an anti-Lie algebra homomorphism is a consequence of the formal conjugation formula due to 
\cite[Lemma 5.2.3]{FHL}:
\begin{lemma}\label{lem:fhl}
{\rm (1)} As elements of $\frakg(V)[[z,z^{-1}]]$
\[
\mu^{T(0)}J(v,z)\mu^{-T(0)}=J(\mu^{T(0)}v,\mu z)\quad\mbox{for}\quad\mu\in\C^\times.
\]
\noindent
{\rm (2)} As elements of $\frakg(V)[[z,z^{-1},t\,]]$
\[
e^{tT(1)}J(v,z)e^{-tT(1)}=J\left(e^{t(1-tz)T(1)}(1-tz)^{-2T(0)}v,z/(1-tz)\right),
\]
where the binomial expression in the right hand side is expanded as 
\[
(1-tz)^n=\sum_{j=0}^\infty(-1)^j\binom{n}{j}t^jz^j\quad\mbox{for all $n\in\Z$}.
\]
\end{lemma}

We will prove the formula
\begin{equation}
\begin{split}\label{eqn:theta11}
[J^\theta(v_2,z_2),&\theta(J_m(v_1))]\\
&=\underset{z_1=z_2}{\operatorname{Res}}J^\theta(J(v_1,z_1-z_2)v_2,z_2)z_1^{m+\Delta_1-1}\,dz_1\\
&=\sum_{j=0}^{\Delta_1+\Delta_2-1}\binom{m+\Delta_1-1}{j}J^\theta(J_{j-\Delta_1+1}(v_1)v_2,z_2)z_2^{m+\Delta_1-j-1}.
\end{split}
\end{equation}
Then applying the commutator formula to the right hand side we see that 
\[
[\theta(J_n(v_2)),\theta(J_m(v_1))]=\theta([J_m(v_1), J_n(v_2)]),
\]
i.e., the linear map $\theta$ is an anti-homomorphism.

To prove (\ref{eqn:theta11}) we first note that that the commutator formula is equivalent to the following formal equation in
$\frakg(V)^f[[z_1,z_1^{-1},z_2,z_2^{-1}]]$
\[
[J_m(v_1),J(v_2,z_2)]
=\underset{z_1=z_2}{\operatorname{Res}}J(J(v_1,z_1-z_2)v_2,z_2)z_1^{m+\Delta_1-1}\,dz_1,
\]
where the right hand side is defined by $\underset{w=0}{\operatorname{Res}}J(J(v_1,w)v_2,z_2)(z_2+w)^{m+\Delta_1-1}\,dw$.
We now see that
\[
\begin{split}
&[J^\theta(v_2,z_2),\theta(J_m(v_1))]\\
&=(-1)^{\Delta_1+\Delta_2}\sum_{i,\,j=0}^\infty\frac{1}{i!j!}
z_2^{j-2\Delta_2}[J(T(1)^jv_2,z_2^{-1}),J_{-m}(T(1)^iv_1)]\\
&=(-1)^{\Delta_1+\Delta_2+1}\sum_{i,\,j=0}^\infty\frac{z_2^{j-2\Delta_2}}{i!j!}
\underset{z_1=z_2^{-1}}{\operatorname{Res}}J(J(T(1)^iv_1,z_1-z_2^{-1})T(1)^jv_2,z_2^{-1})
z_1^{-m+\Delta_1-i-1}\,dz_1\\
&=(-1)^{\Delta_1+\Delta_2+1}z_2^{-2\Delta_2}
\underset{z_1=z_2^{-1}}{\operatorname{Res}}J(J(e^{z_1^{-1}T(1)}v_1,z_1-z_2^{-1})e^{z_2T(1)}v_2,z_2^{-1})
z_1^{-m+\Delta_1-1}\,dz_1.
\end{split}
\]
Using the conjugation formula ($z=z_2,\,t=z_2^{-1}(1-z_1/z_2) $) we see that
\[
\begin{split}
J^\theta(J(v_1,z_1-z_2)v_2,z_2)&=(-1)^{\Delta_1+\Delta_2}z_2^{-2\Delta_1-2\Delta_2}
J(e^{z_2T(1)}J(v_1,z_2^{-1}(1-z_1/z_2))v_2,z_2^{-1})\\
&=(-1)^{\Delta_1+\Delta_2}z_1^{-2\Delta_1}z_2^{-2\Delta_2}
J(J(e^{z_1T(1)}v_1,z_1^{-1}-z_2^{-1})e^{z_2T(1)}v_2,z_2^{-1}).
\end{split}
\]
We now get the desired formula after the change of variables $z_1\rightarrow z_1^{-1}$.

By very definition of $J^\theta(v,z)$ we see that
\[
(J^\theta)^\theta(v,z)=J(e^{zT(1)}(-z^{-2})^{T(0)}e^{z^{-1}T(1)}(-z^{2})^{T(0)}v,z),
\]
and that $(J^\theta)^\theta(v,z)=J(v,z)$ by $(-z^2)^{T(0)}e^{zT(1)}(-z^2)^{-T(0)}=e^{-z^{-1}T(1)}$ (see \cite[5.3.1]{FHL}),
so $\theta^2(J_n(v))=J_n(v)$, i.e., $\theta$\, is an involution.

\subsection{Duality operation for modules}
Let $V$ be a $C_2$-finite chiral vertex operator algebra. Recall that any $V$-module $M$ has the decomposition 
$M=\bigoplus_{h\in P(M)}M_{(h)}$ where $\dim\,M_{(h)}<\infty$ for all $h\in P(M)$. We set $D(M)=\bigoplus_{h\in P(M)}M_{(h)}^*$ and
introduce a $\frakg(V)^f$-module structure on $D(M)$ by
$\langle J_n(v)\varphi,m\rangle=\langle\,\varphi,\theta(J_n(v))m\rangle$ for all $\varphi\in D(M)$ and $m\in M$. 

\begin{proposition}\label{prop:contra2}
Let $V$ be a $C_2$-finite chiral vertex operator algebra.

\vskip 1.5ex
\noindent
{\rm (1)} The $\frakg(V)^f$-module structure on $D(M)$ is uniquely extended to a $\calu(V)$-module structure,
and $D(M)$ is a $V$-module.

\vskip 1ex
\noindent
{\rm (2)} $D:\calm\!\operatorname{od}(V)\rightarrow\calm\!\operatorname{od}(V)$ is a contravariant functor such that 
$D(D(M))=M$ for any $M\in\ob(\calm\!\operatorname{od}(V))$.
\end{proposition}

By Remark \ref{rem:extension} it suffices to verify the $\cals_2$-symmetry, the operator product property and $J_0(\ket{0})=\id$\,.
Then $D(M)$ is a $V$-module by Theorem \ref{th:module} (6).
Let $v_1\in V_{\Delta_1}$ and $v_2\in V_{\Delta_2}$. 
Using the $\cals_2$-symmetry of a $V$-modules we see that the formal power series
$\bra{m^*}J^\theta(v_1,z_1)J^\theta(v_2,z_2)\ket{m}\in\C[[z_1,z_1^{-1},z_2,z_2^{-1}]]$ absolutely converges 
on the domain $|z_2|>|z_1|>0$, and is analytically continued to a
rational function on $\CP^1\times\CP^1$, whose poles are located at $z_1=z_2$ and $z_i=0,\, z_j=\infty\,(i,j=1,2)$, and that
$\mathcal{S}_2$-symmetry follows from the one of the $V$-module $M$. 
The operator product property
\[
\bra{m^*}J^\theta(v_1,z_1)J^\theta(v_2,z_2)\ket{m}
=\bra{m^*}J^\theta(J(v_1,z_1-z_2)v_2,z_2)\ket{m}
\]
on $0<|z_1-z_2|<|z_1|$ also holds by virtue of  
\[
e^{z_2T(1)}(-z_2^{-2})^{T(0)}J(v_1,z_1-z_2)(-z_2^{-2})^{-T(0)}e^{-z_2T(1)}
=J\left(e^{z_1T(1)}(-z_1^{-2})^{T(0)}v_1,z_1^{-1}-z_2^{-1}\right),
\]
which is derived by Lemma \ref{lem:fhl}, and the operator product property for $(J,M)$. 
Finally $\theta(J_0(\ket{0}))=J_0(\ket{0})$ so that $J_0(\ket{0})=\id$ on $D(M)$.

\begin{remark}\label{re:theta-asso}
The operator product property gives rise to the identity
\[
J^\theta(J^{\leq0}(v_1,z_1-z_2)v_2,z_2)=J^{\theta,>0}(v_1,z_1)J^{\theta}(v_2,z_2)+J^{\theta}(v_2,z_2)J^{\theta,\leq0}(v_1,z_1)
\]
on the domain $0<|z_1-z_2|<|z_1|,\,|z_2|$.
\end{remark}

\section{Conformal blocks}\label{sec:conformalblocks}
\subsection{Meromorphic forms and Lie algebras}\label{sub:lie}
Let $\CP^1=\C\cup\{\infty\}$ be the projective line and $z$ be its canonical affine coordinate. 
We denote by $H^0(\CP^1,\Omega^{1-\Delta}(*))$ the vector space of global meromorphic $(1-\Delta)$-forms on $\CP^1$, and
set $\mathcal{H}(V)^{(1)}=\bigoplus_{\Delta=0}^\infty V_\Delta\otimes H^0(\CP^1,\Omega^{1-\Delta}(*))$ and
$\mathcal{H}(V)^{(0)}=\bigoplus_{\Delta=0}^\infty V_\Delta\otimes H^0(\CP^1,\Omega^{-\Delta}(*))$. Define a linear map
$\nabla:\mathcal{H}(V)^{(0)}\rightarrow\mathcal{H}(V)^{(1)}$ by
\[
v\otimes f(z)(dz)^{-\Delta}\longmapsto T(-1)v\otimes f(z)(dz)^{-\Delta}+v\otimes\frac{df(z)}{dz}(dz)^{1-\Delta}\quad (v\in V_\Delta),
\]
and set $\frakg(V;\CP^1,*)=\mathcal{H}(V)^{(1)}/\nabla\mathcal{H}(V)^{(0)}$.
Introduce a bilinear operation
$[\;\,,\;]:\mathcal{H}(V)^{(1)}\rightarrow\mathcal{H}(V)^{(1)}$ by
\begin{multline*}
[v_1\otimes f_1(z)(dz)^{1-\Delta_1},v_2\otimes f_2(z)(dz)^{1-\Delta_2} ]\\
=\sum_{m=0}^\infty\frac{1}{m!}J_{m-\Delta_1+1}(v_1)v_2\otimes\frac{d^mf_1(z)}{dz^m}f_2(z)(dz)^{2-\Delta_1-\Delta_2+m}
\end{multline*} 
Using the exactly same argument given in the proof of
Proposition \ref{prop:lie} we get:
\begin{proposition}
The vector space $\frakg(V;\CP^1,*)$ is a Lie algebra with the bracket $[\;\,,\;]$.
\end{proposition}

For any point $w\in\CP^1$ we define $\xi_w=z-w$ for $w\neq\infty$ and $\xi_\infty=z$, and set  
\[
\widehat{\mathcal{K}}_w
=\begin{cases}
\C((\xi_w)&\quad(w\neq\infty),\\
\C((\xi_\infty^{-1})&\quad(w=\infty).
\end{cases}
\]
We note that the local parameter $\xi_\infty$ has the pole at $z=\infty$ of order $1$.

\begin{definition}
(1) For any $w\in\CP^1$ we define the linear map 
$i_w:H^0(\CP^1,\Omega^\Delta(*))\rightarrow \widehat{\mathcal{K}}_w(d\xi_w)^{\Delta}$ 
by taking the Laurent expansion at $z=w$ in terms of the variable $\xi_w$. 
We often use the notation $i_wf(z)(dz)^{\Delta}=f_w(\xi_w)(d\xi_w)^{\Delta}$.

\vskip 1ex
\noindent
(2) For any $w\in\CP^1\!\setminus\!\{\infty\}$ we define the linear map $j_w:\frakg(V;\CP^1,*)\rightarrow \frakg(V)$ by
$v\otimes f(z)(dz)^{1-\Delta}\mapsto v\otimes f_w(\xi_w)(d\xi_w)^{1-\Delta}$, while for $w=\infty$ we define 
the linear map $j_\infty$ by 
$v\otimes f(z)(dz)^{1-\Delta}\mapsto -\theta\left(v\otimes f_\infty(\xi_\infty)(d\xi_\infty)^{1-\Delta}\right)$,
where $\theta\,$ is the anti-Lie algebra involution of $\frakg(V)^f$ defined in \S\,\ref{sub:theta}. 
More precisely, let $f_\infty(\xi_\infty)=\sum_{n\leq n_0}f_{\infty,\,n}\,\xi_\infty^{n+\Delta-1}$. 
Then we define
\[
j_\infty\left(v\otimes f(z)(dz)^{1-\Delta}\right)=-\sum_{n\leq
n_0}f_{\infty,\,n}\theta\left(J_n(v)\right).
\]
Recall that $\theta\left(J_n(v)\right)=(-1)^{\Delta}J_{-n}(e^{T(1)}v)$ 
so that the infinite sum in the previous text defines an element of $\frakg(V)$. 
\end{definition}

\begin{proposition}\label{prop:lieiso}
For any $w\in\CP^1$ the linear map $j_w:\frakg(V;\CP^1,*)\rightarrow \frakg(V)$ is a Lie algebra homomorphism.
\end{proposition}

\begin{remark}\label{re:residue102}
Let $v\in V_\Delta,\,f(z)(dz)^{1-\Delta}\in H^0(\CP^1,\Omega^{1-\Delta}(*))$ and $M$ be a $V$-module. Then 
\[
j_\infty\left(v\otimes f(z)(dz)^{1-\Delta}\right)u
=\underset{\xi_\infty=\infty}{\operatorname{Res}}J^{\theta}\left(v,\xi_\infty\right)f_\infty(\xi_\infty)u\,d\xi_\infty
\]
for all $u\in M$ since $\underset{\xi_\infty=\infty}{\operatorname{Res}}g(\xi_\infty)d\xi_\infty=-g_{-1}$ for 
$g(\xi_\infty)=\sum_{n\ll0}g_{n}\,\xi_\infty^{n}\in\widehat{\mathcal{K}}_\infty$\,. 
\end{remark}

\subsection{Space of covacua and conformal blocks}
Let $A=\{1,\,2,\,\ldots,\,N,\,\infty\}$, and fix a set of distinct $N+1$ points $w_a\in \CP^1\,(a\in A)$ with $w_\infty=\infty$, 
which is denoted by $w_A=(w_a)_{a\in A}$. We often identify the index set $A$ with the set of distinct points $w_A$ by the
correspondence $a\leftrightarrow w_a$. As in \S\,\ref{sub:lie} we write $\,\xi_a=z-w_a\,(a\neq\infty)$ and $\,\xi_\infty=z$.

Let $M^a\,(a\in A)$ be $V$-modules.
Define $\frakg_{A}(V)=\bigoplus_{a\in A}\frakg(V)_a$, where each $\frakg(V)_a\,(a\in A)$ is a copy of the current Lie algebra
$\frakg(V)$. The vector space $\frakg_{A}(V)$ is naturally a Lie algebra,
and the tensor product vector space $M_A=\bigotimes_{a\in A}M^a$ becomes a $\frakg_A(V)$-module by
$\rho_A:\frakg_A(V)\rightarrow\End\,M_A$, where
$\rho_A=\sum_{a\in A}\rho_a$ and $\rho_a\,(a\in A)$ is the action of $\frakg(V)_a$ on the $a$-th component of $M_A$. 

Let $H^0(\CP^1,\Omega^{1-\Delta}(*{w_A}))$ be the vector space of global meromorphic $(1-\Delta)$-forms on $\CP^1$, 
whose poles are only located at $w_a\,(a\in A)$. We set 
\[
\frakg(V;\CP^1,*w_A)=\bigoplus_{\Delta=0}^\infty V_\Delta\otimes H^0(\CP^1,\Omega^{1-\Delta}(*{w_A}))/
\nabla\left(\bigoplus_{\Delta=0}^\infty V_\Delta\otimes H^0(\CP^1,\Omega^{-\Delta}(*{w_A}))\right).
\]
Then $\frakg(V;\CP^1,*w_A)$ is a Lie subalgebra of $\frakg(V;\CP^1,*)$, and 
the linear map 
\[
j_{w_a}:\frakg(V;\CP^1,*w_A)\rightarrow\frakg(V)_a
\]
is a Lie algebra homomorphism for all $a\in A$.   We denote by $\frakg_{w_A}^{out}(V)$ the image of the Lie algebra homomorphism
$j_{w_A}=\sum_{a\in A}j_{w_a}:\frakg(V;\CP^1,*w_A)\rightarrow\frakg_A(V)$ , which is a Lie subalgebra of $\frakg_A(V)$.

\begin{definition}
We call the vector space $\calv_{w_A}(M_A)=M_A/\frakg_{w_A}^{out}(V)M_A$ the \textit{space of covacua at} $w_A$, and call the dual space 
$\calv_{w_A}^\dagger(M_A)=\Hom_\C(\calv_{w_A}(M_A),\C)$ the \textit{space of conformal blocks at $w_A$}.
\end{definition}

\subsection{Systems of current correlation functions}
We give the definition of system of current correlation functions, and state several properties of them.
Let $A=\{1,\,2,\,\ldots,\,N,\,\infty\}$ and fix
$w_A=(w_a)_{a\in A}\in (\CP^1)^{N+1}$ such that $w_a\neq w_b\,(a\neq b)$ with $w_\infty=\infty$. For any nonnegative integer $m$ we
set
$[0]=\emptyset$, and
$[m]=\{\,1,\,\ldots,\,m\,\}$ for $m\geq1$.   For all $(z_i)_{i\in[m]}\in(\CP^1)^m$ we define divisors 
$D_{a,\,i}=\{\,z_i=w_a\,\}$ and $D_{i,\,j}=\{\,z_i=z_j\,\}\,(i\neq j)$. 
Let $M^a\,(a\in A)$ be $V$-modules and set $M_{[m],\,A}=\overset{m}{\overbrace{V\otimes\cdots\otimes V}}\otimes M_A$.

\begin{definition}
The set $\{\Phi_m\}_{m=0}^\infty$ is called a \textit{system of current correlation functions at
$w_A$}  if it satisfies the following set of axioms:

\vskip 1.5ex
\noindent
(a) $\Phi_m$ is a multilinear functional on $M_{[m],\,A}$ with values on meromorphic forms on $(\CP^1)^m$ , i.e.,
for any $v_i\in V_{\Delta_i}\,(1\leq i\leq m)$ and $u\in M_A$ its value $\Phi_m(\bigotimes_{i=1}^mv_i\otimes u)$ is a meromorphic
$(\Delta_1,\ldots,\Delta_m)$-form on $(z_i)_{i\in[m]}\in(\CP^1)^m$.
In particular $\Phi_0\in\Hom_\C(M_A,\C)$, which is called the \textit{initial term} of the system. 

\vskip 1.5ex
\noindent
We denote $\Phi_m(\bigotimes_{i=1}^mv_i\otimes u)$ by $\Phi_m(v_1,\ldots,v_m;u)_{(z_1,\ldots,z_m)}$ when we need to indicate the
dependence of the variables $z_1,\ldots,z_m$.

\vskip 1.5ex
\noindent
(b) $\Phi_m(v_1,\ldots,v_m;u)_{(z_1,\ldots,z_m)}\in 
H^0((\CP^1)^m,\Omega^{\boxtimes_{i=1}^m\Delta_i}(\sum_{i,\,a}*D_{a,\,i}+\sum_{i\neq j}*D_{i,\,j}))$ for all positive integers $m$.

\vskip 1ex
\noindent
(c) $\Phi_m$ is $\mathcal{S}_m$-invariant, i.e., for any $\sigma\in\mathcal{S}_m$  
\[
\Phi_m(v_{\sigma(1)},\,v_{\sigma(2)},\ldots,v_{\sigma(m)};u)_{(z_{\sigma(1)},\,z_{\sigma(2)}\ldots,z_{\sigma(m)})}
=\Phi_m(v_1,\,v_2,\ldots,v_m;u)_{(z_1,\,z_2,\ldots,z_m)}.
\]
\vskip 1ex
\noindent (d) The meromorphic form $\Phi_{m+1}(v,v_1,\ldots,v_m;u)_{(z,z_1,\ldots,z_m)}\,(v\in V_{\Delta})$ has the
following Laurent expansion at $z=w_a\,(a\in A)$ and $z=z_i\,(1\leq i\leq m)$:
\[
\begin{split}
&\Phi_{m+1}(v,v_1,\ldots,v_m;u)_{(z,z_1,\ldots,z_m)}\\
&\\
&\qquad=\begin{cases}
\displaystyle{\sum_{n\in\Z}\Phi_m\left(v_1,\ldots,v_m;\rho_a(J_n(v))u\right)_{(z_1,\ldots,z_m)}(z-w_a)^{-n-\Delta}(dz)^{\Delta}}&\quad
(a\neq\infty),\\
\displaystyle{\sum_{n\in\Z}\Phi_m\left(v_1,\ldots,v_m;\rho_\infty(\theta(J_{n}(v)))u\right)_{(z_1,\ldots,z_m)}
z^{-n-\Delta}\,(dz)^{\Delta}}&\quad(a=\infty),
\end{cases}
\end{split}
\]
and
\begin{multline*}
\Phi_{m+1}(v,v_1,\ldots,v_m;u)_{(z,z_1,\ldots,z_m)}\\
=\sum_{n\in\Z}\Phi_m\left(v_1,\ldots,J_n(v)v_i,\ldots,v_m;u\right)_{(z_1,\ldots,z_m)}(z-z_i)^{-n-\Delta}(dz)^{\Delta}.
\end{multline*}

\noindent
{\rm (e)} For any $1\leq i\leq m$
\[
\frac{\partial}{\partial z_i}\Phi_m(v_1,\ldots,v_m;u)_{(z_1,\ldots,z_m)}dz_i
=\Phi_m(v_1,\ldots, T(-1)v_i,\ldots,v_m ;u)_{(z_1,\ldots,z_m)}.
\]

\vskip 1ex
We denote by $\calc\!\operatorname{or}_{w_A}(M_A)$ the set of all systems of current correlation functions at $w_A$.
\end{definition}

\begin{note}
We often denote $\Phi_m(v_1,\ldots,v_m;u)_{(z_1,\ldots,z_m)}$ by
\[
\bra{\Phi}J(v_1,z_1)\cdots J(v_m,z_m)\ket{u}(dz_1)^{\Delta_1}\cdots(dz_m)^{\Delta_m}\,.
\]
\end{note}

\begin{theorem}\label{th:syscorr}
The linear map $\calc\!\operatorname{or}_{w_A}(M_A)\rightarrow\Hom_\C(M_A,\C)\,(\{\Phi_m\}_{m=0}^\infty\mapsto\Phi_0)$
induces an isomorphism $\calc\!\operatorname{or}_{w_A}(M_A)\cong\calv_{w_A}^\dagger(M_A)$.
\end{theorem}
\begin{proof}
We first show that for any $\{\Phi_m\}_{m=0}^\infty\in\calc\!\operatorname{or}_{w_A}(M_A)$ the initial term $\Phi_0$  belongs to
the space of conformal blocks at $w_A$, i.e., $\Phi_0(\frakg_{w_A}^{out}(V)M_A)=0$. Let $v\in V_\Delta$ and $f(z)(dz)^{1-\Delta}\in
H^0(\CP^1,\Omega^{1-\Delta}(*w_A))$. 
By the property (d) we see that $\Phi_1(v;u)_z=\Phi_0(\rho_a\left(J(v,\xi_a\right))u)(d\xi_a)^{\Delta}$ and
$\Phi_1(v;u)_z=\Phi_0(\rho_\infty\left(J^{\theta}(v,\xi_\infty\right))u)(d\xi_\infty)^{\Delta}$ in a
neighborhood of $z=w_a\,(a\in A\setminus\{\infty\})$
and  $z=\infty$, respectively. So using Remark \ref{re:residue102} we see that
\[
\begin{split}
&\Phi_0(\rho_A\left(j_{w_A}(v\otimes f(z)(dz)^{1-\Delta})\right)u)\\
&=\sum_{a\in A\setminus\{0\}}\Res{\xi_a}\Phi_0(\rho_a\left(J(v,\xi_a)\right)u)f_a(\xi_a)\,d\xi_a
+\underset{\xi_\infty=\infty}{\operatorname{Res}}
\Phi_0\left(\rho_\infty\left(J^{\theta}\left(v,\xi_\infty\right)\right)u\right)f_\infty(\xi_\infty)\,d\xi_\infty\\
&=\sum_{a\in A}\underset{z=w_a}{\operatorname{Res}}\,\Phi_1(v;u)_zf(z)\,(dz)^{1-\Delta}=0.
\end{split}
\]
since $\Phi_1(v,z)f(z)(dz)^{1-\Delta}$ is a meromorphic $1$-form whose poles are located only at
$z=w_a\,(a\in A)$, and the total sum of residues is zero. Thus the linear map
$\calc\!\operatorname{or}_{w_A}(M_A)\rightarrow\calv_w^\dagger(M_A)$ is well-defined, and is obviously injective;
the Laurent expansion uniquely determines the meromorphic form.

In order to prove that the linear map is surjective we need the notions of $1$-point and $2$-point current correlation functions.
The proof of the surjectivity shall be given in \S\, \ref{sub:cor-con} after we study these current correlation functions.
\end{proof}

\subsection{$1$-point current correlation functions}
For any $\Phi\in\calv_{w_A}^\dagger(M_A)$ we will construct $\Phi_1(v;u)_{(z)}$ for the system of
current correlation functions with the initial term $\Phi$.

We set $\widehat{\mathcal{K}}_a=\widehat{\mathcal{K}}_{w_a}\,(a\in A)$, and define the pairing between
$\widehat{\mathcal{K}}_A^\Delta=\bigoplus_{a\in A}\widehat{\mathcal{K}}_a(d\xi_a)^\Delta$ and 
$\widehat{\mathcal{K}}_A^{1-\Delta}=\bigoplus_{a\in A}\widehat{\mathcal{K}}_a(d\xi_a)^{1-\Delta}$
by 
\[
\langle
\sum_{a\in A}f_a(\xi_a)(d\xi_a)^\Delta,\sum_{a\in A}g_a(\xi_a)(d\xi_a)^{1-\Delta}
\rangle
=\sum_{a\in A}\operatorname{Res}_af_a(\xi_a)g_a(\xi_a)d\xi_a
\]
where
\[
\operatorname{Res}_af_a(\xi_a)d\xi_a=\begin{cases}
\Res{\xi_a}f_a(\xi_a)d\xi_a&\quad (a\neq\infty),\\
\underset{\xi_\infty=\infty}{\operatorname{Res}}f_\infty(\xi_\infty)d\xi_\infty&\quad(a=\infty).
\end{cases}
\]
Recall that
\[
i_{w_a}:H^0(\CP^1,\Omega^\Delta(*{w_A}))\rightarrow\widehat{\mathcal{K}}_a(d\xi_a)^\Delta\mbox{ and }
i_{w_a}:H^0(\CP^1,\Omega^{1-\Delta}(*{w_A}))\rightarrow\widehat{\mathcal{K}}_a(d\xi_a)^{1-\Delta}
\]
are linear maps defined by the Laurent expansion. We set $i_{w_A}=\bigoplus_{a\in A}i_{w_a}$.
The following is well-known as the \textit{strong residue theorem}:
\begin{proposition}\label{prop:orth1}
{\rm (1)} The pairing $\pairing{\;\,}{\;}$ is a complete dual pairing between complete topological vector spaces
$\widehat{\mathcal{K}}_A^\Delta$ and $\widehat{\mathcal{K}}_A^{1-\Delta}$.

\vskip 1ex
\noindent
{\rm (2)} $f\in\widehat{\mathcal{K}}_A^\Delta$ and $g\in\widehat{\mathcal{K}}_A^{1-\Delta}$ belongs to
$i_{w_A}H^0(\CP^1,\Omega^\Delta(*{w_A}))$ and $i_{w_A}H^0(\CP^1,\Omega^{1-\Delta}(*{w_A}))$ 
if and only if 
$\langle\,f,i_{w_A}H^0(\CP^1,\Omega^{1-\Delta}(*{w_A}))\rangle=0$ and $\langle
\,i_{w_A}H^0(\CP^1,\Omega^\Delta(*{w_A})),g\,\rangle=0$, respectively.
\end{proposition}

\begin{proposition}\label{prop:expansion1}
For any $\Phi\in\calv_{w_A}^\dagger(M_A),\,v\in V_\Delta$ and $u\in M_A$  
there exists a unique element $\Phi_1(v;u)_{(z)}\in H^0(\CP^1,\Omega^{\Delta}(*w_A))$ such that
\[
i_{w_a}\left(\Phi_1(v;u)_{(z)}\right)
=\begin{cases}
\displaystyle{\sum_{n\in\Z}\Phi\left(\rho_a\left(J_n(v)\right)u\right)\xi_a^{-n-\Delta}}(d\xi_a)^\Delta&\quad (a\neq\infty),\\
\displaystyle{\sum_{n\in\Z}\Phi\left(\rho_\infty\left(\theta(J_{n}(v))\right)u\right)\xi_\infty^{-n-\Delta}(d\xi_\infty)^\Delta}
&\quad(a=\infty).
\end{cases}
\]
\end{proposition}
\begin{proof}
Let us denote the right hand sides of the equality by $g_a(\xi_a)$ for $a\in A$. 
Using Remark \ref{re:residue102} we see that for any $f(z)(dz)^{1-\Delta}\in H^0(\CP^1,\Omega^{1-\Delta}(*w_A))$ 
\[
\begin{split}
&\pairing{\sum_{a\in A}g_a(\xi_a)}{i_{w_A}(f(z)(dz)^{1-\Delta})}\\
&=\sum_{a\in A\setminus\{\infty\}}\sum_{n\in\Z}\,\Res{\xi_a}\Phi(\rho_a(J_n(v))u)\xi_a^{-n-\Delta}f_a(\xi_a)\,d\xi_a\\
&\phantom{\sum_{a\in A\setminus\{\infty\}\Res{\xi_a}\Phi(\rho_a(J_n(v)u)}
\sum_{n\in\Z}}+\sum_{n\in\Z}\,\underset{\xi_\infty=\infty}{\operatorname{Res}}
\Phi\left(\rho_\infty\left(\theta(J_{n}(v))\right)u\right)\xi_\infty^{-n-\Delta}\,f_\infty(\xi_\infty)\,d\xi_\infty\\
&=\sum_{a\in A\setminus\{\infty\}}\!\!\Res{\xi_a}\Phi\left(\rho_a(J(v,\xi_a))uf_a(\xi_a)\right)
+\Res{\xi_\infty}\Phi\left(\rho_\infty(J^\theta(v,\xi_\infty))uf_\infty(\xi_\infty)\right)\\
&=\Phi\left(\rho_A\left(j_{w_A}\left(v\otimes
f(z)(dz)^{1-\Delta}\right)\right)u\right)=0.
\end{split}
\]
We now see $\sum_{a\in A}g_a(\xi_a)\in i_{w_A}H^0(\CP^1,\Omega^{1-\Delta}(*w_A))$ by Proposition \ref{prop:orth1} (2).
\end{proof}

\begin{remark}
(1) The $1$-point current correlation function $\Phi_1(v;u)_{(z)}$ associated to a conformal block $\Phi$  is a meromorphic
$\Delta$-form on $\CP^1$,  and its poles are located at $z=w_a\,(a\in A)$. 

\vskip 1ex
\noindent
(2) In terms of the nonabelian current $J(v,z)$ the proposition is also stated that in a neighborhood of $z=w_a\,(a\in A)$
\begin{equation}\label{eqn:1point-expansion}
\Phi_1(v;u)_{(z)}
=\begin{cases}
\Phi\left(\rho_a\left(J(v,z-w_a)\right)u\right)(dz)^\Delta&\quad (a\neq\infty),\\
\Phi\left(\rho_\infty\left(J^\theta(v,z)\right)u\right)(dz)^\Delta&\quad
(a=\infty).
\end{cases}
\end{equation}
\end{remark}

\begin{definition}
Let $\Phi$ be an element of $\calv_{w_A}^\dagger(M_A)$, and $v\in V_\Delta,\,u\in M_A$. The form $\Phi_1(v;u)_{(z)}\in
H^0(\CP^1,\Omega^{\Delta}(*w_A))$ in the proposition is called  the \textit{$1$-point current correlation function 
associated to the conformal block} $\Phi$.  To indicate that $\Phi_1(v;u)_{(z)}$ is uniquely determined by $\Phi$ 
we often use the notation
\[
\Phi_1(v;u)_{(z)}=\bra{\Phi}J(v,z)\ket{u}(dz)^\Delta.
\]
\end{definition}

\begin{remark}
Let $\{\Phi_m\}_{m=0}^\infty$ be a system of current correlation functions  at $w_A$. Then the initial term $\Phi_0$ is a
conformal block,  i.e., $\Phi_0\in\calv_{w_A}^\dagger(M_A)$ by the proof of Theorem \ref{th:syscorr}, 
and by definition $\Phi_1$ is the $1$-point current correlation function associated to the conformal block $\Phi_0$;
for the property (e) see the proof of Proposition \ref{th:w01}.
\end{remark}

\begin{proposition}\label{prop:1-point}
Let $w_A=(w_a)_{a\in A}\in (\CP^1)^{N+1}$ be a set of distinct points with $w_\infty=\infty$,
and $\{\Phi_m\}_{m=0}^\infty$ be a system of current correlation functions at $w_A$. 
Then for any $v\in V_{\Delta}$ and $u\in M_A$
\begin{multline*}
\Phi_1(v;u)_{(z)}=\sum_{\substack{a\in
A\setminus\{\infty\}}}\sum_{n\geq-\Delta+1}\Phi_0\left(\rho_a\left(J_n(v)\right)u\right)(z-w_a)^{-n-\Delta}(dz)^{\Delta}\\
+\sum_{n=0}^\infty\Phi_0\left(\rho_\infty\left(\theta(J_{-n-\Delta}(v))\right)u\right)z^{n}(dz)^{\Delta},
\end{multline*}
where the sum for the index $n$ is finite, in particular, $\Phi_1(v;u)_{(z)}(dz)^{-\Delta}$ is a rational function of
$z$.
\end{proposition}
\begin{proof}
The singular parts of the Laurent  expansions of both hand sides are same by the property (d).
Then we conclude that both hand sides being multiplied by $(dz)^{-\Delta}$ give the same rational function 
because the difference of them is holomorphic on $\CP^1$ and zero at $z=\infty$.
\end{proof}

\subsection{2-point current correlation functions}
Our next aim is to construct $\Phi_2$ for the system of current correlation functions 
with the initial term $\Phi$.

\begin{definition}
Let $\Phi\in\calv_{w_A}^\dagger(M_A)$ be a conformal block at $w_A$.
The \textit{2-point current correlation function} $\Phi_2(v_1,v_2;u)_{(z_1,z_2)}$ associated to $\Phi$ is defined by 
\[
\begin{split}
\Phi_2(v_1,v_2;u)_{(z_1,z_2)}
&=\sum_{n\geq-\Delta_1+1}\Phi_1(J_{n}(v_1)v_2;u)_{(z_2)}(z_1-z_2)^{-n-\Delta_1}(dz_1)^{\Delta_1}\\
&+\sum_{a\in
A\setminus\{\infty\}}\sum_{n\geq-\Delta_1+1}\Phi_1(v_2;\rho_a(J_{n}(v_1))u)_{(z_2)}(z_1-w_a)^{-n-\Delta_1}(dz_1)^{\Delta_1}\\
&+\sum_{n=0}^\infty\Phi_1(v_2;\rho_\infty(\theta(J_{-n-\Delta_1}(v_1)))u)_{(z_2)}z_1^n(dz_1)^{\Delta_1}
\end{split}
\]
for any $v_1\in V_{\Delta_1}$ and $v_2\in V_{\Delta_2}$.
The 2-point correlation function $\Phi_2(v_1,v_2;u)_{(z_1,z_2)}$ is also written as
\begin{equation}\label{eqn:2point-current}
\begin{split}
&\Phi_2(v_1,v_2;u)_{(z_1,z_2)}=\Phi_1(J^{>0}(v_1,z_1-z_2)v_2;u)_{(z_2)}(dz_1)^{\Delta_1}\\
&\phantom{\Phi_2(v_1,v_2;u)_{(z_1,z_2)}}+\sum_{a\in
A\setminus\{\infty\}}\Phi_1(v_2;\rho_a(J^{>0}(v_1,z_1-w_a))u)_{(z_2)}(dz_1)^{\Delta_1}\\
&\phantom{\Phi_2(v_1,v_2;u)_{(z_1,z_2)}}+\Phi_1(v_2;\rho_\infty(J^{\theta,\,\leq0}(v_1,z_1))u)_{(z_2)}(dz_1)^{\Delta_1},
\end{split}
\end{equation}
where $J^{>0}(v,z)=\sum_{n\geq-\Delta+1}J_n(v)z^{-n-\Delta}$ and $J^{\leq0}(v,z)=\sum_{n\leq-\Delta}J_n(v)z^{-n-\Delta}$.
\end{definition}
\vskip 1.5ex
\begin{remark}
The $2$-point current correlation function
$\Phi_2(v_1,v_2;u)_{(z_1,z_2)}$ is a meromorphic $(\Delta_1,\Delta_2)$-form of $z_1$ and $z_2$, 
and its poles are located at $z_1=z_2,\,z_i=w_a\,(i=1,2, a\in A)$.
\end{remark}

\begin{proposition}\label{pro:expan1}
For $v_1\in V_{\Delta_1}$ and $v_2\in V_{\Delta_2}$ the Laurent expansion of\,
the $2$-point current correlation function $\Phi_2(v_1,v_2;u)_{(z_1,z_2)}$ 
associated to a conformal block $\Phi\in\calv_{w_A}^\dagger(M_A)$ is as follows\,{\rm :}

\vskip 1.5ex
\noindent
{\rm (1)} In a neighborhood of $z_1=z_2$
\[
\Phi_2(v_1,v_2;u)_{(z_1,z_2)}=\sum_{n\in\Z}\Phi_1(J_n(v_1)v_2;u)_{(z_2)}(z_1-z_2)^{-n-\Delta_1}(dz_1)^{\Delta_1}.
\]
\vskip 1ex
\noindent
{\rm (2)} In a neighborhood of $z_1=w_a\,(a\in A)$ 
\[
\Phi_2(v_1,v_2;u)_{(z_1,z_2)}=
\begin{cases}
\displaystyle{\sum_{n\in\Z}\Phi_1(v_2;\rho_a(J_n(v_1))u)_{(z_2)}(z_1-w_a)^{-n-\Delta_1}(dz_1)^{\Delta_1}}&\quad(a\neq\infty),\\
\displaystyle{\sum_{n\in\Z}\Phi_1\left(v_2;\rho_\infty
\left(\theta(J_{n}(v_1))\right)u\right)_{(z_2)}\,z_1^{-n-\Delta_1}\,(dz_1)^{\Delta_1}}&\quad(a=\infty).
\end{cases}
\]
\end{proposition}

\begin{remark}\label{re:opeforblock}
The statements (1) and (2) in Proposition \ref{pro:expan1} respectively means:

\vskip 1.5ex
\noindent
(1) In a neighborhood of $z_1=z_2$ 
\[
\Phi_2(v_1,v_2;u)_{(z_1,z_2)}=\Phi_1(J(v_1,z_1-z_2)v_2;u)_{(z_2)}(dz_1)^{\Delta_1}.
\]

\vskip 1ex
\noindent
(2) In a neighborhood of $z_1=w_a\, (a\in A)$ 
\[
\Phi_2(v_1,v_2;u)_{(z_1,z_2)}=
\begin{cases}
\Phi_1(v_2;\rho_a(J(v_1,z_1-w_a))u)_{(z_2)}(dz_1)^{\Delta_1}&\quad(a\neq\infty),\\
&\\
\Phi_1(v_2;\rho_\infty(J^\theta(v_1,z_1))u)_{(z_2)}(dz_1)^{\Delta_1}&\quad(a=\infty).
\end{cases}
\]
\end{remark}

We now give a proof of the statement (1) of the proposition. By (\ref{eqn:2point-current})  the first statement is equivalent to 
\begin{multline}\label{eqn:2-point01}
\Phi_1(J^{\leq0}(v_1,z_1-z_2)v_2;u)_{(z_2)}\\
=\Phi_1(v_2;\rho_\infty(J^{\theta,\,\leq0}(v_1,z_1))u)_{(z_2)}
+\sum_{a\in A\setminus\{\infty\}}\Phi_1(v_2;\rho_a(J^{>0}(v_1,z_1-w_a))u)_{(z_2)}.
\end{multline}
For any fixed $z_1$ both hand sides of (\ref{eqn:2-point01}) are meromorphic functions of $z_2$, whose poles
are located at $z_2=w_a\,(a\in A)$. 

In order to prove (\ref{eqn:2-point01}) we show that both hand sides
of (\ref{eqn:2-point01}) have the same singular part at the pole $z_2=w_a\,(a\in A)$.  
Using (\ref{eqn:1point-expansion}) and (\ref{eqn:asso01}) we see that the singular part of (\ref{eqn:2-point01}) at
$z_2=w_a\,(a\neq\infty)$ is
\begin{multline}\label{eqn:singular102}
\Phi\left(\rho_a\left(J^{\leq0}(v_1,z_1-w_a)\right)u_a\right)\\
=\Phi\left(\rho_\infty\left(J^{\theta,\,\leq0}(v_1,z_1)\right)u_a\right)
+\sum_{b\in A\setminus\{a,\infty\}}\Phi\left(\rho_b\left(J^{>0}(v_1,z_1-w_b)\right)u_a\right),
\end{multline}
where $u_a=\rho_a\left(J^{>0}(v_2,z_2-w_a)\right)u$, while using (\ref{eqn:1point-expansion}) and Remark \ref{re:theta-asso} we
see that the nonnegative power term of the Laurent expansion at $z_2=\infty$ is
\begin{equation}\label{eqn:singular103}
\Phi\left(\rho_\infty\left(J^{\theta,>0}(v_1,z_1)\right)u_\infty\right)
=\sum_{b\in A\setminus\{\infty\}}\Phi\left(\rho_b\left(J^{>0}(v_1,z_1-w_b)\right)u_\infty\right),
\end{equation}
where $u_\infty=\rho_\infty\left(J^{\theta,\leq0}(v_2,z_2)\right)u$. 
Now suppose that (\ref{eqn:singular102}) and (\ref{eqn:singular103}) hold. 
Then the both hand sides of (\ref{eqn:2-point01}) have the same singular part at the pole $z_2=w_a\,(a\in A)$, and the difference of
them is zero at $z_2=\infty$, which shows (\ref{eqn:2-point01}). So it suffice to see:

\begin{lemma}\label{lem:inv}
Let $\Phi\in\calv_{w_A}^\dagger(M_A)$ and $u\in M_A$. Then 
\begin{align*}
&\Phi\left(\rho_\infty\left(J^{\theta,>0}(v,z)\right)u\right)
=\sum_{b\in A\setminus\{\infty\}}\Phi\left(\rho_b\left(J^{>0}(v,z-w_b)\right)u\right),\mbox{ and}\\
&\Phi\left(\rho_a\left(J^{\leq0}(v,z-w_a)\right)u\right)\\
&\phantom{\Phi(\rho_a(J^{\leq0}}=\Phi\left(\rho_\infty\left(J^{\theta,\,\leq0}(v,z)\right)u\right)
+\sum_{b\in A\setminus\{a,\infty\}}\Phi\left(\rho_b\left(J^{>0}(v,z-w_b)\right)u\right)\mbox{for $a\neq\infty$.}
\end{align*}
\end{lemma}
\begin{proof}
Recall the definition of the conformal block, i.e.,
\begin{equation}\label{eqn:block-identity}
\Phi(j_a(v\otimes\varphi)u)+\sum_{b\neq a,\infty}\Phi(j_b(v\otimes\varphi)u)+\Phi(j_\infty(v\otimes\varphi)u)=0
\end{equation}
for all $v\in V_\Delta$ and $\varphi\in H^0(\CP^1,\Omega^{1-\Delta}(*w_A))$. 
Then the first and second equality follows from (\ref{eqn:block-identity}) for $\varphi=z^{n+1}(dz)^{1-\Delta}\,(n\geq0)$ and
$\varphi=(z-w_a)^{-n-1}(dz)^{1-\Delta}\,(n\geq0)$, respectively.
\end{proof}

The statement (2) of Proposition \ref{pro:expan1} can be proved in a very similar way by Lemma \ref{lem:inv}; here we use 
\[
J(J^{>0}(v_1,z_1-z_2)v_2,z_2)= 
\begin{cases}
[J(v_2,z_2),J^{\leq0}(v_1,z_1)]&\quad\mbox{for}\quad 0<|z_1|<|z_2|,\\
{}[J^{>0}(v_1,z_1),J(v_2,z_2)]&\quad\mbox{for}\quad 0<|z_2|<|z_1|,
\end{cases}
\]
and 
\[
J^{\theta}(J^{>0}(v_1,z_1-z_2)v_2,z_2)= 
\begin{cases}
[J^{\theta,\leq0}(v_1,z_1),J^{\theta}(v_2,z_2)]&\quad\mbox{for}\quad 0<|z_1|<|z_2|,\\
{}[J^{\theta}(v_2,z_2),J^{\theta,>0}(v_1,z_1)]&\quad\mbox{for}\quad 0<|z_2|<|z_1|.
\end{cases}
\]

\begin{remark}
Let $A$ be a finite set. 
In order to avoid complicated description of 1-point and 2-point functions we have so far considered sets of
points $w_A\subset\CP^1$ which contain the point $\infty$. 
In the case that $w_a\neq\infty$ for all $a\in A$ we just drop the terms
$z^n$  from its definition.  More precisely the $1$-point and $2$-point current correlation function associated to
$\Phi$ is respectively given by $\Phi_1(\ket{0};u)_{(z)}=\Phi(u),
\Phi_2(\ket{0},v_2;u)_{(z_1,z_2)}=\Phi_1(v_2;u)_{(z_2)}$, and for any $v\in V_{\Delta}\,(\Delta\neq0),\,v_1\in
V_{\Delta_1}\,(\Delta_1\neq0)$
\begin{align}
\Phi_1(v;u)_{(z)}&
=\sum_{a\in A}\sum_{n\geq-\Delta+1}\Phi(\rho_a(J_n(v))u)(z-w_a)^{-n-\Delta}(dz)^{\Delta}\quad \mbox{and}\notag\\
\Phi_2(v_1,v_2;u)_{(z_1,z_2)}&
=\sum_{n\geq-\Delta_1+1}\Phi_1(J_{n}(v_1)v_2;u)_{(z_2)}(z_1-z_2)^{-n-\Delta_1}(dz_1)^{\Delta_1}\notag\\
&+\sum_{a\in A}\sum_{n\geq-\Delta_1+1}\Phi_1(v_2;\rho_a(J_{n}(v_1))u)_{(z_2)}(z_1-w_a)^{-n-\Delta_1}(dz_1)^{\Delta_1}.
\label{eqn:2point-noinfty}
\end{align}
\end{remark}

\subsection{Propagation of vacua}\label{sub:propvacua}
Let $\tilde{A}=\{0\}\sqcup A$ and $w_{\tilde{A}}=(w_0,w_A),\,w_A=(w_a)_{a\in A}$ where $w_0\neq w_a$ for all $a\in A$,
and set $M_{\tilde{A}}=V\otimes M_A$. 
The main aim of this subsection is to prove the following theorem which is called \textit{the propagation of vacua}\,.

\begin{theorem}
\label{th:propvacua}
The linear map $M_A\rightarrow M_{\tilde{A}}\,(u\mapsto\ket{0}\otimes u)$ induces a vector space isomorphism
$\calv_{w_A}(M_A)\cong\calv_{w_{\tilde{A}}}(M_{\tilde{A}})$.
\end{theorem}
\begin{proof}
We first show that the linear map $M_A\rightarrow M_{\tilde{A}}\,(u\mapsto\ket{0}\otimes u)$ induce a linear map
$\iota:\calv_{w_A}(M_A)\rightarrow\calv_{w_\tilde{A}}(M_{\tilde{A}})$. To see this it suffice to prove that $\ket{0}\otimes
u\in\frakg_{w_{\tilde{A}}}^{out}(V)M_{\tilde{A}}$ for any $u\in\frakg_{w_A}^{out}(V)M_A$; using the fact that $j_{w_0}(v\otimes
f)\ket{0}=0$ if $f(z)(dz)^{1-\Delta}$ is holomorphic at $z=w_0$ 
we see that 
$\ket{0}\otimes j_{w_A}(v\otimes f)u=j_{w_{\tilde{A}}}(v\otimes f)(\ket{0}\otimes u)\in\frakg_{w_{\tilde{A}}}^{out}(V)M_{\tilde{A}}$
for any $v\in V_\Delta$  and $f\in H^0(\CP^1,\Omega^{1-\Delta}(*w_A))$.

In order to prove that the induced map $\iota$ is surjective we will show that for any $v\in V_\Delta$ and $u\in M_A$ there exists 
$f\in H^0(\CP^1,\Omega^{1-\Delta}(*w_{\tilde{A}}))$ such that $v\otimes u\equiv -\ket{0}\otimes j_{w_A}(v\otimes f)u$.
Take $f$ such that $j_{w_0}(f)=\left(\xi_0^{-1}+\mbox{higher order terms of $\xi_0$}\right)(d\xi_0)^{1-\Delta}.$
Then
\[
J_{-\Delta}(v)\ket{0}\otimes u=j_{w_0}(v\otimes f)\ket{0}\otimes u
=j_{w_{\tilde{A}}}(v\otimes f)(\ket{0}\otimes u)-\ket{0}\otimes j_{w_A}(v\otimes f)u,
\]
in other words, 
$v\otimes u=J_{-\Delta}(v)\ket{0}\otimes u\equiv -\ket{0}\otimes j_{w_A}(v\otimes
f)u\mod\frakg_{w_{\tilde{A}}}^{out}(V)M_{\tilde{A}}$.

We now prove that the map $\iota$ is injective, 
which is equivalent to showing that the dual map
$\iota^*:\calv_{w_{\tilde{A}}}^\dagger(M_{\tilde{A}})\rightarrow\calv_{w_A}^\dagger(M_A)$ is surjective, i.e., for any
$\Phi\in\calv_{w_A}^\dagger(M_A)$ there exists $\tilde{\Phi}\in\calv_{w_{\tilde{A}}}^\dagger(M_{\tilde{A}})$ such that
$\tilde{\Phi}(\vac\otimes u)=\Phi(u)$ for all $u\in M_A$. We define $\tilde{\Phi}\in\Hom_\C(M_{\tilde{A}},\C)$ by
\[
\tilde{\Phi}(u_0\otimes
u)=\Phi_1(u_0;u)_{(w_0)}(dw_0)^{-\Delta}=\bra{\Phi}J(u_0,w_0)\ket{u}
\]
for any $u_0\in V_\Delta$ and $u\in M_A$,where
$\Phi_1(u_0;u)_{(z)}$ is the $1$-point current correlation function associated to the conformal block $\Phi$.
It suffices to show that $\tilde{\Phi}$ induces an element in $\calv_{w_{\tilde{A}}}^\dagger(M_{\tilde{A}})$, i.e., 
$\tilde{\Phi}(j_{w_0}(v\otimes f)u_0\otimes u)+\tilde{\Phi}(u_0\otimes j_{w_A}(v\otimes f)u)=0$
for all $v\in V_\Delta$ and $f(z)(dz)^{1-\Delta}\in H^0(\CP^1,\Omega^{1-\Delta}(*w_{\tilde{A}}))$, which follows from 
Proposition \ref{pro:expan1} and Remark \ref{re:opeforblock} as
\[
\begin{split}
\bra{\Phi}(J(j_0(v\otimes f)u_0,w_0)\ket{u}&=\underset{z=w_0}{\operatorname{Res}}\,\bra{\Phi}J(J(v,z-w_0)u_0,w_0)\ket{u}f(z)\,dz\\
&=\underset{z=w_0}{\operatorname{Res}}\,\bra{\Phi}J(v,z)J(u_0,w_0)\ket{u}f(z)\,dz\\ 
&=-\sum_{a\in A}\underset{z=w_a}{\operatorname{Res}}\,\bra{\Phi}J(v,z)J(u_0,w_0)\ket{u}f(z)\,dz\\ 
&=-\sum_{a\in A}\underset{z=w_a}{\operatorname{Res}}\,\bra{\Phi}J(u_0,w_0)\ket{J(v,z-w_a)u}f(z)\,dz\\
&=-\sum_{a\in A}\bra{\Phi}J(u_0,w_0)\ket{j_a(v\otimes f)u}.
\end{split}
\]
\end{proof}

As a corollary we get the $\mathcal{S}_2$-symmetry of the $2$-point current correlation functions:

\begin{corollary}\label{cor:sym}
\textit{Let $v_i\in V_{\Delta_i}\,(i=1,2)$. Then $\Phi_2(v_1,v_2;u)_{(z_1,z_2)}=\Phi_2(v_2,v_1;u)_{(z_2,z_1)}.$
}
\end{corollary}

The $1$-point current correlation function is a meromorphic $\Delta$-form on $\CP^1$, whose poles are located only 
at $z=w_a\,(a\in A)$. If $\infty\in A$ the $2$-point function has poles only at $z_1=z_2$ and $z_i=w_a\,(i=1,2,a\in A)$, 
while if $\infty\not\in A$ the point $z_1=\infty$ may be also a pole by (\ref{eqn:2point-noinfty}).
However, since the 2-point current correlation function is holomorphic at $z_2=\infty$ 
it is also holomorphic at $z_1=\infty$ by the $\mathcal{S}_2$-symmetry:

\begin{proposition}\label{prop:divisors02}
The $2$-point current correlation function $\Phi_2(v_1,v_2;u)$ associated to a conformal block
$\Phi\in\calv_{w_A}^\dagger(M_A)$
 belongs to $H^0((\CP^1)^2,\Omega^{\Delta_1\boxtimes\Delta_2}(*D_{1,\,2}+\sum_{i,\,a}*D_{i,\,a})$.
\end{proposition}

\subsection{Construction of systems of current correlation functions}\label{sub:cor-con}
Using the propagation of vacua we will obtain the whole system of current correlation function $\{\Phi_m\}_{m=0}^\infty$ 
with a given initial term $\Phi_0=\Phi\in\calv_{w_A}(M_A)$.

Let $m$ be a nonnegative integer.  For any set of \textit{distinct} $m$ points $z_{[m]}=(z_i)_{i\in[m]}\in\C^m$ such that $z_i\neq
w_a\,(i\in[m],a\in A)$ we assign the space of covacua $\calv_{w_A,[m]}(M_{[m],A})$.
By the propagation of vacua there exists a canonical vector space isomorphism
$\iota_m:\calv_{z_{[m]},w_A}(M_{[m],A})\rightarrow\calv_{z_{[m+1]},w_A}(M_{[m+1],A})$.
Let $\bar{\Phi}_0=\Phi\in\calv_{w_A}(M_A)$. We define a system $\{\bar{\Phi}_m\}_{m=0}^\infty$ by
$\bar{\Phi}_{m+1}=\iota_m^{*-1}(\bar{\Phi}_m)$ for $m\geq 0$.

\begin{proposition}\label{th:w01}
Let $u\in M_A,\,v_i\in V_{\Delta_i}\,(i\in [m])$ and $\Phi\in\calv_{w_A}^\dagger(M_A)$. Then 
\[
\Phi_m=\bar{\Phi}_m(v_1,\ldots, v_m;u)_{(z_1,\ldots,z_m)}(dz_1)^{\Delta_1}\cdots(dz_m)^{\Delta_m},\quad m=0,1,\ldots
\]
is the system of current correlation functions at $w_A$ with the initial term $\Phi$.
\end{proposition}
\begin{proof}
The properties (a) and (d) are immediate consequences of the inductive definition of $\{\Phi_m\}_{m=0}^\infty$. The property (b)
follows from Proposition \ref{prop:divisors02}. 

We now prove the $\mathcal{S}_m$-symmetry. On the one hand, using Corollary \ref{cor:sym} and the inductive definition of $\Phi_m$ we
see that
\[
\Phi_m(v_1,v_2,v_3,\ldots,v_m;u)_{(z_1,z_2,z_3,\ldots,z_m)}=\Phi_m(v_2,v_1,v_3,\ldots,v_m;u)_{(z_2,z_1,z_3,\ldots,z_m)}.
\]
On the other hand, at a neighborhood of $z_1=z_j\,(j\neq1)$ and $z_1=w_a\,(a\in A)$ we have the Laurent expansions
\[
\begin{split}
&\Phi_m(v_1,v_2,v_3,\ldots,v_m;u)_{(z_1,z_2,z_3,\ldots,z_m)}\\
&\phantom{\langle\Phi|J(v_1}=\sum_{n\in\Z}\Phi_{m-1}(v_2,\ldots,J_n(v_1)v_j,\ldots,v_m;u)_{(z_2,\ldots,z_m)}(z_1-z_j)^{-n-\Delta_1}
(dz_1)^{\Delta_1},\quad\mbox{and}\\
&\Phi_m(v_1,v_2,v_3,\ldots,v_m;u)_{(z_1,z_2,z_3,\ldots,z_m)}\\ 
&\phantom{\langle\Phi|J(v_1}=\sum_{n\in\Z}
\Phi_{m-1}(v_2,\ldots,v_m;\rho_a(J_{n}(v_1))u)_{(z_2,\ldots,z_m)}(z_1-w_a)^{-n-\Delta_1}(dz_1)^{\Delta_1}.
\end{split}
\]
The $\mathcal{S}_m$-symmetry now follows by induction on $m$.

To verify (e) it suffices to show the case $m=1$, i.e.,
\[
\partial/\partial z\,\Phi_1(v;u)_{(z)}\,dz=\Phi_1(T(-1)v;u)_{(z)}.
\]
Since both hand sides are meromorphic $(\Delta+1)$-forms on
$\CP^1$ with poles only at $z=w_a\,(a\in A)$ it is enough to see that these two have the same Laurent expansions at the poles;
in a neighborhood of
$z=w_a\,(a\in A)$
\[
\begin{split}
\Phi_1(T(-1)v;u)_{(z)}&=\sum_{n\in\Z}\Phi(\rho_a(J_n(T(-1)v))u)(z-w_a)^{-n-\Delta-1}(dz)^{\Delta+1}\\
&=-\sum_{n\in\Z}(n+\Delta)\Phi(\rho_a(J_n(v))u)(z-w_a)^{-n-\Delta-1}(dz)^{\Delta+1}\\
&=\frac{\partial}{\partial z}\Phi_1(v;u)_{(z)}\,dz.
\end{split}
\] 
\end{proof}

\subsection{Finiteness of conformal blocks}
\begin{theorem}\label{th:finiteblock}
Let $V$ be a $C_2$-finite chiral vertex operator algebra. 
Let $A$ be the finite set $A=\{\,1,\,2,\,\ldots,\,N,\,\infty\,\}$ and $w_A=(w_a)_{a\in A}$ be a set of distinct points in
$\CP^1$ with $w_\infty=\infty$. Let $M_A=\bigotimes_{a\in A}M^a$ where $M^a$ is a $V$-module. 
Then the space of covacua $\calv_{w_A}(M_A)$ is finite dimensional.
\end{theorem}
\begin{proof}
Let $M^a=\bigoplus_{i=1}^r\bigoplus_{n=0}^\infty M_{(h_i+n)}^a$ for all $a\in A$.
We set 
\[
F_pM^a=
\begin{cases}
0&\quad (p\leq-1),\\
\bigoplus_{i=1}^r\bigoplus_{n=0}^p M_{(h_i+n)}^a&\quad (p\geq0),
\end{cases}
\]
which define an increasing filtration on $M^a$ and the induced filtration on $C_2(M^a)$.
We note that $\gr_\bullet^F M^a\cong M^a$ and $\gr_\bullet^F C_2(M^a)\cong C_2(M^a)$.
Since $V$ is $C_2$-finite we see that $\gr_\bullet^F M^a/\gr_\bullet^F C_2(M^a)$ is finite dimensional
by Corollary \ref{cor:finiteproperty}.

Let us define an increasing filtration on $M_A$ by 
\[
F_{p}M_{A}=\sum_{\sum p_a=p}\bigotimes_{a\in A}F_{p_a}M^{a}.
\]
The filtration $F_pM_A\,(p\in\Z)$ on $M_A$ induces a filtration on $\frakg_{w_A}^{out}(V)M_A$,
and then a filtration on $\calv_{w_A}(M_A)$.
There exists a canonical exact sequence as vector spaces
\[
0\rightarrow\gr_\bullet^F(\frakg_{w_A}^{out}(V)M_A)\rightarrow\gr_\bullet^F M_A\rightarrow\gr_\bullet^F\calv_{w_A}(M_A)\rightarrow0\,,
\]
in other words, $\gr_\bullet^F\calv_{w_A}(M_A)\cong\gr_\bullet^F M_A/\gr_\bullet^F(\frakg_{w_A}^{out}(V)M_A)$.
Let us denote by $\pi$ the canonical projection 
$\gr_\bullet^F M_A=\bigotimes_{a\in A}\gr_\bullet^F M^a\rightarrow\gr_\bullet^F\calv_{w_A}(M_A)$. We will show that $\pi$ induces a
surjective linear map $\bigotimes_{a\in A}\gr_\bullet^F M^a/\gr_\bullet^F C_2(M^a)\rightarrow\gr_\bullet^F\calv_{w_A}(M_A)$ so that
$\gr_\bullet^F\calv_{w_A}(M_A)$ is finite dimensional. To this end it suffices to see that 
\begin{equation}\label{eqn:red}
\gr_\bullet^F M^1\otimes\cdots\otimes\gr_\bullet^F C_2(M^a)\otimes\cdots\otimes\gr_\bullet^F
M^\infty\subseteq \operatorname{ker}\,\pi\quad\mbox{for all $a\in A$.}
\end{equation}

For any $a\in A$ we set
\[
\varphi^a=
\begin{cases}
(z-w_{a})^{-n+\Delta-1}(dz)^{1-\Delta}&\quad\mbox{if $a\neq\infty$},\\
z^{n+\Delta-1}(dz)^{1-\Delta}&\quad \mbox{if $a=\infty$}.
\end{cases}
\]
If $n+\Delta-1\geq0$ then $\varphi^a\in H^0(\CP^1,\Omega^{1-\Delta}(*w_A))$ for all $a\in A$.
For any $v\in V_\Delta$ and $m_A=\otimes_{b\in A}m_b\in M_A$ we see that
\begin{multline}\label{eqn:finiteblock}
\rho_A(j_A(v\otimes\varphi^a))(m_A)
=m_{1}\otimes\cdots\otimes J_{-n}(v_a)m_{a}\otimes\cdots\otimes m_{\infty}\\
+\sum_{\substack{b\in A\\b\neq a}}m_{1}\otimes\cdots\otimes
j_{b}(v\otimes\varphi^a)m_{b}\otimes\cdots\otimes m_{\infty}\,,
\end{multline}
where $v_a=v$ and $v_a=e^{T(1)}v$ for $a\neq\infty$ and $a=\infty$, respectively.
We suppose that $n\geq \Delta$ and $m_A\in F_{p-n}M_A$. Then the second term of the right hand side of (\ref{eqn:finiteblock})
belongs to $F_{p-1}M_A$ so that
that 
\[
m_{1}\otimes\cdots\otimes J_{-n}(v_a)m_{a}\otimes\cdots\otimes m_{\infty}\in F_{p-1}M_{A}+F_p\frakg_{w_A}^{out}(V)M_A.
\]
So (\ref{eqn:red}) immediately follows for $a\neq\infty$; the case $a=\infty$ is proved by using the identity $v=e^{T(1)}e^{-T(1)}v$.
\end{proof}

\begin{remark} 
From the proof of the theorem we see that $\dim_\C\calv_{w_A}(M_A)<\infty $ if $\dim_\C M^a/C_1(M^a)<\infty$ for all $a\in A$ where
$C_1(M)$ is the vector subspace which is linearly spanned by vectors $J_{-n}(v)m\,(v\in V_\Delta,\,n\geq\Delta)$ for all $m\in M$.
However in the case of the space of conformal blocks over a higher genus compact Riemann surface even $C_2$-finite condition is 
not enough for the finiteness; one has to suppose $C_n$-finiteness for all $n\geq2$, which is proved in Corollary
\ref{cor:finiteproperty}.   
\end{remark}

\subsection{$2$-point conformal blocks}
\begin{proposition}\label{prop:homblock}
Let $A=\{0,\,\infty\}$, and $w_A=(0,\infty)$.  Let $V$ be a $C_2$-finite chiral vertex operator algebra and $M^0,\,M^\infty$ be
$V$-modules.  Then \[
\calv_{w_A}^\dagger(M_A)\cong\Hom_{\,\calu(V)}(M^0,D(M^\infty))
\]
where $M_A=M^0\otimes M^\infty$.
\end{proposition}
\begin{proof}
We first note that $\varphi_n=z^{n+\Delta-1}\,(dz)^{1-\Delta}\,(n\in\Z)$ is a topological basis of $H^0(\CP^1,\Omega^{1-\Delta}(*[0]+*[\infty]))$, 
and that $j_0\left(v\otimes\varphi_n\right)u_0=J_n(v)u_0$ and $j_\infty\left(v\otimes\varphi_n\right)u_\infty=-\theta\left(J_n(v)\right)u_\infty$
for all $u_0\in M^0$ and $u_\infty\in M^\infty$. We identify $\calv_{w_A}^\dagger(M_A)$ with the subspace of $\Hom_\C(M_A,\C)$, which
consists of elements $\Phi\in\Hom_\C(M_A,\C)$ such that $\Phi(\frakg_A^{out}(V)M_A)=0$, i.e.,
\begin{equation}
\Phi(u_0\otimes\theta(J_n(v))u_\infty)=\Phi(J_n(v)u_0\otimes u_\infty).
\label{eqn:dual}
\end{equation}

Let $f\in\Hom_{\,\calu(V)}(M^0,D(M^\infty))$ and define $\Phi_f\in\Hom_\C(M_A,\C)$ by $\Phi_f(u_0\otimes
u_\infty)=\pairing{\,f(u_0)}{u_\infty}$. Then we wee that 
\begin{multline*}
\Phi_f\left(u_0\otimes\theta\left(J_n(v)\right)u_\infty\right)=\pairing{\,f(u_0)}{\theta\left(J_n(v)\right)u_\infty}\\
=\pairing{\,J_n(v)f(u_0)}{u_\infty}=\pairing{\,f(J_n(v)u_0)}{u_\infty}
=\Phi_f\left(J_n(v)u_0\otimes u_\infty\right),
\end{multline*}
i.e., $\Phi_f\in\calv_{w_A}(M_A)^\dagger$. Conversely for any $\Phi\in\calv_{w_A}^\dagger(M_A)$ and $u_0\in M^0$ we define
$f(u_0)\in\Hom_\C(M^\infty,\C)$ by $\pairing{\,f(u_0)}{u_\infty}=\Phi(u_0\otimes u_\infty)$. Using (\ref{eqn:dual}) we see that
$\pairing{\,f(T(0)u_0)}{u_\infty}=\pairing{\,f(u_0)}{T(0)u_\infty}$ so that if $u_0\in M^0_{(h_0)}$
then $f(u_0)|M^0_{(h)}=0$ for all $h\neq h_0$ and $f(u_0)\in D(M^\infty)$. We now conclude that
$f\in\Hom_{\,\calu(V)}(M^0,D(M^\infty))$ by (\ref{eqn:dual}).
\end{proof}

\section{Integrable connections}\label{sec:connection}
\subsection{Integrable connection on sheaves of covacua}
\label{sub:covacua}
Let $A=\{\,1,\,\ldots,\,N,\,\infty\,\}$ and
$X_A=\{w_A=(w_a)_{a\in A}\}$ be the configuration space of distinct points on $\CP^1$ such that $w_\infty=\infty$.
Since we fix the finite set $A$ throughout this subsection we simply write $X=X_A$.

In the previous section we study the space of covacua associated to a point $w_A\in X$ in the configuration space. Here we vary the point $w_A$, and
consequently we get the sheaf of covacua; the $w_A$ dependence of the space of covacua will be described by an integrable connection on the sheaf.

Let $\pi:\calc=X\times\CP^1\rightarrow X$ be the trivial $\CP^1$-bundle over $X$, and $z$ be the canonical affine coordinate  of $\CP^1$;
the pair $(w_A,z)$ is a local coordinate for $\calc$.  We define cross-sections $s_a:X\rightarrow\calc\,(a\in A)$ by
$s_a(w)=(w_A,w_a)$, and introduce divisors on $\calc$ by $S_a=s_a(X)$ and $S_A=\bigcup_{a\in A}S_a$.
For $a\neq\infty$ the holomorphic function $\xi_a=z-w_a$ defines a local coordinate system $(w_A,\xi_a)$ of $\calc$ around the divisor $S_a$,
while the meromorphic function  $\xi_\infty$ defines a local meromorphic coordinate system $(w_A,\xi_\infty)$ of $\calc$ around the divisor
$S_\infty$ 

Let $\calo_X$ be the structure sheaf of $X$, i.e., the sheaf of germs of holomorphic functions on $X$. 
For any $a\in A$ we define the Lie algebra $\frakg_{X,a}(V)$ over $\calo_X$ by
\[
\frakg_{X,a}(V)=\bigoplus_{\Delta=0}^\infty V_\Delta\otimes\calo_X((\xi_a))(d\xi_a)^{1-\Delta}
/\nabla(\bigoplus_{\Delta=0}^\infty V_\Delta\otimes\calo_X((\xi_a))(d\xi_a)^{-\Delta}),
\] 
and set $\frakg_X(V)=\oplus_{a\in A}\frakg_{X,a}(V)$, which is also a Lie algebra over $\calo_X$.
Let $M^a\,(a\in A)$ be $V$-modules and set $M_X=\calo_X\otimes M_A$ where $M_A=\bigotimes_{a\in A}M^a$.
Then $M_X$ becomes a $\frakg_X(V)$-module over $\calo_X$. 

Let $\Omega_{\calc/X}^1$ be the sheaf of relative $1$-forms for the surjective holomorphic map $\pi:\calc\rightarrow X$, and
$\Omega^{\Delta}_{\calc/X}(*S_A)$ be the sheaf of relative $\Delta$-forms which have possible poles on the divisor $S_A$. For any element in
$V_\Delta\otimes\pi_*\Omega^{1-\Delta}_{\calc/X}(*S_A)$, where $\pi_*\Omega^{1-\Delta}_{\calc/X}(*S_A)$ is the direct image sheaf,
we assign an element in $\frakg_{X,a}(V)$ by taking Laurent expansion at $S_a$ in terms of coordinates $\xi_a$, i.e., for any 
$v\otimes f(w,z)(dz)^{1-\Delta}\in V_\Delta\otimes\pi_*\Omega^{1-\Delta}_{\calc/X}(*S_A)$ we define the map
$j_a:V_\Delta\otimes\pi_*\Omega^{1-\Delta}_{\calc/X}(*S_A)\rightarrow\frakg_{X,a}(V)$ by
\[
j_a\left(v\otimes f(w,z)(dz)^{1-\Delta}\right)= J\left(v,f_a(w,\xi_a)\right),
\]
where $f_a(w,\xi_a)$ is the Laurent expansion of $f(w,z)$ at $z=w_a$ in terms of the coordinate $\xi_a$.
We set $j_A=\sum_{a\in A}j_a$.

\begin{definition}
We define the $\calo_X$-submodule $\frakg_X^{out}(V)$ of $\frakg_X(V)$ by
\[
\frakg_X^{out}(V)=j_A\left(\bigoplus_{\Delta=0}^\infty V_\Delta\otimes\pi_*\Omega^{1-\Delta}_{\calc/X}(*S_A)\right).
\]
\end{definition}

\begin{proposition}
\label{prop:sub}
The $\calo_X$-submodule $\frakg_X^{out}(V)$ is a Lie subalgebra of $\frakg_X(V)$ over $\calo_X$.
\end{proposition}
\begin{proof}
Let $v_1\in V_{\Delta_1}$ and $v_2\in V_{\Delta_2}$, and let 
$f_1(w,z)(dz)^{1-\Delta_1}\in\pi_*\Omega^{1-\Delta_1}_{\calc/X}(*S_A)$ and $f_2(w,z)(dz)^{1-\Delta_2}\in\pi_*\Omega^{1-\Delta_2}_{\calc/X}(*S_A)$.
Then for any $a\in A$ we see that
\[
\begin{split}
&\left[j_a\left(v_1\otimes f_{1,a}(w,\xi_a)(dz)^{1-\Delta_1}\right),j_a\left(v_2\otimes f_{2,a}(w,\xi_a)(dz)^{1-\Delta_2}\right)\right]\\
&=\left[J\left(v_1,f_{1,a}(w,\xi_a)\right),J\left(v_2,f_{2,a}(w,\xi_a)\right)\right]\\
&=\sum_{m=0}^{\Delta_1+\Delta_2-1}\frac{1}{m!}
J\left(J_{m-\Delta_1+1}(v_1)v_2,\left(\frac{d}{d\xi_a}\right)^m\!\!\!f_{1,a}(w,\xi_a)\cdot f_{2,a}(w,\xi_a)\right)\\
&=\sum_{m=0}^{\Delta_1+\Delta_2-1}j_a\left(\frac{1}{m!}
J_{m-\Delta_1+1}(v_1)v_2\otimes
\left(\frac{d}{dz}\right)^m\!\!\!f_1(w,z)\cdot f_2(w,z)(dz)^{2-\Delta_1-\Delta_2+m}\right).
\end{split}
\]
\end{proof}

Let $\bar{A}=A\setminus\{\infty\}$. For any $a\in\bar{A}$ we define the
derivation $\nabla_{\partial/\partial w_a}^{\frakg_X(V)}:\frakg_X(V)\rightarrow\frakg_X(V)$ by
$\nabla_{\partial/\partial w_a}^{\frakg_X(V)}=\partial/\partial w_a+[\rho_a(T(-1)),\,*]$, and the linear map
$\nabla_{\partial/\partial w_a}^{M_X}:M_X\rightarrow M_X$ by $\nabla_{\partial/\partial w_a}^{M_X}=\partial/\partial
w_a+\rho_a(T(-1))$.

\begin{proposition}\label{prop:connection}
The linear maps $\nabla_{\partial/\partial w_a}^{\frakg_X(V)}$ and $\nabla_{\partial/\partial w_a}^{M_X}\,(a\in\bar{A}\,)$ satisfy the
following.

\vskip 1.5ex
\noindent
{\rm (1)} $[\nabla_{\partial/\partial w_a}^{\frakg_X(V)},\nabla_{\partial/\partial w_b}^{\frakg_X(V)}]=0$ and
$[\nabla_{\partial/\partial w_a}^{M_X},\nabla_{\partial/\partial w_b}^{M_X}]=0$.

\vskip 1ex
\noindent
{\rm (2)} $\nabla_{\partial/\partial w_a}^{\frakg_X(V)}(f\cdot J)=\partial_af\cdot J+ f\nabla_{\partial/\partial w_a}^{\frakg_X(V)}(J)$
and $\nabla_{\partial/\partial w_a}^{M_X}(f\cdot m)=\partial_af\cdot m+ f\nabla_{\partial/\partial w_a}^{M_X}(m)$ for all
$f\in\calo_X,\,J\in\frakg_X(V)$ and $m\in M_X$ where $\partial_a=\partial/\partial w_a\,(a\in\bar{A})$.

\vskip 1ex
\noindent
{\rm (3)} $\nabla_{\partial/\partial w_a}^{\frakg_X(V)}([J_1,J_2])=[\nabla_{\partial/\partial
w_a}^{\frakg_X(V)}(J_1),J_2]+[J_1,\nabla_{\partial/\partial w_a}^{\frakg_X(V)}(J_2)]$ for all
$J_1,\,J_2\in\frakg_X(V)$ and $\nabla_{\partial/\partial w_a}^{M_X}(J\cdot m)=\nabla_{\partial/\partial w_a}^{\frakg_X(V)}(J)\cdot
m+J\nabla_{\partial/\partial w_a}^{M_X}(m)$ for all $J\in\frakg_X(V),\,m\in M_X$.

\vskip 1ex
\noindent
{\rm (4)} $\nabla_{\partial/\partial w_a}^{\frakg_X(V)}(\frakg_X^{out}(V))\subset\frakg_X^{out}(V)$.
\end{proposition}
\begin{proof}
The statements (1) and (2) are immediate consequences of the definition. To show (3), 
it suffices to consider $J_1,\,J_2\in\frakg(V)$ by (2);
\begin{multline*}
\nabla_{\partial/\partial w_a}^{\frakg_X(V)}([J_1,J_2])
=[T(-1),[J_1,J_2]]\\=[[T(-1),J_1],J_2]]+[[J_1,[T(-1),J_2]]=[\nabla_{\partial/\partial
w_a}^{\frakg_X(V)}(J_1),J_2]+[J_1,\nabla_{\partial/\partial w_a}^{\frakg_X(V)}(J_2)].
\end{multline*}

To see (4) we note that that $[T(-1),J(v,f)]=J(T(-1)v,f)$ and $J(v,df/d\xi)+J(T(-1)v,f)=0$. 
Then for any $v\in V_\Delta$ and $f(w,z)(dz)^{1-\Delta}\in\pi_*\Omega^{1-\Delta}_{\calc/X}(*S_A)$ we see that
\[
\begin{split}
&\nabla_{\partial/\partial w_a}^{\frakg_X(V)}j_A\left(v\otimes f(w,z)(dz)^{1-\Delta}\right)\\
&=\sum_{b\in A}\nabla_{\partial/\partial w_a}^{\frakg_X(V)}J(v,f_b(w,\xi_b))\\
&=\sum_{b\in A}J\left(v,\frac{\partial}{\partial w_a}f_b(w,\xi_b)\right)+\left[T(-1),J(v,f_a(w,\xi_a))\right]\\
&=\sum_{b\in A}J\left(v,\left(\frac{\partial f}{\partial w_a}\right)_b(w,\xi_b)\right)
+J\left(v,\frac{\partial f_a}{\partial \xi_a}(w,\xi_a)\right)
+J(T(-1)v,f_a(w,\xi_a))\\ 
&=j_A\left(v\otimes\partial_af(w,z)(dz)^{1-\Delta}\right).
\end{split}
\]
\end{proof}

\begin{definition}
The quasi-coherent $\calo_X$-module 
$\calv_X(M_A)=M_X/\frakg_X^{out}(V)M_X$ is called the \textit{sheaf of covacua}.
\end{definition}

\subsection{Structure of sheaves of covacua}

\begin{theorem}\label{th:free} 
Let $V$ be a $C_2$-finite chiral vertex operator algebra. 

\vskip 1.5ex
\noindent
{\rm (1)} The sheaf of covacua $\calv_X(M_A)$ is a coherent $\calo_X$-module.

\vskip 1ex
\noindent
{\rm (2)} The linear map $\nabla_{\partial/\partial w_a}^{M_X}\,(a\in A)$ induce an integrable connection 
$\nabla_{\partial/\partial w_a}$ on the sheaf of covacua $\calv_X(M_A)$.

\vskip 1ex
\noindent
{\rm (3)} The sheaf of covacua $\calv_X(M_A)$ is a locally free $\calo_X$-module of finite rank.
\end{theorem}
\begin{proof}
By the proof of Theorem \ref{th:finiteblock} we see that there exists a finite dimensional vector space $U$ such that
$\calo_X\otimes U\rightarrow \gr_\bullet^F\calv_X(M_A)\rightarrow0$ is exact.
So we see that $\calv_X(M_A)$ is finitely generated $\calo_X$-module, and is then $\calo_X$-coherent since $\calo_X$ is noetherian.

By Proposition \ref{prop:connection} (4)  the linear map $\nabla_{\partial/\partial w_a}^{M_X}\,(a\in A)$ induce 
an integrable connection on $\calv_X(M_A)$. It is well-known that a coherent $\calo_X$-module with an integrable connection 
is a locally free $\calo_X$-module of finite rank (cf. \cite[(VI, 1.7)]{Bor}).
\end{proof}

\section{Regularity conditions}\label{sec:regularity}
\subsection{Condition I\!I and I\!I\!I}
\begin{definition}
We set $M(0)=\frakU(V)/\frakU(V)F^{1}\frakU(V)$, which is canonically $\frakU(V)$-module,
and denote by $\unit$ the element in $M(0)$ with the representative $1\in\frakU(V)$.
\end{definition}

\begin{remark}
There exists a canonical $\calu(V)$-module isomorphism 
\[
M(0)\cong\calu(V)\otimes_{F^0\calu(V)}A_0(V).
\]
\end{remark}

We have so far worked under the condition that the chiral vertex operator algebra is $C_2$-finite, 
which we call the condition I. 

\vskip 2ex
\noindent
\textbf{Condition I.}
The chiral vertex operator algebra $V$ is $C_2$-finite.
\vskip 2ex

We further assume that the zero-mode algebra $A_0(V)$ is semisimple, which we call the  condition I\!I:

\vskip 2ex
\noindent
\textbf{Condition I\!I.}
The chiral vertex operator algebra $V$ is $C_2$-finite and its zero-mode algebra $A_0(V)$ is \textit{semisimple}\,.
\vskip 2ex

The zero-mode algebra $A_0(V)$ is finite dimensional if $V$ is $C_2$-finite.
The category of finite dimensional $A_0(V)$-modules has finitely many 
simple objects upto equivalence; let $W_0,\,W_1,\,\ldots,\,W_r$ be the complete list of simple objects of this  category, where we set
$W_0=\C\ket{0}$ (see Proposition \ref{prop:zero-trivial}).

Any $A_0(V)$-module $W_i\,(0\leq i\leq r)$ is an $F^0\frakU(V)$-module by the natural projection $F^0\frakU(V)\rightarrow
A_0(V)\rightarrow 0$. Then we can define the induced $\frakU(V)$-module
$L_i=\frakU(V)\otimes_{F^0\frakU(V)}W_i\quad (0\leq i\leq r)$, which are in fact $V$-modules.

\vskip 2ex
\noindent
\textbf{Problem.} If a chiral vertex operator algebra $V$ is a $C_2$-finite simple $V$-module and its zero-mode algebra
$A_0(V)$ is semisimple then induced modules
$L_i\;(0\leq i\leq r)$ are simple $V$-modules.
\vskip 2ex

At the moment we are not able to prove this problem, instead we assume:

\vskip 2ex
\noindent
\textbf{Condition I\!I\!I.}
The chiral vertex operator algebra $V$ is $C_2$-finite, its zero-mode algebra $A_0(V)$ is semisimple, and the induced $V$-modules
$L_i\;(0\leq i\leq r)$ are simple.
\vskip 2ex

\subsection{Structure of the space of covacua}
In this subsection we assume the condition \III\!.
Recall that the zero-mode algebra $A_0(V)$ is linearly spanned by $[J_0(v)]\,(v\in V_\Delta)$.
We define a linear map $\theta:A_0(V)\rightarrow A_0(V)$ by
\[
J_0(v)\longmapsto(-1)^\Delta J_0(e^{T(1)}v)\quad\mbox{for all $v\in V_\Delta\,(\Delta\geq0)$}.
\]
Then $\theta$ is well-defined and is an anti-algebra involution of $A_0(V)$ (see Proposition \ref{prop:theta-involution}).

Let $W$ be a finite dimensional $A_0(V)$-module and set $D(W)=\Hom_\C(W,\C)$, which is a left $A_0(V)$-module by
$\pairing{J_0(v)\varphi}{u}=\pairing{\,\varphi}{\theta(J_0(v))u}$ for all $\varphi\in D(W),\,u\in W$ and $[J_0(v)]\in A_0(V)$.
The map $D:\calm\!\operatorname{od}(A_0(V))\rightarrow\calm\!\operatorname{od}(A_0(V))$ is a contravariant functor such that
$D^2=\id$. For any simple $A_0(V)$-module $W_i$ its dual module
$D(W_i)$ is also simple, and is isomorphic to one of the simple
$A_0(V)$-modules, which we denote by $W_{i^\dagger}$, i.e., $D(W_i)\cong W_{i^\dagger}$.

\begin{proposition}
Let $V$ be a chiral vertex operator algebra which satisfies conditions \III\!, and $W$ be a simple $A_0(V)$-module. 
Then there exists a canonical $V$-module isomorphism
$D(L(W))\cong L(D(W))$ where $L(W)=\mathcal{U}(V)\bigotimes_{F^0\mathcal{U}(V)}W$ is the induced $V$-module. 
\end{proposition}

Let $A=\{\,1,\,2,\,\ldots,\,N,\,\infty\,\}$ and set $\bar{A}=A\setminus\{\infty\}$.
Let $w_A=(w_a)_{a\in A}$ with $w_\infty=\infty$.
For $V$-modules $M^a\,(a\in \bar{A})$ we set $M_{A}=\bigotimes_{a\in A}M^a=M_{\bar{A}}\otimes M(0)$ where
$M_{\bar{A}}=\bigotimes_{a\in\bar{A}}M^a$ and $M^\infty=M(0)$. Let
$\frakg_{w_{\bar{A}}}^{out}(\infty)$ be the image of
$\bigoplus_{\Delta=0}^\infty V_\Delta\otimes H^0(\CP^1,\Omega^{1-\Delta}(*S_{\bar{A}}-\Delta[\infty]))$ in 
$\frakg_{\bar{A}}(V)$ under the map $j_{\bar{A}}$.
We denote by $\unit_{\infty}$ the canonical generator of $M(0)$. 

\begin{proposition}\label{prop:fusion}
Let $\iota:M_{\bar{A}}\rightarrow M_{\bar{A}}\otimes M(0)=M_A$ be the linear map defined by $u\mapsto u\otimes\unit_{\infty}$.
Then $\iota$ induces a vector space isomorphism
\[
M_{\bar{A}}/\frakg_{w_{\bar{A}}}^{out}(\infty)M_{\bar{A}}\cong\calv_{w_{A}}(M_{A}).
\]
\end{proposition}
\begin{proof}
We set $\frakg_\infty^{out}(\infty)=j_\infty\left(\bigoplus_{\Delta=0}^\infty V_\Delta\otimes
H^0\left(\CP^1,\Omega^{1-\Delta}(*S_{\bar{A}}-\Delta[\infty])\right)
\right)$, which is a Lie subalgebra of $\frakg(V)_\infty$. 
Since $\frakg_\infty^{out}(\infty)\subset F^{1}\frakg(V)_\infty$ and $F^{1}\frakg(V)_\infty\unit_\infty=0$ 
we have $j_{\bar{A}}(v\otimes\varphi)u\otimes\unit_{\infty} =j_{A}(v\otimes\varphi)(u\otimes\unit_{\infty})$ for any
$v\in V_\Delta$ and
$\varphi\in H^0(\CP^1,\Omega^{1-\Delta}(*S_{\bar{A}}-\Delta[\infty]))$, i.e., 
$\iota(\frakg_{w_{\bar{A}}}^{out}(\infty)M_{\bar{A}})\subset\frakg_{w_{A}}^{out}(V)M_{A}$;
in particular $\iota$ induces the linear map
$M_{\bar{A}}/\frakg_{w_{\bar{A}}}^{out}(\infty)M_{\bar{A}}\rightarrow\calv_{w_{A}}(M_{A})$.

In order to show that $\iota$ is surjective it suffices to prove that $M(0)$ is generated by $\unit_\infty$ as a
$\frakg_\infty^{out}(V)$-module where $\frakg_\infty^{out}(V)=j_\infty\left(\bigoplus_{\Delta=0}^\infty V_\Delta\otimes
H^0\left(\CP^1,\Omega^{1-\Delta}(*S_A)\right)\right)$; this follows from the fact that
$\varphi_n=z^{n+\Delta-1}(dz)^{1-\Delta}$ belongs to 
$H^0(\CP^1,\Omega^{1-\Delta}(*S_A))$ for any nonnegative integer $n$ and
$j_\infty(v\otimes\varphi_n)=(-1)^{\Delta+1}J_{-n}(e^{T(1)}v)$ or equivalently
$j_\infty(e^{-T(1)}v\otimes\varphi)=(-1)^{\Delta+1}J_{-n}(v)$ for all $v\in V_\Delta$.

The injectivity of the map $\iota$ is proved in the same spirit appeared in the proof of Theorem
\ref{th:propvacua}, i.e., the surjectivity of the dual map
\[
\iota^*:\calv_{w_{A}}^\dagger(M_{A})\rightarrow\Hom_\C(M_{\bar{A}}/\frakg_{w_{\bar{A}}}^{out}(\infty)M_{\bar{A}},\C)
\]
is proved by using the one to one correspondence between
$\Hom_\C(M_{\bar{A}}/\frakg_{w_{\bar{A}}}^{out}(\infty)M_{\bar{A}},\C)$ and the space of systems of current correlation
functions $\{\Phi_m\}_{m=0}^\infty$ such that
\[
\Phi_m(v_1,\ldots,v_m;u)_{(z_1,\ldots,z_m)}\in 
H^0((\CP^1)^m,\Omega^{\boxtimes_{i=1}^m\Delta_i}(\sum_{a\in\bar{A},\,i}*D_{a,\,i}+\sum_{i\neq
j}*D_{i,\,j}+\sum_{i}\Delta_iD_{i,\infty})).
\]
Let $\Phi\in\Hom_\C(M_{\bar{A}}/\frakg_{w_{\bar{A}}}^{out}(\infty)M_{\bar{A}},\C)$. By Proposition \ref{prop:orth1} there exists a
meromorphic form $\bra{\Phi}J(v,z)\ket{u_{\bar{A}}}(dz)^{\Delta}\in H^0(\CP^1,\Omega^{\Delta}(*w_{\bar{A}}+\Delta[\infty]))$ whose
Laurent expansion at $z=w_a\,(a\in\bar{A})$ is
\[
\bra{\Phi}J(v,z)\ket{u_{\bar{A}}}(dz)^{\Delta}
=\sum_{n\in\Z}\Phi\left(\rho_a\left(J_n(v)\right)u_{\bar{A}}\right)(z-w_a)^{-n-\Delta}(dz)^{\Delta}.
\]
There is another expression of this $1$-point current correlation function;
\begin{equation}\label{eqn:m1-point}
\bra{\Phi}J(v,z)\ket{u_{\bar{A}}}(dz)^{\Delta}
=\sum_{a\in\bar{A}}\sum_{n\geq-\Delta+1}\Phi\left(\rho_a\left(J_n(v)\right)u_{\bar{A}}\right)(z-w_a)^{-n-\Delta}(dz)^{\Delta}.
\end{equation}
To see this it suffices to show that the right hand side of (\ref{eqn:m1-point}) belongs to the space
$H^0(\CP^1,\Omega^{\Delta}(*w_{\bar{A}}+\Delta[\infty]))$, which follows from the fact that for all $0\leq n\leq
\Delta-2$ the form $z^n(dz)^{1-\Delta}$ belongs to $H^0(\CP^1,\Omega^{1-\Delta}(*w_{\bar{A}}-\Delta[\infty]))$  and then
$\Phi(j_{\bar{A}}(v\otimes z^n(dz)^{1-\Delta})u_{\bar{A}})=0$. Since the second equality of Lemma \ref{lem:inv} holds for
$\Phi\in\Hom_\C(M_{\bar{A}}/\frakg_{w_{\bar{A}}}^{out}(\infty)M_{\bar{A}},\C)$ we see that (\ref{eqn:2point-noinfty}) defines
the $2$-point current correlation function 
\[
\Phi_2(v_1,v_2;u_{\bar{A}})_{(z_1,z_2)}\in H^0((\CP^1)^2,
\Omega^{\Delta_1\boxtimes\Delta_2}(\sum_{\substack{i=1,2\\a\in\bar{A}}}*D_{i,a}+*D_{1,2}+\Delta_1D_{1,\infty}+\Delta_2
D_{2,\infty}))
\]
by using the same argument given in the proof of Proposition \ref{pro:expan1}.
Then by the same calculation in Theorem \ref{th:propvacua} we see that $\bra{\Phi}J(v,z)\ket{\,*}\in
\Hom_\C(M_{\bar{A}}/\frakg_{w_{\bar{A}}}^{out}(\infty)M_{\bar{A}},\C)$ for any fixed $v\in V_\Delta$ and $z\in\C$.
We now can inductively construct the whole system of current correlation functions associated to $\Phi$,
which we denote by
\[
\bra{\Phi}J(v_1,z_1)\cdots J(v_m,z_m)\ket{u_{\bar{A}}}(dz_1)^{\Delta_1}\cdots(dz_m)^{\Delta_m}.
\]
We now define $\tilde{\Phi}\in\Hom_\C(M_{\bar{A}}\otimes M(0),\C)$ by
\begin{multline*}
\tilde{\Phi}(u_{\bar{A}}\otimes\theta\left(J_{n_1}(v_1)\right)\cdots\theta\left(J_{n_k}(v_k)\right)\unit_\infty)
=
\left(\frac{1}{2\pi\sqrt{-1}}\right)^k\\
\times\oint_{|z_1|=r_1}\!\!\!\!\!\cdots\;\oint_{|z_k|=r_k}\bra{\Phi}J(v_k,z_k)\cdots
J(v_1,z_1)\ket{u_{\bar{A}}} z_k^{n_k+\Delta_k-1}\cdots z_1^{n_1+\Delta_1-1}\,dz_k\cdots dz_1,
\end{multline*}
where $r_k>r_{k-1}>\ldots>r_1>|w_a|\,(a\in\bar{A})$. 
Using the operator product property and the $\cals_m$-symmetry of the system $\{\Phi_m\}_{m=0}^\infty$ it is proved that
$\tilde{\Phi}$ defines an element in $\Hom(M_{\bar{A}}\otimes M(0),\C)$, i.e., it can be shown that the definition of $\tilde{\Phi}$
is compatible with the commutator formula, the associativity, etc.  The operator product property also shows that  $\tilde{\Phi}$
induce an element of $\calv_{w_{A}}^\dagger(M_{A})$, i.e., $\tilde{\Phi}(\frakg_A^{out}(V)M_A)=0$ (see the proof of Theorem
\ref{th:propvacua}).
\end{proof}

\begin{theorem}\label{th:decomptow}
Let $V$ be a chiral vertex operator algebra which satisfies conditions \III\!, and $W_i\,(0\leq i\leq r)$ be the complete list of
simple objects of $\calm\!\operatorname{od}(A_0(V))$ and $L_i$ be the induced $V$-module.
Let $M^a\,(a\in\bar{A})$ be $V$-modules and set $M_{\bar{A}}=\bigotimes_{a\in\bar{A}}M^a$.
Then there exists a canonical vector space isomorphism
\[
\calv_{w_{A}}(M_{A})\cong\bigoplus_{i=0}^{r}\calv_{w_{A}}(M_{\bar{A}}\otimes L_i)\otimes W_i^*,
\]
where $M_A=M_{\bar{A}}\otimes M(0)$.
\end{theorem}
\begin{proof}
Since $A_0(V)\cong\bigoplus_{i=0}^{r}W_i\otimes W_i^*$ by Wedderburn's theorem it follows that 
$M(0)\cong \frakU(V)\otimes_{F^0\frakU(V)}A_0(V)\cong\bigoplus_{i=0}^{r}L_i\otimes W_i^*$
as $V$-modules.
We now conclude that 
\[
\begin{split}
M_A/\frakg_{w_{A}}^{out}(V)M_{A}&\cong\bigoplus_{i=0}^{r}M_{\bar{A}}\otimes L_i\otimes W_i^*
/\frakg_{w_{A}}^{out}(V)\bigoplus_{i=0}^{r}M_{\bar{A}}\otimes L_i\otimes W_i^*\\
&\cong\bigoplus_{i=0}^{r}\calv_{w_A}(M_{\bar{A}}\otimes L_i)\otimes W_i^*\;.
\end{split}
\]
\end{proof}

\section{Factorization}\label{sec:factorization}
The aim of this section is to prove the factorization theorem; one of the immediate applications of the theorem is the formula of the
dimension of conformal blocks for $N\geq4$ in terms of fusion rules, i.e., the dimension of the conformal blocks for $N=3$.  The
theorem is proved by pinching the $N$-pointed projective line into two copies of the projective line, each of which is suitably
pointed.  We first define the fiber bundle which describes the pinching procedure, and study the sheaf of covacua associated to the
fiber bundle. 

\subsection{Gauge condition}
Let $A=\{\,0,\,1,\,\ldots,N,\,\infty\,\}$ and $w_A=(w_a)_{a\in A}$ be the set of distinct points in $\CP^1$
with $w_\infty=\infty$. We denote by $X_A$ the set of such $w_A$'s, which is canonically embedded in $(\CP^1)^{N+2}$.
Let $Z_A=\{w_A\in X_A\;|\; w_0=0,\,w_1=1\,\}$. We write $X=X_A$ and $Z=Z_A$ for short. For $V$-modules 
$M^a\,(a\in A)$ we set $M_A=\bigotimes_{a\in A}M^a$.
We define the sheaf of covacua over $X$ and $Z$ by 
\begin{align*}
&\calv_X(M_A)=\calo_X\otimes M_A/\frakg_X^{out}(V)\left(\calo_X\otimes M_A\right),\mbox{ and}\\
&\calv_Z(M_A)=\calo_Z\otimes M_A/\frakg_Z^{out}(V)\left(\calo_Z\otimes M_A\right),
\end{align*}
respectively. There is a $\Theta_X$ and $\Theta_Z$-action (an integrable connection) on $\calv_X(M_A)$ and $\calv_Z(M_A)$, which is
respectively given by
\begin{align*}
\nabla_{\partial/\partial w_a}&=\frac{\partial}{\partial w_a}+\rho_a(T(-1))\quad\mbox{for all $a\in A$, and}\\
\nabla_{\partial/\partial w_a}&=\frac{\partial}{\partial w_a}+\rho_a(T(-1))\quad\mbox{for all $a\in A\setminus\{0,\,1\}$}.
\end{align*}

\begin{proposition}
There is a canonical $\calo_X$-module isomorphism 
\[
\calv_X(M_A)\cong\calo_X\otimes_{\calo_Z}\calv_Z(M_A).
\]
\end{proposition}

\subsection{Geometric settings}
Let $A$ and $B$ be finite sets with cardinality $|A|,\,|B|\geq3$, which has $3$ distinct fixed points
$0_A,\,1_A,\,\infty_A\in A$ and $0_B,\,1_B,\,\infty_B\in B$, respectively. 
We set $\bar{A}=A\setminus\{0_A\}$ and $\bar{B}=B\setminus\{\infty_B\}$. We define  $C=\bar{A}\sqcup\bar{B}=A\setminus\{0_A\}\sqcup
B\setminus\{\infty_B\}$ with fixed 3 points being assigned by $0_C=0_B,\,1_C=1_A,\,\infty_C=\infty_A$.
We set 
\[
U_A=\{w_A\in Z_A\;|\; |w_{A,a}|>1/2,\, a\neq0_A\;\},\quad
U_B=\{w_B\in Z_B\;|\; |w_{B,b}|<2,\, b\neq\infty_B\;\},
\]
which is an open set in $Z_A$ and $Z_B$ respectively.

We define an open subset in $\C$ with the coordinate $q$ by
$U_q^\times=\{\;q\in\C\;|\; 0<|q|<1/16\;\}$ and set $U^\times=U_q^\times\times U_A\times U_B$.
The manifold $U^\times$ is embedded into $Z_C$ by
$w_{C,a}=w_{A,a}\;(a\in\bar{A})$ and $w_{C,b}=qw_{B,b}\;(b\in\bar{B})$.
There are two coordinate systems $(w_{C,c})_{c\in C^\circ}$ and $(q, w_{A^\circ},w_{B^\circ})$ for $U^\times$ where
$A^\circ=A\setminus\{0_A,1_A,\infty_A\},\,B^\circ=B\setminus\{0_B,1_B,\infty_B\}$ and $C^\circ=C\setminus\{0_C,1_C,\infty_C\}$. 

In order to define sheaves of covacua over $U^\times$ we introduce the trivial $\CP^1$-bundle 
$\pi:\calc_{U^\times}=U^\times\times\CP^1\rightarrow U^\times$; we regard $\CP^1$ as $\C\cup\{\infty\}$ and let $z_C$ be 
the affine coordinate of $\C$. For any $c\in C$ we define the section $s_c:U^\times\rightarrow\calc_{U^\times}$ by
$s_c(w_C)=(w_C,w_{C,c})$  and divisors over $\calc_{U^\times}$ by $S_c^\times=s_c(U^\times)$ and $S_C^\times=\bigcup_{c\in
C}S_c^\times$. We define a local coordinate around $S_c^\times$ by 
\[
\xi_{\,C,c}=\begin{cases}
z_C-w_{C,c}&\quad(c\neq\infty),\\
z_C&\quad(c=\infty).
\end{cases}
\]

Let $M^c\,(c\in C)$ be $V$-modules and set $M_C=\bigotimes_{c\in C}M^c$. For all $c\in C$ we set 
\[
\frakg_{U^\times,c}(V)=\bigoplus_{\Delta=0}^\infty V_\Delta\otimes\calo_{U^\times}((\xi_{C,c}))(d\xi_{C,c})^{1-\Delta}/
\nabla\left(\bigoplus_{\Delta=0}^\infty V_\Delta\otimes\calo_{U^\times}((\xi_{C,c}))(d\xi_{C,c})^{-\Delta}\right),
\]
which is a Lie algebra over $\calo_{U^\times}$. The Lie algebra $\frakg_{U^\times}(V)=\bigoplus_{c\in C}\frakg_{U^\times,c}(V)$ 
over $\calo_{U^\times}$ canonically acts on $\calo_{U^\times}\otimes M_C$. We denote this action  by $\rho_C=\sum_{c\in
C}\rho_{C,c}$.  Let
$j_C:\bigoplus_{\Delta=0}^\infty
V_\Delta\otimes\pi_*\Omega_{\calc_{U^\times}/U^\times}^{1-\Delta}(*S_C^\times)\rightarrow\frakg_{U^\times}(V)$ be the linear map
defined by the Laurent expansion at the divisors $S_c^\times\,(c\in C)$ in terms of the coordinates $\xi_{C,c}$ and denote the image
by
$\frakg_{U^\times}^{out}(V)$. We set
\[
\calv_{U^\times}(M_C)=\calo_{U^\times}\otimes M_C/\frakg_{U^\times}^{out}(V)\left(\calo_{U^\times}\otimes M_C\right)
\]
The sheaf $\calv_{U^\times}(M_C)$ is a coherent $\calo_{U^\times}$-module and is equipped with the
integrable connection
$\nabla_{\partial/\partial w_{C,c}}\,(c\in C)$ so that  $\calv_{U^\times}(M_C)$ is a locally free $\calo_{U^\times}$-module of finite
rank.

We set $U_q=\{\;q\in\C\;|\; |q|<1/16\;\}\supset U_q^\times$ and $U=U_q\times U_A\times U_B\supset U^\times$,
and define the trivial $\CP^1\times\CP^1$-bundle over $U$ by $\tilde{\calc}_U=U\times\CP^1\times\CP^1\rightarrow U$.
There exists a canonical commutative embedding 
\[
\begin{array}{ccc}
\calc_{U^\times}&\longrightarrow&\tilde{\calc}_U\\
\Big\downarrow&&\Big\downarrow\\
U^\times&\longrightarrow&U
\end{array}
\quad\mbox{by}\quad 
(q,w_A,w_B,z_C)\mapsto (q,w_A,w_B,z_C,q^{-1}z_C).
\]
Using the homogeneous coordinates $([z_A^0:z_A^1],[z_B^0:z_B^1])$ of the fiber $\CP^1\times\CP^1$
the manifold $\calc_{U^\times}$ is described in $\tilde{\calc}_U$ by
\[
\calc_{U^\times}=\{\;(q,w_A,w_B,[z_A^0:z_A^1],[z_B^0:z_B^1])\;|\;z_A^0z_B^1=qz_A^1z_B^0,\,q\neq0\;\}.
\] 
We now define the algebraic set $\calc_U$ of $\tilde{\calc}_U$ by $z_A^0z_B^1=qz_A^1z_B^0$,
and denote the natural projection by $\pi:\calc_U\rightarrow U$.
\begin{proposition}
{\rm (1)} $\calc_U$ is a codimension $1$ closed complex submanifold in $\tilde{\calc}_U$.

\vskip 1ex
\noindent
{\rm (2)} The map $\pi:\calc_U\rightarrow U$ is flat and proper.
\end{proposition}

The space $\calc_U$ is non-singular, while the fibers over the points of the domain $D=U\setminus U^\times$ are singular.
In order to see more detail of the manifold $\calc_U$ we introduce the open covering of $\calc_U$ and the local coordinate systems;
$\calc_U=V_0\cup V_1\cup V_\infty$ where 
\begin{align*}
&V_{0}=\{\,(q,w_A,w_B,z_B)\;|\;|z_B|<4\} \quad z_A=qz_B,\\
&V_{1}=\{\, (w_A,w_B,z_A,t_B)\;|\;|z_A|<1/2,\, |t_B|<1/2\;\}\quad q=z_At_B,\\
&V_{\infty}=\{\,(q,w_A,w_B,t_A)\;|\;|t_A|<4\,\}\quad t_B=qt_A.
\end{align*}
where $z_A=z_A^0/z_A^1$ and $t_A=z_A^1/z_A^0$, etc.

Let $\Sigma$ is the set of all points $p\in\calc_U$ at which the differential 
$d\pi_p:T_p\calc_U\rightarrow T_{\pi(p)}U$ is not surjective.

\begin{proposition}
{\rm (1)} The set $\Sigma$ is a codimension $2$ submanifold of $\calc_U$ such that $\Sigma\cap(V_0\cup V_\infty)=\emptyset$ and
$\Sigma=\{\;(w_A,w_B,z_A,t_B)\;|\; z_A=t_B=0\;\}\subset V_{1}$.

\vskip 1ex
\noindent
{\rm (2)} The projection $\pi$ gives an isomorphism $\Sigma\cong D$.
\end{proposition}

For any $c\in C$ we define the section $s_c:U\rightarrow\calc_U$ by
\[
(q,w_A,w_B)\longmapsto 
\begin{cases}
(q,w_A,w_B,w_{A,a},q^{-1}w_{A,a})&\quad (a\in \bar{A}),\\
(q,w_A,w_B,qw_{B,b},w_{B,b})&\quad (b\in \bar{B}),
\end{cases}
\]
and the divisor $S_c=s_c(U)$ on $\calc_U$. We set $S_C=S_{\bar{A}}\cup S_{\bar{B}}$ where
$S_{\bar{A}}=\bigcup_{a\in\bar{A}}S_a$ and $S_{\bar{B}}=\bigcup_{b\in\bar{B}}S_b$.
\begin{lemma}\label{lem:div10}
{\rm (1)} 
Let $a\in\bar{A}$ and $b\in\bar{B}$. Then $S_a\subset V_\infty,\,S_a\cap(V_1\cup V_0)=\emptyset$
and $S_b\subset V_0,\,S_b\cap(V_1\cup V_\infty)=\emptyset$. In particular $S_c\cap\Sigma=\emptyset$ for any $c\in C$.

\vskip 1ex
\noindent 
{\rm (2)} The following is a local coordinates system of $\calc_U$ around $S_c\,(c\in C)$;
\begin{align*}
(q,w_A,w_B,\xi_{A,a}),&\quad\xi_{A,a}=z_A-w_{A,a}\quad(a\in A\setminus\{0_A,\infty_A\}),\\
(q,w_A,w_B,\xi_{A,\infty}),&\quad\xi_{A,a}=z_A=t_A^{-1}\quad(a=\infty_A),\\
(q,w_A,w_B,\xi_{B,b}),&\quad\xi_{B,b}=z_B-w_{B,b}\quad(b\in \bar{B}).
\end{align*}

\vskip 1ex
\noindent
{\rm (3)} $\xi_{C,a}=\xi_{A,a}$ on $S_a\cap\calc_{U^\times}$ for any $a\in A\setminus\{0_A,\infty_A\}$
and $q^{-1}\xi_{C,b}=\xi_{B,b}$ on $S_b\cap\calc_{U^\times}$ for any $b\in \bar{B}$.
\end{lemma}

The fibers of $\calc_U$ over $D$ are given by $z_At_B=0$, i.e., each fiber consists of two copies of the projective line 
with the identification of $z_A=0$ and $z_B=\infty$ as a simple double point; $\Sigma$ is the set of all these simple double points.
We now define open subsets $V_1^0$ and $V_1^\infty$ of $V_1$ by
\[
V_1^0=\{\, (w_A,w_B,z_A,t_B)\in V_1\;|\;t_B\neq0\;\}\;\mbox{and}\;
V_1^\infty=\{\, (w_A,w_B,z_A,t_B)\in V_1\;|\;z_A\neq0\;\}.
\]
Then we see that $V_1\setminus\Sigma=V_1^0\cup V_1^\infty$ and $\calc_U\setminus\Sigma=V_0\cup V_1^0\cup V_1^\infty\cup V_\infty$.

The rest of this subsection is devoted to defining the sheaf of cavacua $\calv_U(M_C)$ attached to the family of curves
$\pi:\calc_U\rightarrow U$ such that $\calv_U(M_C)|U^\times=\calv_{U^\times}(M_C)$. 

Let $\omega_{\,\calc_U/U}$ be the relative dualizing sheaf of $\pi:\calc_U\rightarrow U$. Note that $\calc_U$ and $U$ are
nonsingular and $\pi$ is flat. Then there exists an $\calo_{\calc_U}$-module isomorphism 
$\omega_{\,\calc_U/U}\cong\omega_{\calc_U}\otimes\pi^*\omega_U^{-1}$ 
where $\omega_{\,Y}$ is the canonical sheaf of a complex manifold $Y$. 
In particular $\omega_{\,\calc_U/U}$ is a locally free $\calo_{\calc_U}$-module of rank $1$.  

\begin{proposition}
{\rm (1)} $\omega_{\,\calc_U/U}$ is a locally free $\calo_{\calc_U}$-module of rank $1$.

\vskip 1ex
\noindent
{\rm (2)}
Let $\Omega^1_{\calc_U/U}=\Omega^1_{\calc_U}/\calo_{\calc_U}(d\pi^*\Omega^1_U)$. 
Then 
$\omega_{\,\calc_U/U}|\calc_U\setminus\Sigma=\Omega^1_{\calc_U/U}|\calc_U\setminus\Sigma$
and
\begin{align*}
&\omega_{\,\calc_U/U}|V_0=\calo_{V_0}dz_B,\quad\omega_{\,\calc_U/U}|V_\infty=\calo_{V_\infty}dt_A,\\
&\omega_{\,\calc_U/U}|V_1^0=\calo_{V_1^0}dt_B,\quad\omega_{\,\calc_U/U}|V_1^\infty=\calo_{V_1^\infty}dz_A.
\end{align*}
\end{proposition}

Let $\frakg_{U,\bar{A}}(V)=\bigoplus_{a\in\bar{A}}\frakg_{U,a}(V)$ and $\frakg_{U,\bar{B}}(V)=\bigoplus_{b\in\bar{B}}\frakg_{U,b}(V)$
where 
\begin{align*}
\frakg_{U,a}(V)&=\bigoplus_{\Delta=0}^\infty
V_\Delta\otimes\calo_U((\xi_{A,a}))(d\xi_{A,a})^{1-\Delta}/
\nabla\left(\bigoplus_{\Delta=0}^\infty
V_\Delta\otimes\calo_U((\xi_{A,a}))(d\xi_{A,a})^{-\Delta}\right),\\
\frakg_{U,b}(V)&=\bigoplus_{\Delta=0}^\infty V_\Delta\otimes\calo_U((\xi_{B,b}))(d\xi_{B,b})^{1-\Delta}/
\nabla\left(\bigoplus_{\Delta=0}^\infty
V_\Delta\otimes\calo_U((\xi_{B,b}))(d\xi_{B,b})^{-\Delta}\right).
\end{align*}
The Lie algebra $\frakg_U(V)=\frakg_{U,\bar{A}}(V)\oplus\frakg_{U,\bar{B}}(V)$ over $\calo_U$ canonically acts on
$\calo_U\otimes M_C=\calo_U\otimes M_{\bar{A}}\otimes M_{\bar{B}}$ where $M^a\,(a\in\bar{A}),\,M^b\,(b\in\bar{B})$ be $V$-modules
and $M_{\bar{A}}=\bigotimes_{a\in\bar{A}}M^a,\,M_{\bar{B}}=\bigotimes_{b\in\bar{B}}M^b$. We denote this action by
\[
\rho_{\bar{A},\bar{B}}=\sum_{a\in\bar{A}}\rho_{\bar{A},\bar{B},a}+\sum_{b\in\bar{B}}\rho_{\bar{A},\bar{B},b}.
\] 

Recall that $\xi_{C,a}=\xi_{A,a}\,(a\in\bar{A})$ and $\xi_{C,b}=q\xi_{B,b}\,(b\in\bar{B})$. Let
$g_{\bar{A},\bar{B}}^{\;\;C}:\frakg_{U^\times}(V)\rightarrow\frakg_U(V)$ be the
$\calo_{U^\times}$-module map $J(v,\xi_{C,a}^{n+\Delta-1})\mapsto J(v,\xi_{A,a}^{n+\Delta-1})\,(a\in\bar{A}),\, 
J(v,\xi_{C,b}^{n+\Delta-1})\mapsto q^nJ(v,\xi_{B,b}^{n+\Delta-1})\,(b\in\bar{B})$ for any $v\in V_\Delta$ and $n\in\Z$,
which gives rise to an $\calo_{U^\times}$-module isomorphism $\frakg_{U^\times}(V)\cong\frakg_U(V)$.
The Lie algebra $\frakg_{U^\times}(V)$ acts on $\calv_{U^\times}(M_C)$ by $\rho_C$, while the Lie algebra $\frakg_U(V)$ acts
on $\calv_{U}(M_C)$ by $\rho_{\bar{A},\bar{B}}$. There is an $\calo_{U^\times}$-module isomorphism
$G_{\bar{A},\bar{B}}^{\;\;C}:\calo_{U^\times}\otimes M_{\bar{A}}\otimes M_{\bar{B}}\rightarrow\calo_{U}\otimes
M_{\bar{A}}\otimes M_{\bar{B}}|{U^\times}$ such that
$G_{\bar{A},\bar{B}}^{\;\;C}(\ell_{U^\times}m_{U^\times})
=g_{\bar{A},\bar{B}}^{\;\;C}(\ell_{U^\times})G_{\bar{A},\bar{B}}^{\;\;C}(m_{U^\times})$
for all $\ell_{U^\times}\in\frakg_{U^\times}(V)$ and $m_{U^\times}\in\calo_{U^\times}\otimes M_{\bar{A}}\otimes M_{\bar{B}}$:
in fact let $M^b=\bigoplus_{i}\bigoplus_{n=0}^\infty M^b(h_i^b+n)$. 
Then for any $u_{\bar{A}}\in M_{\bar{A}}$ and $u_{\bar{B}}=\otimes_{b\in\bar{B}}u_b\in M_{\bar{B}},\,u_b\in M^b(h_i^b+d_b)$
the isomorphism $G_{\bar{A},\bar{B}}^{\;\;C}$ is given by 
\[
G_{\bar{A},\bar{B}}^{\;\;C}(u_{\bar{A}}\otimes u_{\bar{B}})=\left(\sum_{b\in\bar{B}}q^{-d_b}\right)u_{\bar{A}}\otimes u_{\bar{B}},
\]
where we remark that $G_{\bar{A},\bar{B}}^{\;\;C}(M_{\bar{A}}\otimes M_{\bar{B}})\not\subset M_{\bar{A}}\otimes M_{\bar{B}}$.

Let $j_{\bar{A},\bar{B}}=j_{\bar{A}}+j_{\bar{B}}$ where $j_{\bar{A}}:\bigoplus_{\Delta=0}^\infty
V_\Delta\otimes\pi_*\omega_{\calc_{U}/U}^{1-\Delta}(*S_{\bar{A}})\rightarrow\frakg_{U,\bar{A}}(V)\subset\frakg_U(V)$  and
$j_{\bar{B}}:\bigoplus_{\Delta=0}^\infty
V_\Delta\otimes\pi_*\omega_{\calc_{U}/U}^{1-\Delta}(*S_{\bar{B}})\rightarrow\frakg_{U,\bar{B}}(V)\subset\frakg_U(V)$ is the linear
map defined by the Laurent expansion at the divisors $S_{\bar{A}}$ and $S_{\bar{B}}$ in terms of coordinates
$\xi_{A,a}\,(a\in\bar{A})$ and
$\xi_{B,b}\,(b\in\bar{B})$, respectively. The Lie subalgebra $\frakg_U^{out}(V)$ of $\frakg_U(V)$ is defined to be the image of
the map $j_{\bar{A},\bar{B}}$. We see that
$g_{\bar{A},\bar{B}}^{\;\;C}(\frakg_{U^\times}^{out}(V))=\frakg_{U}^{out}(V)|U^{\times}$.

\begin{definition}
Let $M^c\,(c\in C)$ be $V$-modules and set $M_C=\bigotimes_{c\in C}M^c$. We define the $\calo_U$-module $\calv_U(M_C)$ by
\[
\calv_U(M_C)=\calo_U\otimes M_C/\frakg_U^{out}(V)(\calo_U\otimes M_C).
\]
\end{definition}

Using the same argument given in the proof of Theorem \ref{th:free} we get:
\begin{proposition}\label{prop:extension}
$\calv_U(M_C)$ is a coherent $\calo_U$-module, and $G_{\bar{A},\bar{B}}^{\;\;C}$ induces an $\calo_{U^\times}$-module
isomorphism $\calv_{U^\times}(M_C)\cong\left.\calv_U(M_C)\right|U^\times$.
\end{proposition}

\subsection{Integrable connections}
Let $\mathfrak{m}_D$ be the defining ideal of $D$ in $U$, i.e., $\mathfrak{m}_D=q\calo_U$ and set 
$\Theta_U(-\log D)=\{\ell\in\Theta_U\,|\,\ell(\mathfrak{m}_D)\subset\mathfrak{m}_D\}$, 
which is a Lie subalgebra of $\Theta_U$ over $\calo_U$, 
where for any complex manifold $Y$ we denote by $\Theta_Y$ the sheaf of germs of vector fields on $Y$.
Using the local coordinate system $(q, w_{A^\circ},\,w_{B^\circ})$ of $U$ we get
\[
\Theta_U(-\log D)=\bigoplus_{a\in A^\circ}\calo_U\frac{\partial}{\partial w_{A,a}}
\oplus\bigoplus_{b\in B^\circ}\calo_U\frac{\partial}{\partial w_{B,b}}
\oplus\calo_Uq\frac{\partial}{\partial q}\;.
\]
Recall that there are two coordinate systems $(w_{\,C,c})_{c\in C^\circ}$ and $(q, w_{A^\circ},w_{B^\circ})$ for $U^\times$; 
the relation between them is 
\begin{equation}\label{eqn:translaw}
\begin{gathered}
\displaystyle{\frac{\partial}{\partial w_{A,a}}=\frac{\partial}{\partial w_{\,C,a}}}\quad(a\in A^\circ),\quad
\displaystyle{\frac{\partial}{\partial w_{B,b}}=w_{C,1_B}\frac{\partial}{\partial w_{\,C,b}}}\quad(b\in B^\circ),\\
\displaystyle{q\frac{\partial}{\partial q}=\sum_{b\in B\setminus\{0_B,\infty_B\}} w_{\,C,b}\frac{\partial}{\partial w_{\,C,b}}}\;,
\end{gathered}
\end{equation}
where $w_{C,1_B}=q$ while $w_{B,1_B}=1$

Recall that the integrable connection on $\calv_{U^\times}(M_C)$ is defined by
\[
\nabla_{\partial/\partial w_{C,c}}=\frac{\partial}{\partial w_{C,c}}+\rho_{C,c}\left(T(-1)\right)\quad(c\in C).
\]
Using (\ref{eqn:translaw}) and the fact that $q\rho_{C,b}(T(-1))=\rho_{\bar{A},\bar{B},b}(T(-1))$ for all $b\in\bar{B}$ we see that
the corresponding action of $\Theta_{U^\times}$ on $\calv_U(M_C)|U^\times$ is
\begin{equation}\label{eqn:conn01}
\begin{gathered}
\nabla_{\partial/\partial w_{A,a}}=\frac{\partial}{\partial
w_{A,a}}+\rho_{\bar{A},\bar{B},a}\left(T(-1)\right)\quad (a\in A^\circ),\\ 
\nabla_{\partial/\partial w_{A,a}}=\frac{\partial}{\partial w_{B,b}}+\rho_{\bar{A},\bar{B},b}\left(T(-1)\right)\;\quad (b\in
B^\circ),\\ 
\nabla_{q\partial/\partial q}=q\frac{\partial}{\partial
q}+\!\!\sum_{b\in B\setminus\{0_B,\infty_B\}}w_{B,b}\rho_{\bar{A},\bar{B},b}\left(T(-1)\right).
\end{gathered}
\end{equation}
It should be remarked that $[q\partial/\partial q,\rho_{\bar{A},\bar{B},b}(J_n(v))]=-n\rho_{\bar{A},\bar{B},b}(J_n(v))$ for all
$b\in\bar{B}$, and in particular, $[q\partial/\partial q,\rho_{\bar{A},\bar{B},b}(T(-1))]=\rho_{\bar{A},\bar{B},b}(T(-1))$.
The formula (\ref{eqn:conn01}) defines an integrable action of $\Theta_U(-\log D)$ on $\calo_U\otimes
M_{\bar{A}}\otimes M_{\bar{B}}$ such that $\nabla_X\frakg_{U}^{out}(V)\subset\frakg_{U}^{out}(V)$ for all $X\in\Theta_U(-\log D)$;
this fact is proved by the exactly same calculation done in the proof of Proposition \ref{prop:connection} (4) since $S_C\cap
V_1=\emptyset$. So there is a
$\Theta_U(-\log D)$ action on
$\calv_U(M_C)$  whose restriction to $\calv_U(M_C)|U^\times\cong\calv_{U^\times}(M_C)$ coincides with the $\Theta_{U^\times}$ action.
In other words;

\begin{theorem}\label{th:extension-connection}
The $\Theta_{U^\times}$ action on $\calv_{U^\times}(M_C)$ is naturally extended to a $\Theta_{U}(-\log D)$ action on $\calv_{U}(M_C)$.
\end{theorem}

\subsection{Local freeness} 
We spend the following several subsections to prove that the sheaf $\calv_U(M_C)$ is a locally free $\calo_U$-module.

Recall that $\Theta_U(-\log D)$ acts on $\calv_U(M_C)$ and induces $\Theta_{U^\times}$ action on
$\calv_U(M_C)|U^\times\cong\calv_{U^\times}(M_C)$, which is a locally free $\calo_{U^\times}$-module.  Let $\calo_D$ be the structure
sheaf of $D$, i.e., $\calo_D=\calo_U/\mathfrak{m}_D$. 
Then $\overline{\calv}_D(M_C)=\calv_U(M_C)\otimes_{\calo_U}\calo_D$  is a coherent $\calo_D$-module, and 
there is a natural $\Theta_U(-\log D)$ action on $\overline{\calv}_D(M_C)$ since
$\Theta_U(-\log D)\left(\mathfrak{m}_D\calv_U(M_C)\right)\subset\mathfrak{m}_D\calv_U(M_C)$.
In particular $\Theta_D$ acts on $\overline{\calv}_D(M_C)$  and then $\overline{\calv}_D(M_C)$ is a locally free $\calo_D$-module of
finite rank. Since $\left.\calv_U(M_C)\right|U^\times\cong\calv_{U^\times}(M_C)$ and $\calv_U(M_C)$ is $\calo_U$-coherent it follows
that
$\rank\,\calv_{U^\times}(M_C)\leq\rank\,\overline{\calv}_D(M_C)$.  So in order to prove the local freeness of $\calv_U(M_C)$ over
$\calo_U$ it suffices to show that $\rank\,\calv_{U^\times}(M_C)=\rank\,\overline{\calv}_D(M_C)$.

We denote the restriction of $\calc_U$ to $D$ by $\pi_D:\calc_D=\pi^{-1}(D)\rightarrow D$, and set $\bar{S}_C=\bigcup_{c\in
C}s_c(D)$, which is a divisor over $\calc_D$. For any $a\in\bar{A}$ and $b\in\bar{B}$ we set 
\begin{align*}
\frakg_{D,a}(V)&=\bigoplus_{\Delta=0}^\infty
V_\Delta\otimes\calo_D((\xi_{A,a}))(d\xi_{A,a})^{1-\Delta}/
\nabla\left(\bigoplus_{\Delta=0}^\infty
V_\Delta\otimes\calo_D((\xi_{A,a}))(d\xi_{A,a})^{-\Delta}\right),\\
\frakg_{D,b}(V)&=\bigoplus_{\Delta=0}^\infty
V_\Delta\otimes\calo_D((\xi_{A,b}))(d\xi_{A,b})^{1-\Delta}/
\nabla\left(\bigoplus_{\Delta=0}^\infty
V_\Delta\otimes\calo_D((\xi_{A,b}))(d\xi_{A,b})^{-\Delta}\right)
\end{align*}
and define 
$\frakg_{D,\bar{A}}(V)=\bigoplus_{a\in\bar{A}}\frakg_{D,a}(V),\,\frakg_{D,\bar{B}}(V)=\bigoplus_{b\in\bar{A}}\frakg_{D,b}(V)$
and $\frakg_D(V)=\frakg_{D,\bar{A}}(V)\oplus\frakg_{D,\bar{B}}(V)(=\frakg_U(V)\bigotimes_{\calo_U}\calo_D)$,
which are Lie algebras over $\calo_D$. 

\begin{definition}
Let 
$\omega_{\,\calc_D/D}=\omega_{\,\calc_U/U}\otimes_{\calo_{\calc_U}}\calo_{\calc_D}$ be the invertible $\calo_{\calc_D}$-module. 
We set $\frakg_D^{out}(V)=j_{\bar{A},\bar{B}}(\bigoplus_{\Delta=0}^\infty
V_\Delta\otimes(\pi_D)_*\,\omega_{\,\calc_D/D}^{1-\Delta}(*\bar{S}_C))$ where the map $j_{\bar{A},\bar{B}}$ is also defined by
Laurent expansion in terms of local coordinates $\xi_{A,a}\,(a\in\bar{A})$ and $\xi_{B,b}\,(b\in\bar{B})$, which is an $\calo_D$-Lie
subalgebra of $\frakg_D(V)$. We define
\[
\calv_D(M_C)=\calo_D\otimes M_C/\frakg_D^{out}(V)(\calo_D\otimes M_C).
\]
\end{definition}

There is an  $\calo_D$-module isomorphism $\overline{\calv}_D(M_A)\cong\calv_D(M_A)$. 
In fact by using
\begin{equation}\label{eqn:cut10}
(\pi_D)_*\omega_{\calc_D/D}^{1-\Delta}(*\bar{S}_C)\cong\pi_*\omega_{\calc_U/U}^{1-\Delta}(*S_C)\otimes_{\calo_U}\calo_D
\end{equation}
(see \cite[Lemma 4.1.3, page 522]{TUY}) one has the $\calo_D$-module isomorphism
$\frakg_D^{out}(V)\cong\frakg_U^{out}(V)\otimes_{\calo_U}\calo_D$. In the following commutative diagram of $\calo_D$-modules
\[
\begin{array}{ccc}
&&0\\
&&\big\downarrow\\
\frakg_U^{out}(V)(\calo_U\otimes M_C)\otimes_{\calo_U}\calo_D&\overset{\beta}{\longrightarrow}&\frakg_D^{out}(V)(M_C\otimes\calo_D)\\
\Big\downarrow\vcenter{\llap{$\alpha\quad$}}&&\Big\downarrow\vcenter{\rlap{$\epsilon$}}\\
(\calo_U\otimes M_C)\otimes_{\calo_U}\calo_D&\overset{\gamma}{\longrightarrow}&M_C\otimes\calo_D\\
\Big\downarrow&&\Big\downarrow\\
\calv_U(M_C)\otimes_{\calo_U}\calo_D&\longrightarrow&\calv_D(M_C)\\
\big\downarrow&&\big\downarrow\\
0&&0
\end{array}
\]
the map $\beta$ is surjective and the map $\gamma$ is an isomorphism. So $\operatorname{Im}(\alpha)\cong\operatorname{Im}(\epsilon)$,
which shows $\calv_D(M_C)\cong M_C\otimes\calo_D/\operatorname{Im}(\epsilon)\cong(\calo_U\otimes
M_C)\otimes_{\calo_U}\calo_D/\operatorname{Im}(\alpha)\cong\overline{\calv}_D(M_C)$. 

\begin{lemma}\label{lem:vbar}
There is a canonical $\calo_D$-module isomorphism $\overline{\calv}_D(M_A)\cong\calv_D(M_A)$.
\end{lemma}

The structure of the sheaf $\calv_D(M_C)$ is described by the following theorem whose proof is given in \S\,\ref{sub:structure}.

\begin{theorem}\label{th:dstr}
There exists a canonical $\calo_D$-module isomorphism 
\[
\calv_D(M_C)\cong\bigoplus_{i=0}^{r}\calv_D^A(M_{\bar{A}}\otimes L_i)\otimes\calv_D^B(L_{i^\dagger}\otimes M_{\bar{B}}).
\]
\end{theorem}

More explanations here should be in order. Sheaves of covacua over $D$ in the right hand sides are defined by using
the trivial $\CP^1$-bundle $D\times \CP^1\rightarrow D$ and the cross-section $s_a:(w_A,w_B)\mapsto (w_A,w_B,w_a)\,(a\in A)$ and
$s_b:(w_A,w_B)\mapsto (w_A,w_B,w_b)\,(b\in B)$, respectively;
\begin{align*}
\calv_D^A(M_{\bar{A}}\otimes L_i)&=\calo_D\otimes M_{\bar{A}}\otimes L_i/\frakg_D^{A,out}(V)
\left(\calo_D\otimes M_{\bar{A}}\otimes
L_i\right),\\
\calv_D^B(L_{i^\dagger}\otimes M_{\bar{B}})&=\calo_D\otimes L_{i^\dagger}\otimes
M_{\bar{B}}/\frakg_D^{B,out}(V)\left(\calo_D\otimes L_{i^\dagger}\otimes M_{\bar{B}}\right),
\end{align*}
where the simple $V$-module $L_i$ and $L_{i^\dagger}$ corresponds to the point $0_A$ and $\infty_B$.
So these are sheaves of covacua attached to a family of nonsingular curves of genus $0$, while $\calv_D(M_C)$ is the sheaf 
of covacua attached to a family of stable curves of genus $0$.

\vskip 1.5ex
Let $\widehat{\calo}_{U}= \underset{\underset{n}{\longleftarrow}}{\lim}\;\calo_{U}/\mathfrak{m}_D^n\calo_{U}\cong
\calo_D[[q]]$ be the completion along $D$ and set
\begin{align*}
&\widehat{\calv}_{U}(M_C)=\underset{\underset{n}{\longleftarrow}}{\lim}\;\calv_{U}(M_C)/\mathfrak{m}_D^n\calv_{U}(M_C),\\
&\widehat{\Theta}_{U}(-\log D)=\underset{\underset{n}{\longleftarrow}}{\lim}\;\Theta_{U}(-\log
D)/\mathfrak{m}_D^n\Theta_{U}(-\log D).
\end{align*}
The $\widehat{\calv}_{U}(M_C)$ is a coherent $\widehat{\calo}_{U}$-module, and is equipped with the integrable connection
$\nabla_X$ for $X\in\widehat{\Theta}_{U}(-\log D)$.
There is a canonical $\calo_D$-module isomorphism
\[
\widehat{\calv}_{U}(M_C)\otimes_{\widehat{\calo}_U}\calo_D\cong\overline{\calv}_D(M_C).
\]

\begin{definition}
Let $L$ be a simple $V$-module such that $L=\bigoplus_{d=0}^\infty L(d\,)$ where $L(d\,)=\{u\in
L\;|\;T(0)u=(\Delta+d\,)u\;\}$. We define a $\widehat{\Theta}_{U}(-\log D)$ action on $\calv_D^A(M_{\bar{A}}\otimes
L)\otimes\calv_D^B(L^\dagger\otimes M_{\bar{B}})\otimes\C[[q]]$ for $\partial/\partial w_{A,a}\,(a\in A^\circ)$ and
$\partial/\partial w_{B,b}\,(b\in B^\circ)$  by (\ref{eqn:conn01}), while for $q\partial/\partial q$ by
\[
\nabla_{q\partial/\partial q}=q\frac{\partial}{\partial
q}+\Delta-\sum_{b\in\bar{B}}\rho_{\bar{A},\bar{B},b}\left(T(0)\right).
\]
This canonically induces a $\widehat{\Theta}_{U}(-\log D)$ action on
$\overline{\calv}_D(M_C)\otimes\C[[q]]$.
\end{definition}

The following theorem yields the local freeness of the sheaf $\calv_U(M_C)$:
\begin{theorem}\label{th:factor}
There exists a canonical $\calo_D[[q]]$-module isomorphism which is compatible with the action $\nabla_{\partial/\partial
w_{A,a}}\,(a\in\bar{A}),\,\nabla_{\partial/\partial w_{B,b}}\,(a\in\bar{B})$ and $\nabla_{q\partial/\partial q}$
such that the following commutative diagram holds:
\[
\begin{array}{ccc}
\widehat{\calv}_{U}(M_C)&\cong&\overline{\calv}_D(M_C)\otimes\C[[q]]\\
\downarrow&&\downarrow\\
\overline{\calv}_D(M_C)&\underset{\id}{\cong}&\overline{\calv}_D(M_C).
\end{array}
\]
\end{theorem}

\begin{theorem}\label{th:locally-free}
$\calv_U(M_C)$ is a locally free $\calo_U$-module.
\end{theorem}

\subsection{Relative dualizing sheaf over the critical discriminant} 
This subsection is devoted to the study of the direct image sheaf $\pi_*\omega_{\,\calc_U/U}^{1-\Delta}(*S_C)$ and
$(\pi_D)_*\,\omega_{\,\calc_D/D}^{1-\Delta}(*\bar{S}_C)$; the contents of the subsection play the essential role in the proofs of 
Theorem \ref{th:dstr} and Theorem \ref{th:factor}.

Let $\varphi$ be a cross-section of the sheaf $\pi_*\omega_{\,\calc_U/U}^{1-\Delta}(*S_C)$.
We denote the restriction of $\varphi$ to the open set $V_0,\,V_1^0,\,V_1^\infty$ and $V_\infty$ by
$\varphi_0,\,\varphi_1^0,\,\varphi_1^\infty$ and $\varphi_\infty$, respectively, so that
$\varphi_0=\varphi_1^0,\,\varphi_1^0=\varphi_1^\infty$ and $\varphi_1^\infty=\varphi_\infty$ 
on $V_0\cap V_1^0,\,V_1^0\cap V_1^\infty$ and $V_1^\infty\cap V_\infty$\,.  
Recall that the critical locus $\Sigma$ has codimension $2$ in $\calc_U$ and  $\pi_*\omega_{\,\calc_U/U}^{1-\Delta}(*S_C)$ is an
invertible sheaf over $U$. So $\varphi$ is uniquely determined by the value on $\calc_U\setminus\Sigma\,(=V_0\cup V_1^0\cup
V_1^\infty\cup V_\infty)$\,. 

By Lemma \ref{lem:div10} (1) the cross-section $\varphi_1^0$ and $\varphi_1^\infty$ is holomorphic on each domain of the definition.
The expansion of $\varphi_1^0$ and $\varphi_1^\infty$ is respectively given by
\begin{align*}
\varphi_1^0=&\sum_{m\geq0,\,n\in\Z}\varphi_{m,n}^B(w_A,w_B)z_A^mt_B^n(dt_B)^{1-\Delta}\quad
(t_B\neq0),\\
\varphi_1^\infty=&\sum_{m\in\Z,\,n\geq0}\varphi_{m,n}^A(w_A,w_B)z_A^mt_B^n(dz_A)^{1-\Delta}\quad (z_A\neq0).
\end{align*}
Since $\varphi_1^0=\varphi_1^\infty$ on $V_1^0\cap V_1^\infty$ and $t_Bdz_A+z_Adt_B=0$ in $\omega_{\,\calc_U/U}|V_1^0\cap
V_1^\infty$ we see that
\begin{align}
&\varphi_{m,n}^A=0\,(m<\Delta-1),\,\varphi_{m,n}^B=0\,(n<\Delta-1)\quad\mbox{and}\label{eqn:nil01}\\
&\varphi^A_{m+\Delta-1,n}=(-1)^{\Delta-1}\varphi^B_{m,n+\Delta-1}\quad(m,n\geq0).\label{eqn:ext}
\end{align}
Summarizing 
\begin{align}
\varphi_1^0=&\sum_{m,n\geq0}\varphi_{m,n+\Delta-1}^B(w_A,w_B)z_A^mt_B^{n+\Delta-1}(dt_B)^{1-\Delta}\notag\\
=&\sum_{m,n\geq0}\varphi_{m,n+\Delta-1}^B(w_A,w_B)q^mt_B^{n-m+\Delta-1}(dt_B)^{1-\Delta},\label{eqn:exp10}\\
\varphi_1^\infty=&\sum_{m,n\geq0}\varphi_{m+\Delta-1,n}^A(w_A,w_B)z_A^{m+\Delta-1}t_B^n(dz_A)^{1-\Delta}\notag\\
=&\sum_{m,n\geq0}\varphi_{m+\Delta-1,n}^A(w_A,w_B)z_A^{m-n+\Delta-1}q^n(dz_A)^{1-\Delta}.\label{eqn:exp11}
\end{align}

Let $\pi^A_D\sqcup\pi^B_D:\calc_D^A\sqcup\calc_D^B\rightarrow \calc_D$ be obtained by the simultaneous normalization of 
$\pi_D:\calc_D\rightarrow D$ along $\Sigma$: 
\[
\begin{array}{ccc}
\calc_D^A\sqcup\calc_D^B&\longrightarrow&\calc_D\\
\Big\downarrow\vcenter{\llap{$\pi^A_D\sqcup\pi^B_D\quad$}}&&\Big\downarrow\vcenter{\rlap{$\pi_D$}}\\
D&\longrightarrow&D
\end{array}
\]
where $\calc_D^A\sqcup\calc_D^B$ is the disjoint union of two copies of the trivial $\CP^1$-bundle over $D$.
Let $s_{A,0}:D\rightarrow \calc_D^A$ and $s_{B,\infty}:D\rightarrow \calc_D^B$ be defined
by $(w_A,w_B)\mapsto(w_A,w_B,0)$ and $(w_A,w_B)\mapsto(w_A,w_B,\infty)$, which correspond to the normalized double points over $D$,
and set $S_{A,0}=s_{A,0}(D)\subset\calc_D^A$ and $S_{B,\infty}=s_{B,\infty}(D)\subset\calc_D^B$. 
We set $S_{\bar{A}}=\cup_{a\in\bar{A}}s_a(D)\subset\calc_D^A$ and $S_{\bar{B}}=\cup_{b\in\bar{B}}s_b(D)\subset\calc_D^B$.
Note that we use same notations for divisors over $U$ and $D$, for instance, recall that we set
$S_{\bar{A}}=\cup_{a\in\bar{A}}s_a(U)$ as a divisor over $U$. 

By (\ref{eqn:cut10}) any cross-section of the sheaf $(\pi_D)_*\omega_{\,\calc_D/D}^{1-\Delta}(*\bar{S}_C)$ is represented by a
quadruplet 
$(\varphi_0,\,\varphi_1^0,\,\varphi_1^\infty,\,\varphi_\infty)$ where $\varphi_0$
and $\varphi_\infty$ is the restriction of some $\varphi\in\pi_*\omega_{\,\calc_U/U}^{1-\Delta}(*S_C)$ to $V_0\cap\Sigma$ and
$V_\infty\cap\Sigma$, respectively,  while $\varphi_1^0$ and
$\varphi_1^\infty$ is the unique extension of $\varphi|V_1^0$ and $\varphi|V_1^\infty$ to
$\Sigma$. Then by (\ref{eqn:nil01}) we see
that $\varphi_D^A=(\varphi_1^\infty,\varphi_\infty)\in(\pi_D^A)_*\Omega_{\,\calc_D^A/D}^{1-\Delta}(*S_{\bar{A}}-(\Delta-1)S_{A,0})$ and
$\varphi_D^B=(\varphi_0,\varphi_1^0)\in(\pi_D^B)_*\Omega_{\,\calc_D^B/D}^{1-\Delta}(*S_{\bar{B}}-(\Delta-1)S_{B,\infty})$.
The expansion of $\varphi_D^A$ and $\varphi_D^B$ in a neighborhood of $z_A=0$ and $t_B=0$ is respectively given by
\begin{align*}
\varphi_D^A=&\sum_{m\geq0}\varphi_{m+\Delta-1,0}^A(w_A,w_B)z_A^{m+\Delta-1}(dz_A)^{1-\Delta}\quad (z_A\neq0)\quad\mbox{and}\\
\varphi_D^B=&\sum_{n\geq0}\varphi_{0,n+\Delta-1}^B(w_A,w_B)t_B^{n+\Delta-1}(dt_B)^{1-\Delta}\quad (t_B\neq0).
\end{align*}
The leading terms $\varphi_{\Delta-1,0}^A,\,\varphi_{0,\,\Delta-1}^B\in\calo_D$ satisfy
$\varphi_{\Delta-1,0}^A+(-1)^{\Delta}\varphi_{0,\,\Delta-1}^B=0$ by (\ref{eqn:ext}).

\begin{proposition}\label{prop:embedding10}
The sequence of $\calo_D$-module homomorphisms
\begin{multline*}
0\rightarrow(\pi_D)_*\omega_{\,\calc_D/D}^{1-\Delta}(*\bar{S}_C)\rightarrow\\
(\pi_D^A)_*\Omega^{1-\Delta}_{\calc_D^A/D}(*S_{\bar{A}}-(\Delta-1)S_{A,0})
\oplus(\pi_D^B)_*\Omega^{1-\Delta}_{\calc_D^B/D}(*S_{\bar{B}}-(\Delta-1)S_{B,\infty})
\rightarrow\calo_D\rightarrow0
\end{multline*}
defined by $\varphi\mapsto(\varphi_D^A,\varphi_D^B)$ and
$(\varphi_D^A,\varphi_D^B)\mapsto\varphi_{\Delta-1,0}^A+(-1)^{\Delta}\varphi_{0,\,\Delta-1}^B$ is exact.
\end{proposition}
\begin{proof}
Let $\varphi_{\Delta-1,0}^A,\,\varphi_{0,\,\Delta-1}^B\in\calo_D$ such that $\varphi_{\Delta-1,0}^A+(-1)^{\Delta}\varphi_{0,\,\Delta-1}^B=0$. 
Then by (\ref{eqn:cut10}) there exists
$\varphi\in\pi_*\omega_{\,\calc_U/U}^{1-\Delta}(*S_C)$ such that
$\varphi|D=(\varphi_{D}^A,\varphi_{D}^B)$ and the leading term of $\varphi_{D}^A$ and $\varphi_{D}^B$ is
$\varphi_{\Delta-1,0}^A$ and $\varphi_{0,\,\Delta-1}^B$, respectively. 
\end{proof}

\subsection{Proof of the structure theorem for the sheaf over critical discriminant}
\label{sub:structure}
In this subsection we prove Theorem \ref{th:dstr}. 
By using the exact sequence in Proposition \ref{prop:embedding10} we identity
$(\pi_D)_*\omega^{1-\Delta}_{\calc_D/D}(*\bar{S}_C)$ as an $\calo_D$-submodule of
$(\pi_D^A)_*\Omega^{1-\Delta}_{\calc_D^A/D}(*S_{\bar{A}}-(\Delta-1)S_{A,0})
\oplus(\pi_D^B)_*\Omega^{1-\Delta}_{\calc_D^B/D}(*S_{\bar{B}}-(\Delta-1)S_{B,\infty})$.
We introduce the vector subspace 
\[
(\pi_D)_*\omega^{1-\Delta}_{\calc_D/D}(*\bar{S}_C)'=
(\pi_D^A)_*\Omega^{1-\Delta}_{\calc_D^A/D}(*S_{\bar{A}}-\Delta S_{A,0})
\oplus(\pi_D^B)_*\Omega^{1-\Delta}_{\calc_D^B/D}(*S_{\bar{B}}-\Delta S_{B,\infty}),
\]
which has the codimension $1$ in $\pi_*\omega^{1-\Delta}_{\calc_D/D}(*S_C)$. 
We define the Lie subalgebra $\frakg_D^{out\,'}(V)$ of $\frakg_D(V)$ by
\[
\frakg_D^{out\,'}(V)=j_{\bar{A},\bar{B}}\left(\bigoplus_{\Delta=0}^\infty
V_\Delta\otimes(\pi_D)_*\omega^{1-\Delta}_{\calc_D/D}(*\bar{S}_C)'\right),
\]
and set 
\[
\calv_D(M_C)'=\calo_D\otimes M_C/\frakg_D^{out\,'}(V)(\calo_D\otimes M_C).
\]
In order to prove Theorem \ref{th:dstr} it suffices to determine $\calv_D(M_C)'$ and the kernel of the canonical
$\calo_D$-module homomorphism
$\calv_D(M_C)'\rightarrow\calv_D(M_C)\rightarrow 0$. 
With this point in view we introduce the $\calo_D$-module $\calv_D^A(M_{\bar{A}})'=\calo_D\otimes
M_{\bar{A}}/\frakg_D^{A,out\,'}(V)(\calo_D\otimes M_{\bar{A}})$ and
$\calv_D^B(M_{\bar{B}})'=\calo_D\otimes M_{\bar{B}}/\frakg_D^{B,out\,'}(V)(\calo_D\otimes M_{\bar{B}})$ where 
\begin{align*}
&\frakg_D^{A,out\,'}(V)=j_{\bar{A}}\big(\bigoplus_{\Delta=0}^\infty
V_\Delta\otimes(\pi_D^A)_*\Omega^{1-\Delta}_{\calc_D^A/D}(*S_{\bar{A}}-\Delta  S_{A,0})\big),\\
&\frakg_D^{B,out\,'}(V)=j_{\bar{B}}\big(\bigoplus_{\Delta=0}^\infty
V_\Delta\otimes(\pi_D^B)_*\Omega^{1-\Delta}_{\calc_D^B/D}(*S_{\bar{B}}-\Delta S_{B,\infty})\big).
\end{align*}

Using Proposition \ref{prop:fusion} we get:
\begin{lemma}\label{lem:iso'}
The linear map $M_{\bar{A}}\rightarrow M_{\bar{A}}\otimes M(0)\,(u_{\bar{A}}\mapsto u_{\bar{A}}\otimes\unit_{A,0})$ and
$M_{\bar{B}}\rightarrow  M(0)\otimes M_{\bar{B}}\,(u_{\bar{B}}\mapsto\unit_{B,\infty}\otimes u_{\bar{B}})$ respectively induces an
$\calo_D$-module isomorphism 
\[
\calv_D^A(M_{\bar{A}})'\cong\calv_D^A(M_{\bar{A}}\otimes M(0))\quad\mbox{and}\quad
\calv_D^B(M_{\bar{B}})'\cong\calv_D^B(M(0)\otimes M_{\bar{B}}).
\]
\end{lemma}

\begin{corollary}
There exists a canonical $\calo_D$-module isomorphism
\[
\begin{array}{ccc}
\calv_D(M_C)'&\cong&\calv_D^A(M_{\bar{A}}\otimes M(0))\otimes\calv_D^B(M(0)\otimes M_{\bar{B}})\\
&&\\
u_{\bar{A}}\otimes u_{\bar{B}}&\mapsto&(u_{\bar{A}}\otimes\unit_{A,0})\otimes(\unit_{B,\infty}\otimes u_{\bar{B}}).
\end{array}
\]
\end{corollary}

We now determine the kernel of $\calo_D$-module homomorphism
\begin{equation}\label{eqn:seq}
\calv_D^A(M_{\bar{A}}\otimes M(0))\otimes\calv_D^B(M_{\bar{B}}\otimes M(0))
\longrightarrow\calv_D(M_C)\longrightarrow0.
\end{equation}
Note that for any $\Delta\geq0$ the pair
\begin{align*}
&\varphi_D^{A,\,\Delta}=(z_A^{-1}dz_A)^{1-\Delta}\in(\pi_D^A)_*\Omega^{1-\Delta}_{\calc_D^A/D}(*S_{\bar{A}}-(\Delta-1)S_{A,0})
\quad\mbox{and}\\
&\varphi_D^{B,\,\Delta}=(-t_B^{-1}dt_B)^{1-\Delta}\in(\pi_D^B)_*\Omega^{1-\Delta}_{\calc_D^B/D}(*S_{\bar{B}}-(\Delta-1)S_{B,\infty})
\end{align*}
defines an element $\varphi_D^\Delta\in(\pi_D)_*\omega^{1-\Delta}_{\calc_D/D}(*\bar{S}_C)$ such that
\[
(\pi_D)_*\omega^{1-\Delta}_{\calc_D/D}(*\bar{S}_C)
=(\pi_D)_*\omega^{1-\Delta}_{\calc_D/D}(*\bar{S}_C)'\oplus\calo_D\,\varphi^\Delta_D.
\]
This gives rise to the splitting $\frakg_D^{out}(V)=\frakg_D^{out\,'}(V)\oplus\frakg_{D}^{out\,''}(V)$ where
$\frakg_{D}^{out\,''}(V)$ is the $\calo_D$-module generated by 
$j_{\bar{A},\bar{B}}\left(v\otimes\varphi_D^\Delta\right)$ for all $v\in V_\Delta\,(\Delta\geq0)$.

Since $[j_{\bar{A},\bar{B}}(v\otimes\varphi_D^\Delta),\frakg_D^{out\,'}(V)]\subset\frakg_D^{out\,'}(V)$ for all $v\in V_\Delta$
the element $j_{\bar{A},\bar{B}}(v\otimes\varphi_D^\Delta)$ acts on $\calv_D(M_C)'$, which action is also denoted by
$j_{\bar{A},\bar{B}}(v\otimes\varphi_D^\Delta)$; more precisely $j_{\bar{A}}(v\otimes\varphi_D^{A,\Delta})$ and 
$j_{\bar{B}}(v\otimes\varphi_D^{B,\Delta})$ acts on $\calv_D^A(M_{\bar{A}})'$ and $\calv_D^B(M_{\bar{B}})'$, respectively.  
The corresponding action on $\calv_D^A(M_{\bar{A}}\otimes M(0))$ and $\calv_D^B( M(0)\otimes M_{\bar{B}})$
through the isomorphism given in Lemma \ref{lem:iso'} is described as follows: 
recall that these isomorphisms are induced from linear maps
$u_{\bar{A}}\mapsto u_{\bar{A}}\otimes\unit_{A,0}\,(u_{\bar{A}}\in M_{\bar{A}})$ and $u_{\bar{B}}\mapsto\unit_{B,\infty}\otimes
u_{\bar{B}}\,(u_{\bar{B}}\in M_{\bar{B}})$. The linear map
$j_{\bar{A}}(v\otimes\varphi_D^{A,\Delta})$ on $\calv_D^A(M_{\bar{A}}\otimes M(0))$ sends $u_{\bar{A}}\otimes\unit_{A,0}$ to 
$j_{\bar{A}}(v\otimes\varphi_D^{A,\Delta})u_{\bar{A}}\otimes\unit_{A,0}
\equiv -u_{\bar{A}}\otimes J_0(v)\unit_{A,0}\mod\frakg_D^{A,\,out}(V)\left(M_{\bar{A}}\otimes M(0)\right)$. 
By Theorem \ref{th:decomptow} there exists the canonical isomorphism 
\[
\calv_D^A(M_{\bar{A}}\otimes M(0))\cong\bigoplus_{i=0}^{r}\calv_D^A(M_{\bar{A}}\otimes L_i)\otimes W_i^*\,,
\]
on which the linear map $j_{\bar{A}}(v\otimes\varphi^{A,\,\Delta})$ sends $u_{\bar{A}}\otimes\unit_{A,0}$ to
$-u_{\bar{A}}\otimes J_0(v)\unit_{A,0}=-u_{\bar{A}}\otimes\unit_{A,0}J_0(v)$ where $\unit_{A,0}J_0(v)$ indicates the right action of $A_0(V)$ to
$\unit_{A,0}$.
More precisely the linear map $j_{\bar{A}}(v\otimes\varphi^{A,\,\Delta})$ coincides with the  right $A_0(V)$ action on each component
$W_i^*\,(0\leq i\leq r)$; the set $\frakg_D^{out\,''}(V)\calv_D(M_{\bar{A}}\otimes M(0))$ is linearly spanned by
\[
-\bigoplus_{i=0}^{r}\calv_D^A(M_{\bar{A}}\otimes L_i)\otimes w_i^*J_0(v)\quad\mbox{for all}\quad w_i^*\in W_i^*,\,v\in
V_\Delta,\,\Delta\geq0.
\]
Applying the exactly same argument to $\calv_D^B(M_{\bar{B}})'$ we see that $\frakg_D^{out\,''}(V)\calv_D( M(0)\otimes
M_{\bar{B}})$ is linearly spanned by
\[
\bigoplus_{i=0}^{r}\calv_D^B(L_i\otimes M_{\bar{A}})\otimes w_i^*\theta\left(J_0(v)\right)\quad\mbox{for all}\quad w_i^*\in
W_i^*,\,v\in V_\Delta,\,\Delta\geq0.
\]
So $\calv_D(M_C)$ is isomorphic to 
\[
\bigoplus_{i,\,j=0}^{r}\calv_D^A(M_{\bar{A}}\otimes L_i)\otimes\calv_D^B(L_j\otimes M_{\bar{B}})
\otimes\left(W_i^*\otimes W_j^*/\mathcal{I}(W_i^*,W_j^*)\right)
\]
where $\mathcal{I}(W_i^*,W_j^*)$ is the vector space which is linearly spanned by 
$-w_i^*J_0(v)\otimes w_j^*+w_i^*\otimes w_j^*\cdot\theta\left(J_0(v)\right)$ for all
$v\in V_\Delta\,(\Delta\geq0)$ and $w_i^*\in W_i\,(0\leq i\leq r)$. 
We now see that
\[
W_i^*\otimes W_j^*/\mathcal{I}(W_i^*,W_j^*)\cong\Hom_{A_0(V)}(W_i,D(W_j))
=
\begin{cases}
\C&\quad (i=j^\dagger),\\
0&\quad (i\neq j^\dagger)
\end{cases}
\]
by Schur's lemma. 

\subsection{Proof of the glueing isomorphism}
We now prove Theorem \ref{th:factor}; by Lemma \ref{lem:vbar} and Theorem \ref{th:dstr} it suffices to construct the
$\calo_D[[q]]$-module isomorphism
$\Phi$ such that the following diagram commutes
\[
\begin{array}{ccc}
\Phi:\widehat{\calv}_{U}(M_C)&\longrightarrow&\displaystyle{\bigoplus_{i=0}^{r}}\calv_D^A(M_{\bar{A}}\otimes
L_i)\otimes_{\calo_D}\calv_D^B(L_{i^\dagger}\otimes M_{\bar{B}})\otimes\C[[q]]\\
\big\downarrow&&\hskip -7ex\big\downarrow\\
\bar{\Phi}:\overline{\calv}_D(M_C)&\longrightarrow&\hskip -9ex\displaystyle{\bigoplus_{i=0}^{r}}\calv_D^A(M_{\bar{A}}\otimes
L_i)\otimes_{\calo_D}\calv_D^B(L_{i^\dagger}\otimes M_{\bar{B}}),
\end{array}
\]
and $\Phi$ is compatible with the $\Theta(-\log D)$ action, where $\overline{\Phi}$ is $\calo_D$-module isomorphism given in Theorem
\ref{th:dstr}.

For any simple $V$-module $L$ there exists a complex number $\Delta$ such that
$L=\bigoplus_{d=0}^\infty L(d),\,L(d)=\{\,u\in L_i\;|\; T(0)u=(\Delta+d)u\}$ and $\dim\,L(d)<\infty$.
We denote  such a complex number for the simple $V$-modules $L_i$ by $\Delta_i$.
Let $L^\dagger$ be the contragredient dual of $L$.  
For each homogeneous space $L(d)$ and $L^\dagger(d)$ we fix the dual basis $\{u_{d,\,1},\,\ldots,\, u_{d,\,\ell_d}
\}\subset L(d)$ and $\{u_{d}^1,\,\ldots,\, u_{d}^{\ell_d}\}\subset L^\dagger(d)$ with respect to the canonical dual pairing.
We now define the element in $L\otimes L^\dagger\otimes\C[[q]]$ by
\[
\Omega(L)=\sum_{d=0}^\infty\left(
\sum_{j}u_{d,\,j}\otimes u_d^j\otimes q^{d}
\right).
\]
\begin{lemma}\label{lem:Omega}
For any $v\in V_\Delta$ and $n\in\Z$
\[
\left(J_n(v)\otimes 1 -1\otimes\theta\left(J_{n}(v)\right)q^{n}\right)\Omega(L)=0.
\]
\end{lemma}
\begin{proof}
It suffices to show that 
\[
\pairing{\,(J_n(v)\otimes 1)\Omega(L)}{u_{d_1,j_1}\otimes u_{d_2}^{j_2}\,}
=\pairing{\,(1\otimes\theta\left(J_{n}(v)\right)q^{n})\Omega(L)}{u_{d_1,j_1}\otimes u_{d_2}^{j_2}\,}.
\]
By definition
\begin{align*}
&\pairing{\,(J_n(v)\otimes 1)\Omega(L)}{u_{d_1,j_1}\otimes u_{d_2}^{j_2}\,}=\pairing{\,J_{n}(v)u_{d_1,j_1}}{u_{d_2}^{j_2}\,}q^{d_1}
\quad\mbox{and}\\
&\pairing{\,(1\otimes\theta\left(J_{n}(v)\right)q^{n})\Omega(L)}{u_{d_1,j_1}\otimes u_{d_2}^{j_2}\,}
=\pairing{\,u_{d_1,j_1}}{\theta\left(J_{n}(v)\right)u_{d_2}^{j_2}\,}q^{d_2+n},
\end{align*}
so that the right hand sides are zero if $d_1\neq d_2+n$ by definition of the dual pairing,
and those two coincide if $d_1=d_2+n$ by definition of the contragredient module.
\end{proof}

We set $\Omega_i=\Omega(L_i)$ for $0\leq i\leq r$ and define the linear map
\[
\begin{array}{ccc}
M_{\bar{A}}\otimes M_{\bar{B}} &\longrightarrow&M_{\bar{A}}\otimes
M_{\bar{B}}\otimes L_i\otimes
L_{i^\dagger}\otimes\C[[q]]\\
u_{\bar{A}}\otimes u_{\bar{B}} &\longmapsto&u_{\bar{A}}\otimes
u_{\bar{B}}\otimes\Omega_i\,.
\end{array}
\]
In order to see that this linear map induces an $\calo_D[[q]]$-module homomorphism 
\[
\widehat{\calv}_{U}(M_C)\rightarrow\calv_D^A(M_{\bar{A}}\otimes
L_i)\otimes\calv_D^B(L_{i^\dagger}\otimes M_{\bar{B}})\otimes\calo_D[[q]],
\]
it suffices to show that  
$j_{\bar{A},\bar{B}}(v\otimes\varphi)(u_{\bar{A}}\otimes u_{\bar{B}})\otimes\Omega_i
=j_{A\sqcup B}(v\otimes\varphi)(u_{\bar{A}}\otimes u_{\bar{B}}\otimes \Omega_i)$
for any $v\in V_\Delta$ and $\varphi\in\pi_*\omega^{1-\Delta}_{\calc_U/U}(*S_C)$, or equivalently
\begin{equation}\label{eqn:zero-sum}
\left(j_{A,0}(v\otimes\varphi)+j_{B,\infty}(v\otimes\varphi)\right)\Omega_i=0\quad\mbox{for all $i$.}
\end{equation}
Using the expansion (\ref{eqn:exp10}) and (\ref{eqn:exp11}) of $\varphi$ on $V_1^0$ and $V_1^\infty$, respectively,
we see that
\begin{align*}
j_{A,0}(v\otimes\varphi)&=\displaystyle{\sum_{n=0}^\infty\left(\sum_{j\geq-n}\varphi_{j+n+\Delta-1,n}^AJ_j(v)\right)}q^n
\quad\mbox{and}\\
j_{B,\infty}(v\otimes\varphi)&=(-1)^\Delta
\displaystyle{\sum_{n=0}^\infty\left(\sum_{j\geq-n}\varphi_{j+n,n+\Delta-1}^B\theta(J_j(v))q^j\right)}q^n.
\end{align*}
Then the  relation (\ref{eqn:zero-sum}) follows from  (\ref{eqn:ext}) and Lemma \ref{lem:Omega}.

Now the $\calo_D[[q]]$-module homomorphism $\Phi$ is induced by the map
$u_{\bar{A}}\otimes u_{\bar{B}}\longmapsto\sum_{i=0}^{r}u_{\bar{A}}\otimes u_{\bar{B}}\otimes\Omega_i$.
By the proof of Theorem \ref{th:dstr} it follows that $\left.\Phi\right|_{q=0}=\bar{\Phi}$ so that $\Phi$ is an isomorphism since
$\bar{\Phi}$ is an isomorphism.

Finally we show that the map $\Phi$ is compatible with the $\widehat{\Theta}_{U}(-\log D)$ action, i.e.,
\[
\sum_{b\in\bar{B}}w_{B,b}\,\rho_{B,b}(T(-1))u_{\bar{A}}\otimes u_{\bar{B}}\otimes\Omega
=\left(q\frac{\partial}{\partial q}+\Delta-\sum_{b\in\bar{B}}\rho_{B,b}(T(0))\right)u_{\bar{A}}\otimes u_{\bar{B}}\otimes\Omega.
\]
To this end it suffices to prove that for any $u_{\bar{B}}\in M_{\bar{B}}$ and $u_d\in L^\dagger(d\,)$
\begin{equation}\label{eqn:final}
\sum_{b\in\bar{B}}w_{B,b}\,\rho_{B,b}(T(-1))u_{\bar{B}}\otimes u_d
\equiv\left(\Delta+d-\sum_{b\in\bar{B}}\rho_{B,b}(T(0))\right)u_{\bar{B}}\otimes u_d
\end{equation}
on $\calv_D^B(L^\dagger\otimes M_{\bar{B}})$. Since $z_B(dz_B)^{-1}\in\Theta_{D\times\CP^1/D}$ it follows that
$j_B(T\otimes z_B(dz_B)^{-1})u_{\bar{B}}\otimes u_d\equiv0$ on $\calv_D^B(L^\dagger\otimes M_{\bar{B}})$, which shows
(\ref{eqn:final}) since
\[
j_{B,b}(T\otimes z_B(dz_B)^{-1})=
\begin{cases}
w_{B,b}T(-1)+T(0)&\quad(b\neq\infty),\\
-T(0)&\quad(b=\infty).
\end{cases}
\]

\subsection{Remark}
The definition of sheaves of conformal blocks over the partial compactification of the moduli space of
$N$-pointed projective lines is generalized to the one over the moduli space of $N$-pointed stable curves of genus zero,
and the local freeness of sheaves of conformal blocks can be proved. Then it follows that the category of $V$-modules is
rigid braided tensor category (\cite{BK}, \cite{KL1}, \cite{KL2}). Details  will be discussed in the subsequent papers.

\section{Semisimplicity of the category $\calm\!\operatorname{od}(V)$}\label{sec:category}
In this section we will show:
\begin{theorem}\label{th:semisimple02}
Let $V$ be a chiral vertex operator algebra which satisfies the condition \III\!. Then $\mathcal{M}od\,(V)$ is a semisimple category
and the functor $\mathcal{HW}:\mathcal{M}od\,(V)\rightarrow\mathcal{M}od\,(A_0(V))$ gives an equivalence of categories.
\end{theorem}

The rest of this section is devoted to the proof of Theorem \ref{th:semisimple02}. 
For this purpose we prepare lemma and propositions. 

Recall that $M(0)=\mathcal{U}(V)\otimes_{F^0\mathcal{U}(V)}A_0(V)$ is canonically a right
$A_0(V)$-module, which induces a left $A_0(V)$-action on $\Hom_{\,\calu(V)}(M(0),M)$.
\begin{lemma}\label{lem:hw01}
Let $M$ be a $V$-module.

\vskip 1.5ex
\noindent
{\rm (1)} Let $\varphi\in\Hom_{\,\calu(V)}(M(0),M)$. Then $\varphi(\unit)\in\mathcal{HW}(M)$ and the linear map
\[
\Hom_{\,\calu(V)}(M(0),M)\rightarrow\mathcal{HW}(M)\quad(\varphi\mapsto\varphi(\unit))
\]
is an $A_0(V)$-module isomorphism.

\vskip 1ex
\noindent
{\rm (2)}
The canonical pairing $\pairing{\;\,}{\;}:D(M)\times M\rightarrow\C$ induces a nondegenerate pairing 
$\pairing{\;\,}{\;}:\mathcal{CHW}(D(M))\times\mathcal{HW}(M)\rightarrow\C$.
\end{lemma}
\begin{proof}
The statement (1) follows from the definition of $M(0)$.
For the statement (2) we see that $\pairing{\,\frakg(V)^f_{<0}D(M)}{\mathcal{HW}(M)}=0$ as
$\theta(\frakg(V)^f_{<0})=\frakg(V)^f_{>0}$, so that the canonical pairing induces a pairing between $\mathcal{CHW}(D(M))$ and
$\mathcal{HW}(M)$. We now fix $m\in\mathcal{HW}(M)$ and suppose that
$\pairing{\mathcal{CHW}(D(M))}{m}=0$. Then we see that $\pairing{\,\varphi}{m}=0$ for all
$\varphi\in D(M)$ and that $m=0$ since the dual pairing between $D(M)$ and $M$ is nondegenerate.
\end{proof}

\begin{proposition}\label{prop:hw-chw}
Let $V$ be a chiral vertex operator algebra which satisfies the condition \III\!, and $M$ be a $V$-module.

\vskip 1.5ex
\noindent
{\rm (1)} If $M$ is semisimple then $M=\mathcal{U}(V)\mathcal{HW}(M)$.

\vskip 1ex
\noindent
{\rm (2)} There exists a unique maximal semisimple submodule $M_0$ of $M$. The submodule $M_0$ satisfies $\mathcal{HW}(M_0)\cong\mathcal{HW}(M)$.

\vskip 1ex
\noindent
{\rm (3)} There exists a unique maximal semisimple quotient module $\bar{M}$ of $M$. The quotient module $\bar{M}$ satisfies
$\mathcal{CHW}(\bar{M})\cong\mathcal{CHW}(M)$
\end{proposition}
\begin{proof}
(1) Note that $\calhw(M)\neq0$ by Lemma \ref{lem:hwe}. If $M$ is a simple $V$-module then $M=\calu(V)\calhw(M)$.
We now suppose that $M$ is semisimple and let $M=\bigoplus_{a}M^a$ be the irreducible decomposition. 
Then we see that $\calhw(M)=\bigoplus_{a}\calhw(M^a)$ and that
$\calu(V)\calhw(M)=\bigoplus_{a}\calu(V)\calhw(M^a)=M$.

(2) Since $A_0(V)$ is semisimple the $A_0(V)$-module $\mathcal{HW}(M)$ is completely reducible;
$\mathcal{HW}(M)=\bigoplus_{i}W_i$\,. 
We now set $M_0=\mathcal{U}(V)\mathcal{HW}(M)=\bigoplus_{i}\mathcal{U}(V)\otimes_{F^0\mathcal{U}(V)}W_i$\,, 
which is a semisimple $V$-submodule of $M$ by the condition I\!I\!I and $\mathcal{HW}(M_0)=\mathcal{HW}(M)$.
Let $N$ be a semisimple $V$-submodule of $M$. Then $\calhw(N)\subset\calhw(M)$ 
and $N=\calu(V)\calhw(N)\subset\calu(V)\calhw(M)=M_0$ by (1).

(3) Let $D(M)_0$ be the unique maximal semisimple $V$-submodule of $D(M)$ and set $\bar{M}=D(D(M)_0)$.
Then the $V$-module exact sequence
$0\rightarrow D(M)_0\rightarrow D(M)$ yields the $V$-modules exact sequence $0\leftarrow\bar{M}\leftarrow M$, 
which shows that $\bar{M}$ is a semisimple quotient module of $M$.
By (2) we see that $\calhw(D(M)_0)\cong\calhw(D(M))$, and that $\calch(\bar{M})\cong\calch(M)$ by Lemma \ref{lem:hw01} (2).
Now let $M\rightarrow\bar{N}\rightarrow0$ be a semisimple quotient of $M$. Then we see that
$\calch(M)\rightarrow\calch(\bar{N})\rightarrow0$, and by Lemma \ref{lem:hw01} (2) that
$0\rightarrow\calhw(D(\bar{N}))\rightarrow\calhw(D(M))\cong\calhw(D(M)_0)$. So $0\rightarrow
D(\bar{N})\rightarrow D(M)_0$ by (1), which shows $\bar{M}\rightarrow \bar{N}\rightarrow0$,  i.e., $\bar{M}$ is the  maximal
semisimple quotient.
\end{proof}

\begin{proposition}\label{prop:chw}
Let $V$ be a chiral vertex operator algebra which satisfies the condition \III\!, and $M$ be a $V$-module.

\vskip 1.5ex
\noindent
{\rm (1)} The linear map $\iota:M\rightarrow\calv_{(\infty,\,0)}(M(0)\otimes M)$
defined by $m\mapsto\unit_\infty\otimes m$ for all $m\in M$ induces an isomorphism of vector spaces
$\mathcal{CHW}(M)\cong\calv_{(\infty,\,0)}(M(0)\otimes M)$.

\vskip 1ex
\noindent
{\rm (2)} The canonical projection $\pi:M\rightarrow\bar{M}$ induces an isomorphism
$\calv_{(\infty,\,0)}(M(0)\otimes M)\cong\calv_{(\infty,\,0)}(M(0)\otimes\bar{M})$.

\vskip 1ex
\noindent
{\rm (3)}
Let $L$ be a simple $V$-module.
Then $\calv_{(\infty,\,0)}(L\otimes M)\cong\calv_{(\infty,\,0)}(L\otimes\bar{M})$ for any $V$-module $M$.

\vskip 1ex
\noindent
{\rm (4)}
Let $M^\infty,\,M^0$ be $V$-modules and $\bar{M}^0$
be the maximal semisimple quotient of
$M^0$. Then there exist vector space isomorphisms
\begin{align*}
&\calv_{(\infty,\,0)}(M^\infty\otimes M^0)\cong\calv_{(\infty,\,0)}(M^\infty\otimes\bar{M}^0),\\
&\Hom_{\,\calu(V)}(M^0,M^\infty)\cong\Hom_{\,\calu(V)}(\bar{M}^0,M^\infty).
\end{align*}
\end{proposition}
\begin{proof}
(1) By Proposition \ref{prop:fusion} we have the isomorphism 
\[
\begin{array}{ccc}
M/j_{0}\left(\bigoplus_{\Delta=0}^\infty V_\Delta\otimes H^0(\CP^1,\Omega^{1-\Delta}(*[0]-\Delta[\infty])
\right)M&\cong&\calv_{(\infty,\,0)}(M(0)\otimes M)\\
&&\\
m&\mapsto&\unit_\infty\otimes m.
\end{array}
\]
Since $\varphi_n=z^{-n+\Delta-1}(dz)^{1-\Delta}\,(n>0)$ form a topological basis of
$H^0(\CP^1,\Omega^{1-\Delta}(*[0]-\Delta[\infty]))$ and  $j_0(v\otimes\varphi_n)=J_{-n}(v)$ for all
$v\in V_\Delta$ we see that
\[
j_{0}\left(\bigoplus_{\Delta=0}^\infty V_\Delta\otimes H^0(\CP^1,\Omega^{1-\Delta}(*[0]-\Delta[\infty])
\right)M=\frakg(V)_{<0}^fM.
\]
The second statement follows from the first statement and Proposition \ref{prop:hw-chw} (3). 

(3) Let $L_i\,(0\leq i\leq r)$ be the complete list of simple $V$-modules. Since
$M(0)=\calu(V)\otimes_{F^0\calu(V)}A_0(V)\cong\bigoplus_{i}L_i\otimes W_i^*$ as a
$V$-module we see that
$\calv_{(\infty,0)}(M(0)\otimes M)\cong\bigoplus_{i}\calv_{(\infty,0)}(L_i\otimes M)\otimes W_i^*$ and 
$\calv_{(\infty,0)}(M(0)\otimes\bar{M})\cong\bigoplus_{i}\calv_{(\infty,0)}(L_i\otimes\bar{M})\otimes W_i^*$, so that
$\calv_{(\infty,0)}(L_i\otimes M)\cong\calv_{(\infty,0)}(L_i\otimes\bar{M})$ for all $i$.

(4) Let $\alpha,\beta$ be distinct points on $\CP^1$ such that $\alpha,\beta\neq0,\infty$.
Then by Theorem \ref{th:propvacua} there is an isomorphism
$\calv_{(\infty,\,0)}(M^\infty\otimes M^0)\cong\calv_{(\infty,\,\alpha,\,\beta,\,0)}(M^\infty\otimes V\otimes V\otimes M^0)$.
Using the factorization property (Theorem \ref{th:dstr}) and the statement (3) we see that
\[
\begin{split}
&\calv_{(\infty,\,0)}(M^\infty\otimes M^0)\\
&\cong\bigoplus_{i}\calv_{(\infty,\,\alpha,\,0)}(M^\infty\otimes V\otimes D(L_i))\otimes\calv_{(\infty,\,\beta,\,0)}(L_i\otimes V\otimes M^0)
\quad\mbox{(the factorization)}\\
&\cong\bigoplus_{i}\calv_{(\infty,\,\alpha,\,0)}(M^\infty\otimes V\otimes D(L_i))\otimes\calv_{(\infty,\,\beta,\,0)}(L_i\otimes V\otimes\bar{M}^0)
\quad\mbox{by (3)}\\
&\cong\calv_{(\infty,\,0)}(M^\infty\otimes\bar{M}^0)\quad\mbox{(the factorization)}.
\end{split}
\]
The second isomorphism follows from Proposition \ref{prop:homblock}
\end{proof}

To prove the theorem it suffices to show:
\begin{proposition}
Let $V$ be a chiral vertex operator algebra which satisfies the condition \III\! and $M$ be a $V$-module.
Then the maximal semisimple quotient module $\bar{M}$ of $M$ is isomorphic to $M$ as $V$-modules.
\end{proposition}
\begin{proof}
Let $\pi:M\rightarrow\bar{M}$ be the canonical projection. 
It suffices to show that there exists $f\in\Hom_{\,\calu(V)}(\bar{M},M)$ such that
$f\circ\pi=\id_M$.  Then the canonical projection $\pi$ is injective, i.e., $\pi$ is an isomorphism, 
and that $M$ is semisimple. The existence of such a homomorphism $f$ follows from the proof of Proposition \ref{prop:chw} (4).
\end{proof}

\appendix

\section{Appendix}\label{appendix}
\subsection{Examples}
Several chiral vertex operator algebras which satisfy all conditions required in this paper are known. 
By way of example we give affine chiral vertex operator algebras and minimal Virasoro chiral vertex operator algebras.

\begin{example}(Affine chiral vertex operator algebras)
Let $\frakg$ be a finite dimensional simple Lie algebra over $\C$ of rank $n$, and let
$\hat{\frakg}=\C[t,t^{-1}]\otimes\frakg\oplus\C K$ be the corresponding affine Lie algebra. We denote by
$\{\Lambda_{0},\ldots,\Lambda_{n}\}$ the set of fundamental weights for
$\hat{\frakg}$, and by $P^{k}_{+}$ the set of all level $k$ dominant integral weights. 
We denote  by $L(\Lambda)$ the irreducible highest weight module of $\hat{\frakg}$ with highest weight $\Lambda$.
It is known that if $k\neq -h^\vee,0$ then
$L_{\frakg,k}=L(k\Lambda_{0})$ has a structure of chiral vertex operator algebra where $h^\vee$ is the dual Coxeter number of
$\frakg$ (\cite[Theorem 2.9]{TK}). The energy-momentum tensor $T(z)$ of $L_{\frakg,k}$ is the Sugawara form. The $L_{\frakg,k}$ is
called an \textit{affine chiral vertex operator algebra}\,.

The case that the level $k$ is a positive integer has a particular interest.
The affine chiral vertex operator algebra $L_{\frakg,k}$ then satisfies
$C_2$-condition (\cite[Theorem 4.2.4]{TUY}).
Let $e_\theta$ be an element in the root space
$\frakg_\theta$ of the maximal root $\theta$. Then it is known (\cite{FZ}) that
\[
A_0(V)=U(\frakg)/\langle e_\theta^{k+1}\rangle\cong\bigoplus_{(\lambda,\theta)\leq k}L_\lambda\otimes L_\lambda^*
\]
where $\lambda\in\frakh^*$ is an integrable weight  and $L_\lambda$ is the simple $\frakg$-module with the highest weight $\lambda$.
In particular, $A_0(V)$ is semisimple.
Any $V$-module is an integrable $\hat{\frakg}$-module from the category $\calo$ (\cite{TK}).  Then the complete reducibility theorem
for  integrable $\hat{\frakg}$-modules (Theorem 10.7 in \cite{K}) shows that any $V$-module induced from a simple $A_0(V)$-module
is simple. So the affine chiral vertex operator algebra $L_{\frakg,k}$ satisfies the condition \III\!.
The simple $\hat{\frakg}$-module $L(\Lambda)\,(\Lambda\in P_+^k)$ is a simple $V$-module and
any simple $V$-module is isomorphic to $L(\Lambda)$ for a suitable $\Lambda\in P_+^k$.
\end{example}

\begin{example}(The Virasoro vertex operator algebras)
Let $\Vir=\bigoplus_{n\in\Z}\C T(n)\oplus\C C$ be the Virasoro algebra. 
We denote by $M(c,h)$ the Verma module for the Virasoro algebra with a highest weight $h\in\C$ and central charge $c\in\C$.
Let $v_{h,c}$ be the highest weight vector so that $L_{n}v_{h,c}=\delta_{n,0}hv_{h,c}\,(n\geq0)$ and
$Cv_{h,c}=cv_{h,c}$. The Verma module $M(c,h)$ is a free $U(\Vir^-)$-module of rank one with the generator $v_{h,c}$ where 
$\Vir^-=\bigoplus_{n\in\Z_{>0}}\C T(-n)$. Let $L(c,h)$ be the irreducible quotient of $M(c,h)$. 
Then it is known that $L(c,0)$ has a structure of chiral vertex operator algebra (\cite{BPZ}, \cite{BFM}).

In conformal field theory the Virasoro chiral vertex operator algebra with specific central charge is extensively studied.
Let $p,q$ be relatively prime integers such that $1<p<q$ and set
$c_{p,q}=1-6(p-q)^{2}/pq$.  For any integers $r$ and $s$ 
such that $0< r<p,\,0<s<q$ we denote $h_{p,q;r,s}=\{(rq-sp)^{2}-(p-q)^{2}\}/4pq$\;; 
note that $h_{p,q;r,s}=h_{p,q;p-r,q-s}$.
The module $L(c_{p,q},h_{p,q;r,s})$ for the Virasoro algebra is called a \textit{minimal series}, and in particular,
$L(c_{p,q},0)$ is called the \textit{minimal}\, Virasoro chiral vertex operator algebra. The minimal series $L(c_{p,q},0)$
satisfy $C_2$-condition (\cite{BFM}).

Let $V=L(c,0)$ be the Virasoro chiral vertex operator algebra with central charge $c\in\C$. 
For $c\neq c_{p,q}$ the zero-mode algebra is the polynomial ring $\C[x]$ , while for $c=c_{p,q}$ the zero-mode algebra
is a finite dimensional commutative algebra;
\[
A_0(V)=\C[x]/\langle G_{p,q}(x)\rangle,\quad G_{p,q}(x)^2=\prod_{\substack{0<r<p\\0<s<q}}(x-h_{p,q;r,s})
\]
where $G_{p,q}(x)$ is a polynomial such that $\deg G_{p,q}=(1/2)(p-1)(q-1)$ (\cite{FF}, \cite{BFM}).
For $c= c_{p,q}$ any $V$-module is completely reducible (\cite{BFM}) so that $V$ satisfies the condition \III\!.
The Virasoro module $L(c_{p,q},h_{p,q:r,s})$ is a simple $V$-module and any simple $V$-module is isomorphic to one of the minimal
series.
\end{example}

\subsection{Zhu's algebras via sheaves of covacua}\label{sub:zhu-zero}
It is stated in \cite{FZ} without proof that Zhu's algebra is isomorphic to the zero-mode algebra (Remark after Proposition 1.4.2
in \cite[page 13]{FZ}, also see \cite{BN}). Here we give a proof of the isomorphism between Zhu's algebra and the zero-mode algebra
using a space of covacua associated to
$3$ points $\infty,\, -1$ and $0$. 

We first recall definition of Zhu's algebra $A_z(V)$. Let $V$ be a chiral vertex operator algebra 
and $O(V)$ be the vector subspace of $V$, which is linearly spanned by vectors
\[
v_1\circ v_2=\Res{z}J(v_1,z)v_2\frac{(1+z)^{\Delta_1}}{z^2}\,dz=\sum_{n=0}^{\Delta_1}\binom{\Delta_1}{n}J_{-n-1}(v_1)v_2
\]
for all $v_1\in V_{\Delta_1},\,v_2\in V_{\Delta_2}$.

\begin{definition}
For any $v_1\in V_{\Delta_1}$ and $v_2\in V$ we define the binary operation $v_1*v_2$ by
\[
v_1*v_2=\Res{z}J(v_1,z)v_2\frac{(1+z)^{\Delta_1}}{z}\,dz=\sum_{n=0}^{\Delta_1}\binom{\Delta_1}{n}J_{-n}(v_1)v_2,
\]
and extend this to the bilinear operation $V\times V\rightarrow V$.
\end{definition}

The following fundamental results are found in \cite{Z1}.

\begin{proposition}\label{prop:zhu}
Let $V$ be a chiral vertex operator algebra.

\vskip 1.5ex
\noindent
{\rm (1)} $O(V)*V\subset O(V)$ and $V*O(V)\subset O(V)$.

\vskip 1ex
\noindent
{\rm (2)} $(T(-1)+T(0))v\in O(V)$ for all $v\in V$.

\vskip 1ex
\noindent
{\rm (3)} $(v_1*v_2)*v_3-v_1*(v_2*v_3)\in O(V)$ for all $v_1,\,v_2,\,v_3\in V$. 

\vskip 1ex
\noindent
{\rm (4)} For all $v_1\in V_{\Delta_1},\,v_2\in V_{\Delta_2}$ and $m\geq n\geq0$
\[
\Res{z}J(v_1,z)v_2\frac{(1+z)^{\Delta_1+n}}{z^{2+m}}\,dz\in O(V).
\]

\vskip 1ex
\noindent
{\rm (5)} For all $v_1\in V_{\Delta_1}$ and $v_2\in V_{\Delta_2}$
\[
v_1*v_2\equiv \Res{z}J(v_2,z)v_1\frac{(1+z)^{\Delta_2-1}}{z}\,dz\mod O(V).
\]
\end{proposition}

The proposition, in particular, shows:
\begin{proposition}
The quotient vector space $A_z(V)=V/O(V)$ is an associative algebra with a unit $\ket{0}+O(V)$
by the multiplication induced by the operation $*$.
\end{proposition}

\begin{definition}
The associative algebra $A_z(V)$ is called  \textit{Zhu's algebra}.
We use the notation $[v]=v+O(V)$, i.e., $[\;\;]:V\rightarrow A_z(V)$ denotes the canonical projection.
\end{definition}

We next construct a vector space isomorphism between Zhu's algebra and a space of covacua.

\begin{definition}
Let $D_z=[\infty]+[-1]$ be the divisor on $\CP^1$ and $O(D_z)$ be the vector subspace of $V$, which is linearly spanned by
$j_0(v\otimes\varphi)u$ for all $v\in V_\Delta,\,u\in V$ and $\varphi\in H^0(\CP^1, \Omega^{1-\Delta}(*[0]-\Delta D_z))$.
We define a quotient vector space by $A(D_z)=V/O(D_z)$.
\end{definition}

For any $\Delta\in\Z_{\geq0}$ and integer $m,\,n$ we set
\[
\varphi_{m,\,n}=\frac{(z+1)^{\Delta+n}}{z^{2+m}}(dz)^{1-\Delta}.
\]
The $\varphi_{m,\,n}\;(m\geq n\geq0)$ form a topological basis of $H^0(\CP^1,\Omega^{1-\Delta}(*[0]-\Delta D_z))$. 
On the one hand we see that $O(V)\subset O(D_z)$ by definition of $O(V)$.
On the other hand $O(D_z)\subset O(V)$ by Proposition \ref{prop:zhu} (4) so that $O(V)=O(D_z)$:

\begin{proposition}
$A_z(V)=A(D_z)$.
\end{proposition}

We now define a space of covacua which is isomorphic to the vector space $A_z(V)$. 
Let $A=\{\,\infty,\,-1,\,0\,\}$ and $w_A=(\infty,-1,0)$. Let $M(D_z)=M^\infty\otimes M^{-1}\otimes M^0$ where
$M^\infty=M(0),\,M^{-1}=M(0)$ and $M^0=V$.  We set $\calv_A(D_z)=M(D_z)/\frakg_{w_A}^{out}(V)M(D_z)$.
Recall that $M(0)=\mathcal{U}(V)/\mathcal{U}(V)F^{1}\mathcal{U}(V)$. We denote the unit $1+\mathcal{U}(V)F^{1}\mathcal{U}(V)$ 
by $\unit_\infty$ and $\unit_{-1}$ for $M^\infty$ and $M^{-1}$, respectively.

\begin{proposition}\label{prop:equiv}
The linear map
$V\rightarrow M(D_z)\,(u\mapsto \unit_{\infty}\otimes\unit_{-1}\otimes u)$ induces a linear isomorphism $A_z(V)\cong\calv_A(D_z)$.
\end{proposition}
\begin{proof}
We show that any element of $M(D_z)$ is congruent to $\unit_{\infty}\otimes\unit_{-1}\otimes u$ for a suitable $u\in
V$ modulo $\frakg_{w_A}^{out}(V)M(D_z)$ 
and $\unit_{\infty}\otimes\unit_{-1}\otimes u_1\equiv\unit_{\infty}\otimes\unit_{-1}\otimes u_2\mod\frakg_{w_A}^{out}(V)M(D_z)$ 
if and only if $u_1-u_2\in O(V)$.

Using the same argument in the proof of Theorem \ref{th:finiteblock} we see that the linear map $V\rightarrow M(D_z)\,(u\mapsto
\unit_{\infty}\otimes\unit_{-1}\otimes u)$ induces a surjective linear map $\iota:V\rightarrow\calv_A(D_z)$.  
For any $v\in V_\Delta$ and $\varphi\in H^0(\CP^1,\Omega^{1-\Delta}(*[0]-\Delta D_z))$ we see that
$j_{\infty}(v\otimes\varphi),\,j_{-1}(v\otimes\varphi)\in F^1\mathcal{U}(V)$, and that
\[
\unit_{\infty}\otimes\unit_{-1}\otimes j_0(v\otimes\varphi)u
=j_A(v\otimes\varphi)(\unit_{\infty}\otimes\unit_{-1}\otimes u)\in\frakg_{w_A}^{out}(V)M(D_z).
\]
Now the map $\iota$ induces a well-defined surjective linear map $A_z(V)\rightarrow\calv_A(D_z)$, which is also denote by $\iota$.

For the injectivity of $\iota$ we will prove that the dual map $\iota^*:\calv_A(D_z)^*\rightarrow A_z(V)^*$
is surjective. By using the same argument given in the proof of Proposition \ref{prop:fusion}
we can show that for any $\Phi\in A_z(V)^*$ there exists a system of current correlation functions $\{\Phi_m\}_{m=0}^\infty$
associated to $\Phi$ such that 
\begin{multline*}
\Phi_m(v_1,\ldots,v_m;u)_{(z_1,\ldots,z_m)}\\
\in  H^0((\CP^1)^m,\Omega^{\boxtimes_{i=1}^m\Delta_i}
(\sum_{i}*D_i+\sum_{i\neq j}*D_{i,\,j}+\sum_{i}\Delta_i D_{i,-1}+\sum_{i}\Delta_i D_{i,\infty})).
\end{multline*}
where $D_i=\{z_i=0\},\, D_{i,j}=\{z_i=z_j\},\, D_{i,-1}=\{z_i=-1\}$  and $D_{i,\infty}=\{z_i=\infty\}$ are divisors. 
Let us denote
\[
\Phi_m=\bra{\Phi}J(v_1,z_1)\cdots J(v_m,z_m)\ket{v}(dz_1)^{\Delta_1}\cdots(dz_m)^{\Delta_m}.
\]
We now define $\tilde{\Phi}\in\Hom_\C(M(D_z),\C)$ by
\begin{align*}
&\tilde{\Phi}(\theta\left(J_{m_1}(v_1)\right)\cdots\theta\left(J_{m_k}(v_k)\right)\unit_\infty\otimes J_{n_1}(u_1)\cdots
J_{n_\ell}(u_\ell)\unit_{-1}\otimes v)\\
&=\left(\frac{1}{2\pi\sqrt{-1}}\right)^{k+\ell}
\oint_{|z_1|=r_1}\!\!\!\!\!\!\cdots\;\oint_{|z_k|=r_k}\oint_{|w_1+1|=s_1}\!\!\!\!\!\!\cdots\;\oint_{|w_\ell+1|=s_\ell}\\
&\times\bra{\Phi}J(v_k,z_k)\cdots
J(v_1,z_1)J(u_\ell,w_\ell)\cdots J(u_1,w_1)\ket{v}\\ &\times z_k^{n_k+|v_k|-1}\cdots
z_1^{n_1+|v_1|-1}(w_\ell+1)^{m_\ell+|u_\ell|-1}\cdots(w_1+1)^{m_1+|u_1|-1}\,dw_\ell\cdots dw_1dz_k\cdots dz_1\,,
\end{align*}
where $r_k>r_{k-1}>\ldots>r_1>2$ and $1>s_1>s_{2}>\ldots>s_\ell>0$.
Then $\tilde{\Phi}$ induce an element of $\calv_A(D_z)^*$ such that $\iota^*(\tilde{\Phi})=\Phi$.
\end{proof}

\begin{proposition}\label{th:zhu}
Let $V$ be a chiral vertex operator algebra and $A_0(V)$ be its zero-mode algebra.
Let $o:V\rightarrow A_0(V)$ be the linear map defined by $o(v)=\left[J_0(v)\right]$.

\vskip 1.5ex
\noindent
{\rm (1)} The vector space $O(V)$ is in the kernel of the map $o$, and the induced map $o:A_z(V)\rightarrow A_0(V)$ is surjective.

\vskip 1ex
\noindent
{\rm (2)} The linear map $o:A_z(V)\rightarrow A_0(V)$ is an algebra homomorphism.
\end{proposition}
\begin{proof}
For any $v_1\in V_{\Delta_1}$ and $v_2\in V_{\Delta_2}$ we see that
\[
\begin{split}
&J_0(v_1\circ v_2)\\
&=\sum_{n=0}^{\Delta_1}\binom{\Delta_1}{n}\Res{z}J(J_{-n-1}(v_1)v_2,z)z^{\Delta_2+n}\,dz\\
&=\sum_{n=0}^{\Delta_1}\binom{\Delta_1}{n}\Res{z}\underset{w=z}{\operatorname{Res}}\;J(J(v_1,w-z)v_2,z)
(w-z)^{\Delta_1-n-2}z^{\Delta_2+n}\,dwdz\\
&=\Res{z}\underset{w=z}{\operatorname{Res}}\;J(J(v_1,w-z)v_2,z)
(w-z)^{-2}w^{\Delta_1}z^{\Delta_2}\,dwdz\\
&=\sum_{n=0}^\infty\binom{-2}{n}\Res{z}\Res{w}\;J(v_1,w)J(v_2,z)z^{\Delta_2+n}w^{\Delta_1-n-2}\,dwdz\\
&\phantom{\sum_{n=0}^\infty\Res{z}\Res{w}}-\sum_{n=0}^\infty\binom{-2}{n}\Res{w}\Res{z}\;
J(v_2,z)J(v_1,w)z^{\Delta_2-n-2}w^{\Delta_1+n}\,dzdw\\
&=\sum_{n=0}^\infty\binom{-2}{n}
\left(
J_{-n-1}(v_1)J_{n+1}(v_2)-J_{-n-1}(v_2)J_{n+1}(v_1)
\right)
\in \sum_{n\geq1}F^{-n}\mathcal{U}(V)\cdot F^{n}\mathcal{U}(V).
\end{split}
\]
Thus the map $o:A_z(V)\rightarrow A_0(V)$ is well-defined and is surjective by Proposition \ref{prop:zero-trivial}.

We note that
\[
\begin{split}
J_0(v_1*v_2)&=\sum_{n=0}^{\Delta_1}\binom{\Delta_1}{n}\Res{z}J(J_{-n}(v_1)v_2,z)z^{\Delta_2+n-1}\,dz\\
&=\sum_{n=0}^{\Delta_1}\binom{\Delta_1}{n}\Res{z}\underset{w=z}{\operatorname{Res}}\;J(J(v_1,w-z)v_2,z)
(w-z)^{\Delta_1-n-1}z^{\Delta_2+n-1}\,dwdz\\
&=\Res{z}\underset{w=z}{\operatorname{Res}}\;J(J(v_1,w-z)v_2,z)
(w-z)^{-1}z^{\Delta_2-1}w^{\Delta_1}\,dwdz\;(=: I).
\end{split}
\]
On the one hand $J(J(v_1,w-z)v_2,z)=J(v_1,w)J(v_2,z)\,(|w|>|z|>0)$, and on the other hand $J(J(v_1,w-z)v_2,z)=J(v_2,w)J(v_1,z)\,(|z|>|w|>0)$; 
we then see that
\[
\begin{split}
I&=\sum_{n=0}^\infty\Res{z}\Res{w}\;J(v_1,w)J(v_2,z)z^{\Delta_2+n-1}w^{\Delta_1-n-1}\,dwdz\\
&\phantom{\sum_{n=0}^\infty\Res{z}\Res{w}\;J(v_1,w)}-\sum_{n=0}^\infty\Res{w}\Res{z}\;J(v_2,z)J(v_1,w)z^{\Delta_2-n-2}w^{\Delta_1+n}\,dzdw\\
&=\sum_{n=0}^\infty J_{-n}(v_1)J_n(v_2)-\sum_{n=0}^\infty J_{-n-1}(v_2)J_{n+1}(v_1),
\end{split}
\]
and that $J_0(v_1*v_2)\equiv J_0(v_1)J_0(v_2)\mod \sum_{n\geq1}F^{-n}\mathcal{U}(V)\cdot F^{n}\mathcal{U}(V)$.
\end{proof}

Recall that there exists a $\calu(V)$-module isomorphism $M(0)\cong\mathcal{U}(V)\otimes_{F^0\mathcal{U}(V)}A_0(V)$, 
which defines the canonical right $A_0(V)$-module structure on $M(0)$. 

\begin{definition}
We define the linear map $\sigma_{-1}:A_0(V)\rightarrow\calv_A(D_Z)\cong A_z(V)$ and $\sigma_\infty:A_0(V)\rightarrow\calv_A(D_Z)\cong A_z(V)$ by
$a\mapsto \unit_\infty\otimes a\otimes\ket{0}$ and $a\mapsto a\otimes\unit_{-1}\otimes\ket{0}$ for all $a\in A_0(V)$, respectively.
\end{definition}

\begin{proposition}
The map $\sigma_{-1},\,\sigma_\infty:A_0(V)\rightarrow A_z(V)$ is respectively an algebra and an anti-algebra homomorphism
such that $\sigma_{-1}([J_0(v)])=[v]$ and $\sigma_{\infty}([\theta(J_0(v))])=[v]$ for all $v\in V_\Delta$.
\end{proposition}
\begin{proof}
We first prove the second assertion.
For $\varphi=z^{-1}(1+z)^{\Delta-1}(dz)^{1-\Delta}$ and any $v\in V_\Delta$  
we see that $j_{-1}(v\otimes\varphi)\unit_{-1}=-J_0(v)\unit_{-1}$, and that
\[
\begin{split}
\unit_{\infty}\otimes \left[J_0(v)\right]\unit_{-1}\otimes\ket{0}
&=-\unit_{\infty}\otimes j_{-1}(v\otimes\varphi)\unit_{-1}\otimes\ket{0}\\
&\equiv\unit_{\infty}\otimes\unit_{-1}\otimes\Res{z}J(v,z)\ket{0}\frac{(1+z)^{\Delta-1}}{z}\,dz\!\mod\frakg_{w_A}^{out}(V)M(D_z)\\ 
&\equiv\unit_{\infty}\otimes\unit_{-1}\otimes v.
\end{split}
\]
On the other hand, for any $v\in V_\Delta$ take $\varphi=z^{-1}(1+z)^{\Delta}\,(dz)^{1-\Delta}$.  
Then we see that
$j_\infty\left(v\otimes\varphi\right)\unit_\infty=-\theta\left(J_0(v)\right)\unit_\infty$, and that
\[
\begin{split}
\theta\left(J_0(v)\right)\unit_\infty\otimes\unit_{-1}\otimes\ket{0}
&=-j_\infty(v\otimes\varphi)\unit_{\infty}\otimes\unit_{-1}\otimes\ket{0}\\
&\equiv \unit_{\infty}\otimes\unit_{-1}\otimes
\Res{z}J(v,z)\ket{0}\frac{(1+z)^\Delta}{z}\,dz\!\mod\frakg_{w_A}^{out}(V)M(D_z)\\
&=\unit_{\infty}\otimes\unit_{-1}\otimes v.
\end{split}
\]

Similarly we find that 
\begin{align*}
&\unit_{\infty}\otimes \left[J_0(v_1)\right]\left[J_0(v_2)\right]\unit_{-1}\otimes\ket{0}\equiv\unit_{\infty}\otimes\unit_{-1}\otimes
v_1*v_2
\quad\mbox{and}\\
&\left[\theta(J_0(v_1))\right]\left[\theta(J_0(v_2))\right]\unit_{\infty}\otimes\unit_{-1}\otimes\ket{0}\equiv\unit_{\infty}\otimes\unit_{-1}\otimes
v_2*v_1,
\end{align*}
which show the first assertion, where we also use the fact $\theta^2(J_0(v))=J_0(v)$ (see \S\,\ref{sub:theta}).
\end{proof}

\begin{theorem}
Let $V$ be a chiral vertex operator algebra and $A_0(V)$ be its zero-mode algebra.
Then $o:A_z(V)\rightarrow A_0(V)$ is an algebra isomorphism. 
\end{theorem}

\begin{proposition}\label{prop:theta-involution}
The map $\theta:A_0(V)\rightarrow A_0(V)\, (J_0(v)\mapsto \theta(J_0(v)))$ is well-defined, and is anti-algebra involution.
\end{proposition}
\begin{proof}
By the isomorphism $A_0(V)\cong\calv_A(D_Z)$ the map $\theta$ is well-defined. Since $\sigma_{-1}=\sigma_{\infty}\circ\theta$ we see that $\theta$ is
an anti-algebra homomorphism.
\end{proof}

\vskip 5ex
\noindent
\begin{flushleft}
Kiyokazu Nagatomo\footnote{Supported in part by Grant-In-Aid for Scientific Research.}\hskip 5em\\
Department of Pure and Applied Mathematics\\
Graduate School of Information Science and Technology\\
Osaka University,
Osaka, Toyonaka 560-0043, Japan\\
e-mail: nagatomo$@$math.sci.osaka-u.ac.jp
\vskip 2ex
\noindent
Akihiro Tsuchiya\\
Graduate School of Mathematics, Nagoya University\\
Nagoya 464-8602, Japan\\
e-mail:tsuchiya$@$math.nagoya-u.ac.jp
\end{flushleft}
\end{document}